\documentclass[12pt]{amsart}
\usepackage{amssymb}

\input{prepictex}
\input{pictex}
\input{postpictex}

\def\CC{\mathbb C}
\def\NN{\mathbb N}
\def\QQ{\mathbb Q}
\def\TT{\mathbb T}
\def\ZZ{\mathbb Z}

\def\Ff{\mathcal F}
\def\Kk{\mathcal K}
\def\Mm{\mathcal M}

\def\Tt{\mathcal T}

\def\Ad{\operatorname{Ad}}

\def\Ext{\operatorname{Ext}}

\def\id{\operatorname{id}}

\def\Aut{\operatorname{Aut}}

\def\lsp{\operatorname{span}}
\def\clsp{\operatorname{\overline{span\!}\,\,}}
\def\lcm{\operatorname{lcm}}

\def\RR{\operatorname{RR}}

\theoremstyle{plain}
\newtheorem{theorem}{Theorem}[section]
\newtheorem*{theorem*}{Theorem}
\newtheorem*{prop*}{Proposition}
\newtheorem{cor}[theorem]{Corollary}
\newtheorem{lemma}[theorem]{Lemma}
\newtheorem{prop}[theorem]{Proposition}
\theoremstyle{remark}
\newtheorem{rmk}[theorem]{Remark}

\newtheorem{example}[theorem]{Example}

\theoremstyle{definition}
\newtheorem{dfn}[theorem]{Definition}

\newtheorem{notation}[theorem]{Notation}

\def\gbd/{rank-2 Bratteli diagram}
\def\Gbd/{Rank-2 Bratteli diagram}
\def\gbds/{rank-2 Bratteli diagrams}
\def\Gbds/{Rank-2 Bratteli diagrams}

\newcommand\BLUE{f^*_1}
\newcommand\RED{f^*_2}

\def\arr(#1,#2)(#3,#4){\arrow <0.375em> [0.25, 0.75] from {#1} {#2} to {#3} {#4}}

\def\overarrow[#1]{\buildrel #1 \over \longrightarrow}

\def\comment#1{{}}

\numberwithin{equation}{section}

\title[Rank-2 graphs and A$\TT$ algebras]{\boldmath{Rank-two graphs whose $C^*$-algebras are \\ direct limits of circle algebras}}
\author{David Pask}
\address{David Pask, Iain Raeburn and Aidan Sims\\ School of Mathematical and Physical Sciences  \\
University of Newcastle\\ 
NSW  2308\\
AUSTRALIA}
\email{david.pask, iain.raeburn, aidan.sims@newcastle.edu.au}
\author{Iain Raeburn}
\author{Mikael R\o rdam}
\address{Mikael R\o rdam\\ Department of Mathematics and Computer Science\\
University of Southern Denmark, Odense\\ Campusvej 55\\ DK-5230 Odense M\\
DENMARK}
\email{mikael@imada.sdu.dk}
\author{Aidan Sims}

\date{December 13, 2005}
\thanks{This research was supported by the Australian Research Council through the ARC Centre for Complex Dynamic Systems and Control. The fourth author was supported by CDSC and by an Australian Postdoctoral Fellowship.}

\begin{document}

\begin{abstract}
We describe a class of rank-2 graphs whose $C^*$-algebras are A$\TT$ algebras. For a subclass which we call \gbds/, we compute the $K$-theory of the $C^*$-algebra. We identify \gbds/ whose $C^*$-algebras are simple and have real-rank zero, and characterise the $K$-invariants achieved by such algebras. We give examples of \gbds/ whose $C^*$-algebras contain as full corners the irrational rotation algebras and the Bunce-Deddens algebras.
\end{abstract}

\maketitle

\section{Introduction}
The $C^*$-algebras of directed graphs are generalisations of the Cuntz-Krieger algebras of finite $\{0,1\}$-matrices. Graph algebras have an attractive structure theory in which algebraic properties of the $C^*$-algebra are determined by easily visualised properties of the underlying graph. (This theory is summarised in \cite{CBMSbk}.) In particular, it is now easy to decide whether a given graph algebra is simple and purely infinite, and a theorem of Szyma\'nski \cite{Sz} says that every Kirchberg algebra with torsion-free $K_1$ is isomorphic to a corner in a graph algebra. Each AF algebra $A$ can also be realised as a corner in the graph algebra of a Bratteli diagram for $A$ \cite{D, JT}. But this is all we can do: the dichotomy of \cite{KPR} says that every simple graph algebra is either purely infinite or AF, and a theorem of \cite{RSz} says that every graph algebra has torsion-free $K_1$.

Higher-rank analogues of Cuntz-Krieger algebras and graph algebras have been introduced and studied by Robertson-Steger \cite{RobSt2} and Kumjian-Pask \cite{KP}, and are currently attracting a good deal of attention (see  \cite{FMY, Hop, KriPow, RSY2, SZach}, for example). The theory of higher-rank graph algebras mirrors in many respects the theory of ordinary graph algebras, and we have good criteria for deciding when the $C^*$-algebra of a higher-rank graph is simple or purely infinite \cite{KP}. Since these algebras include tensor products of graph algebras, the $K_1$-group of such an algebra can have torsion, so these algebras include more models of Kirchberg algebras than ordinary graph algebras. However, it is not obvious which finite $C^*$-algebras can be realised as the $C^*$-algebras of higher-rank graphs. Indeed, we are not aware of any results in this direction.

Here we discuss a class of rank-$2$ graphs whose $C^*$-algebras are A$\TT$ algebras. We specify a 2-graph $\Lambda$ using a pair of coloured graphs which we call the \emph{blue graph} and the \emph{red graph} with a common vertex set together with a \emph{factorisation property} which identifies each red-blue path of length 2 with a blue-red path. In the 2-graphs which we construct, the blue graph is a Bratteli diagram and the red graph partitions the vertices in each level of the diagram into a collection of disjoint cycles. The $C^*$-algebra of such a \gbd/ $\Lambda$ then has a natural inductive structure. We prove that $C^*(\Lambda)$ is always an A$\TT$ algebra with nontrivial $K_1$-group, and in particular is neither purely infinite nor AF. We compute the $K$-theory of $C^*(\Lambda)$, and produce conditions which ensure that $C^*(\Lambda)$ is simple with real-rank zero. Using these results and Elliott's classification theorem, we identify \gbds/ whose $C^*$-algebras contain as full corners the Bunce-Deddens algebras and the irrational rotation algebras. Under the additional hypothesis that all red cycles in $\Lambda$ have length 1, we improve our analysis of the real rank of $C^*(\Lambda)$, and describe the trace simplex.

\vskip1em
The paper is organised as follows. In Section~\ref{sec:prelim}, we briefly recap the standard definitions and notation for $k$-graphs and their $C^*$-algebras. In Section~\ref{sec:AT 2-graphs}, we describe a class of $2$-graphs $\Lambda$ whose $C^*$-algebras are A$\TT$ algebras. The blue graph of such a $2$-graph $\Lambda$ is a graph with no cycles. The red graph consists of a union of disjoint isolated cycles. Very roughly speaking, the red cycles in $\Lambda$ give rise to unitaries in $C^*(\Lambda)$ while finite collections of blue paths index matrix units in $C^*(\Lambda)$. So carefully constructed finite subgraphs of $\Lambda$ correspond to subalgebras of $C^*(\Lambda)$ which are isomorphic to direct sums of matrix algebras over $C(\TT)$. We write $C^*(\Lambda)$ as the increasing union of these circle algebras to show that $C^*(\Lambda)$ is A$\TT$ (Theorem~\ref{thm:AT alg}).

In Section~\ref{sec:GBDs} we assume further that the blue graph is a Bratteli diagram and that the red graph respects the inductive structure of the diagram, and call the resulting $2$-graphs \emph{\gbds/}. In Theorem~\ref{thm:AT and K-th} we compute the $K$-theory of $C^*(\Lambda)$ for a \gbd/ $\Lambda$. Our arguments are elementary and do not depend on the computations of $K$-theory for general $2$-graph algebras \cite{RobSt3,Ev}. The $K$-theory calculation shows in particular that $K_1(C^*(\Lambda))$ is isomorphic to a subgroup $G$ of $K_0(C^*(\Lambda))$ such that $K_0(C^*(\Lambda))/G$ has rank zero. We also establish a bijection between the gauge-invariant ideals of $C^*(\Lambda)$ and the order ideals of the dimension group $K_0(C^*(\Lambda))$.

Elliott's classification theorem for A$\TT$ algebras says that each A$\TT$ algebra with real-rank zero is determined up to stable isomorphism by its ordered $K_0$-group and its $K_1$-group \cite{Ell}. Thus we turn our attention in Section~\ref{sec:LPF} to identifying \gbds/ whose $C^*$-algebras have real-rank zero. Theorem~\ref{critsimple} establishes a necessary and sufficient condition on $\Lambda$ for $C^*(\Lambda)$ to be simple. We then identify a \emph{large-permutation factorisations} property which guarantees that projections in $C^*(\Lambda)$ separate traces, and deduce from \cite{BBEK} that if $\Lambda$ has large-permutation factorisations and is cofinal in the sense of \cite{KP}, then $C^*(\Lambda)$ is simple and has real-rank zero (Theorem~\ref{thm:LPF -> RR0}). 

In Section~\ref{sec:Achievability} we identify the pairs $(K_0,K_1)$ which can arise as the $K$-theory of the $C^*$-algebra of a \gbd/, and identify among these the pairs which are achievable when $C^*(\Lambda)$ is simple and has real-rank zero. We then construct for each irrational number $\theta \in (0,1)$ a \gbd/ $\Lambda_\theta$ with large-permutation factorisations such that the irrational rotation algebra $A_\theta$ is isomorphic to a full corner of $C^*(\Lambda_\theta)$, and for each infinite supernatural number $\mathbf{m}$ a \gbd/ $\Lambda(\mathbf{m})$ with large-permutation factorisations such that the Bunce-Deddens algebra of type $\mathbf{m}$ is isomorphic to a full corner of $C^*(\Lambda(\mathbf{m}))$. These algebras are examples of simple $2$-graph $C^*$-algebras which are neither AF nor purely infinite, and show that the dichotomy of \cite{KPR} fails for $2$-graphs.

In the last two sections we consider \gbds/ in which all the red cycles have length~$1$. For such \gbds/ we show that the partial inclusions between approximating circle algebras are standard permutation mappings, and identify the associated permutations explicitly in terms of the factorisation property in $\Lambda$. We then consider arbitrary direct limits of circle algebras under such inclusions. We give a sufficient condition for such algebras to have real-rank zero and a related necessary condition. We describe the trace simplex in both cases. When each approximating algebra contains just one direct summand, we obtain a single necessary and sufficient condition for the limit algebra to have real-rank zero.

\section{Preliminaries and notation} \label{sec:prelim}
Our conventions regarding $2$-graphs are largely those of \cite{KP}. By $\NN^2$, we mean the semigroup $\{(n_1,n_2)\in \ZZ^2:n_i\geq 0\}$; we write $e_1=(1,0)$, $e_2=(0,1)$ and $0$ for the identity $(0,0)$. We view $\NN^2$ as a category with one object. We define a partial order on $\NN^2$ by $m \le n$ if and only if $m_1 \le n_1$ and $m_2\leq n_2$. 

A \emph{$2$-graph} is a pair $(\Lambda,d)$ consisting of a countable category $\Lambda$ and a functor $d : \Lambda \to \NN^2$ which satisfies the factorisation property: if $\lambda \in \Lambda$ and $d(\lambda) = m+n$ then there exist unique $\mu$ and $\nu$ in $\Lambda$ such that $d(\mu) = m$, $d(\nu) = n$ and $\lambda = \mu\nu$.
We refer to the morphisms of $\Lambda$ as the \emph{paths} in $\Lambda$; the degree $d(\lambda)$ is the rank-2 analogue of the length of the path $\lambda$, and we write $\Lambda^n:=d^{-1}(n)$. We call the objects of $\Lambda$ vertices, the domain $s(\lambda)$ of $\lambda$ the source of $\lambda$, and  the codomain $r(\lambda)$ the range of $\lambda$. We refer to the paths in $\Lambda^{e_1}$ as blue edges and those in $\Lambda^{e_2}$ as red edges. 

The factorisation property applies with $d(\lambda)=0+d(\lambda)=d(\lambda)+0$, and the uniqueness of factorisations then implies that $v\mapsto \id_v$ is a bijection between the objects of $\Lambda$ and the set $\Lambda^0$ of paths of degree $0$. We use this bijection to to identify the objects of $\Lambda$ with the paths $\Lambda^0$ of degree $0$, and we view $r,s$ as maps from $\Lambda$ to $\Lambda^0$. We say that $\Lambda$ is row-finite if the set $\{\lambda \in \Lambda : r(\lambda) = v, d(\lambda) = n\}$ is finite for every vertex $v \in \Lambda^0$ and every degree $n \in \NN^2$. A $2$-graph $\Lambda$ is locally convex if, whenever $e$ is a blue edge and $f$ is a red edge with $r(e) = r(f)$, there exist a red edge $f'$ and a blue edge $e'$ such that $r(f') = s(e)$ and $r(e') = s(f)$.

For $\lambda \in \Lambda$ and $0 \le m \le n \le d(\lambda)$, we write $\lambda(m,n)$ for the unique path in $\Lambda$ such that $\lambda = \lambda'\lambda(m,n)\lambda''$ where $d(\lambda') = m$, $d(\lambda(m,n)) = n-m$ and $d(\lambda'') = d(\lambda) - n$. We write $\lambda(n)$ for $\lambda(n,n) = s(\lambda(0,n))$. 

For  $E \subset \Lambda$ and $\lambda \in \Lambda$, we denote by $\lambda E$ the collection $\{\lambda\mu : \mu \in E, r(\mu) = s(\lambda)\}$ of paths which extend $\lambda$. Similarly $E\lambda := \{\mu\lambda : \mu \in E, s(\mu) = r(\lambda)\}$. In particular when $\lambda = v \in \Lambda^0$, $E v = E \cap s^{-1}(v)$ and $vE = E \cap r^{-1}(v)$.

As in \cite{RSY1}, we write 
\[
\Lambda^{\le n} := \{\lambda \in \Lambda : d(\lambda) \le n, \lambda\mu \in \lambda\Lambda\setminus\{\lambda\} \implies d(\lambda\mu) \not\le n\}.
\]

The $C^*$-algebra $C^*(\Lambda)$ of a row-finite locally convex $2$-graph $\Lambda$ is the universal algebra generated by a collection $\{s_\lambda : \lambda \in \Lambda\}$ of partial isometries (called a Cuntz-Krieger $\Lambda$-family) satisfying
\begin{itemize}
\item[(CK1)] $\{s_v : v \in \Lambda^0\}$ is a collection of mutually orthogonal projections;\smallskip
\item[(CK2)] $s_\mu s_\nu = s_{\mu\nu}$ whenever $s(\mu) = r(\nu)$;\smallskip
\item[(CK3)] $s^*_\mu s_\mu = s_{s(\mu)}$ for all $\mu$; and\smallskip
\item[(CK4)] $s_v = \sum_{\lambda \in v\Lambda^{\le n}} s_\lambda s^*_\lambda$ for all $v \in \Lambda^0$ and $n \in \NN^2$.
\end{itemize}
Proposition~3.11 of \cite{RSY1} shows that~(CK4) is equivalent to 
\[\textstyle
s_v = \sum_{e \in v\Lambda^{e_i}} s_e s^*_e\text{ for all $v \in \Lambda^0$ and $i = 1,2$.}
\]
There is a strongly continuous action $\gamma : \TT^2 \to \Aut(C^*(\Lambda))$ called the gauge action which satisfies $\gamma_z(s_\lambda) = z^{d(\lambda)} s_\lambda$ for all $\lambda$, where $z^n := z_1^{n_1} z_2^{n_2} \in \TT$ for $z \in \TT^2$ and $n \in \NN^2$.

For $i = 1, 2$, the function $f_i$ is the homomorphism $f_i(n) = n e_i$ of $\NN$ into $\NN^2$. As in \cite{KP}, we write $f^*_i \Lambda$ for the pullback
\[
f^*_i\Lambda := \{(\lambda,n) : n \in \NN, \lambda \in \Lambda, f_i(n) = d(\lambda)\},
\]
which is a $1$-graph with the same vertex set as $\Lambda$. We call $\BLUE\Lambda$ the blue graph and $\RED\Lambda$ the red graph. In this paper, we identify $f^*_i\Lambda$ with $d^{-1}(f_i(\NN)) \subset \Lambda$.

A \emph{cycle} in a $k$-graph is a path $\lambda$ such that $d(\lambda) \not= 0$, $r(\lambda) = s(\lambda)$ and $\lambda(n) \not= s(\lambda)$ for $0 < n < d(\lambda)$. A \emph{loop} is a cycle consisting of a single edge. We say that the cycle $\lambda$ has an \emph{entrance} if there exists $n \le d(\lambda)$ such that $r(\lambda) \Lambda^n \setminus \{\lambda(0,n)\}$ is nonempty. Likewise, $\lambda$ is said to have an \emph{exit} if there exists $n \le d(\lambda)$ such that $\Lambda^n s(\lambda) \setminus \{\lambda(d(\lambda) - n, d(\lambda))\}$ is nonempty. We say that the cycle $\lambda$ is \emph{isolated} if it has no entrances and no exits. 

We say that a $k$-graph $\Lambda$ is \emph{finite} if $\Lambda^n$ is finite for each $n$. If $\Lambda$ is row-finite, this is equivalent to the assumption that $\Lambda^0$ is finite.

\section{Rank-2 graphs with A$\TT$ $C^*$-algebras} \label{sec:AT 2-graphs}

In this section we prove the following theorem which identifies a class of $2$-graphs whose $C^*$-algebras are A$\TT$ algebras. 

\begin{theorem} \label{thm:AT alg}
Let $(\Lambda,d)$ be a $2$-graph such that
\begin{equation}\label{ass:AT}
\ \hfill
\parbox{0.9\textwidth}{
$\Lambda$ is row-finite, $\BLUE\Lambda$ contains no cycles, each vertex $v \in \Lambda^0$
is the range of an isolated cycle $\lambda(v)$ in $\RED\Lambda$.
}%\end{split}%\tag{AT}
\end{equation}
Then $C^*(\Lambda)$ is an A$\TT$ algebra.
\end{theorem}

%The $2$-graphs we explore in this paper satisfy the following assumption.

%\begin{ass} \label{ass:AT}
%$\Lambda$ is row-finite, $\BLUE\Lambda$ contains no cycles, and $\RED\Lambda$ is the path-space of a collection of isolated cycles.
%\end{ass}

If $\Lambda$ satisfies condition~\eqref{ass:AT} then each vertex $v$ lies on a unique isolated cycle $\lambda(v)$ and hence $v$ is the range of exactly one red edge and the source of exactly one red edge. If $e$ is a blue edge there are unique red edges $f,g$ with $r(f) = r(e)$ and $r(g) = s(e)$ and then the factorisation property implies that there is a unique blue edge $e'$ such that $eg = fe'$. Hence every $2$-graph $\Lambda$ which satisfies condition~\eqref{ass:AT} is locally convex.

\begin{figure}[ht]
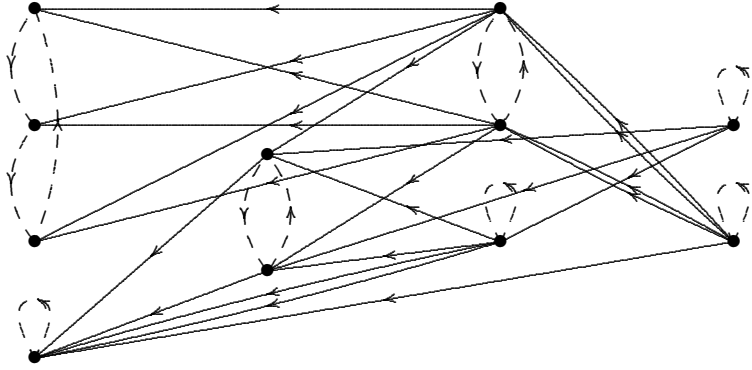

\(
\beginpicture
\setcoordinatesystem units <3.75em,3.75em>
\put{$\bullet$} at 0 0
\put{$\bullet$} at 0 1
\put{$\bullet$} at 0 2
\put{$\bullet$} at 0 3
\put{$\bullet$} at 2 .75
\put{$\bullet$} at 2 1.75
\put{$\bullet$} at 4 1
\put{$\bullet$} at 4 2
\put{$\bullet$} at 4 3
\put{$\bullet$} at 6 1
\put{$\bullet$} at 6 2
\setdashes <0.417em>
\setquadratic
\plot 0 0  -.135 .375 0 .5 /    \plot 0 0  .135 .375 0 .5 /
\arr(.031,.5045)(.03,.505)
\plot 0 1 -.2 1.5 0 2 /     \plot 0 2 -.2 2.5 0 3 /     \plot 0 3 .2 2 0 1 /
\arr(-.2,1.51)(-.2,1.5)    \arr(-.2,2.51)(-.2,2.5)    \arr(.2,2.04)(.2,2.05)
\plot 2 .75 1.8 1.25 2 1.75 /     \plot 2 .75 2.2 1.25 2 1.75 /
\arr(1.8,1.21)(1.8,1.2)           \arr(2.2,1.29)(2.2,1.3)
\plot 4 1  3.865 1.375 4 1.5 /    \plot 4 1  4.135 1.375 4 1.5 /
\arr(4.031,1.5045)(4.03,1.505)
\plot 4 2 3.8 2.5 4 3 /     \plot 4 2 4.2 2.5 4 3 /
\arr(3.8,2.46)(3.8,2.45)           \arr(4.2,2.54)(4.2,2.55)
\plot 6 1  5.865 1.375 6 1.5 /    \plot 6 1  6.135 1.375 6 1.5 /
\arr(6.031,1.5045)(6.03,1.505)
%\plot 6 1  5.9 1.4 6 1.5 /  \plot 6 1  6.1 1.4 6 1.5 /
%\arr(5.951,1.5)(5.95,1.5)
\plot 6 2  5.865 2.375 6 2.5 /    \plot 6 2  6.135 2.375 6 2.5 /
\arr(6.031,2.5045)(6.03,2.505)
%\plot 6 2  5.9 2.4 6 2.5 /  \plot 6 2  6.1 2.4 6 2.5 /
%\arr(5.951,2.5)(5.95,2.5)
%
\setsolid
\setlinear
\plot 0 0 2 .75 /     \plot 0 0 2 1.75 /
\arr(1.008,.378)(1,.375) \arr(1.008,.882)(1,.875)
\setquadratic \plot 0 0 2 .45 4 1 /    \plot 0 0 2 .55 4 1 / \setlinear
\arr(2.004,.451)(2,.45)                \arr(2.004,.551)(2,.55)
\plot 0 0 6 1 /
\arr(3.006,.501)(3,.5)
\plot 0 1 4 2 /     \plot 0 1 4 3 /
\arr(2.004,1.501)(2,1.5)    \arr(2.202,2.101)(2.2,2.1)
\plot 0 2 4 2 /     \plot 0 2 4 3 /
\arr(2.21,2)(2.2,2)         \arr(2.204,2.551)(2.2,2.55)
\plot 0 3 4 2 /     \plot 0 3 4 3 /
\arr(2.204,2.45)(2.2,2.451) \arr(2.1,3)(2,3)
\plot 2 .75 4 2 /   \plot 2 1.75 4 3 /
\arr(3.208,1.505)(3.2,1.5)  \arr(3.208,2.505)(3.2,2.5)
\plot 2 .75 4 1 /   \plot 2 1.75 4 1 /
\arr(3.0008,.8751)(3,.875)  \arr(3.208,1.3)(3.2,1.303)
\plot 2 .75 6 2 /   \plot 2 1.75 6 2 /
\arr(4.208,1.44)(4.2,1.437) \arr(4.016,1.876)(4,1.875)
\plot 4 1 6 2 /
\arr(5.102,1.551)(5.1,1.55)
\setquadratic \plot 4 2 5 1.47 6 1 /    \plot 4 2 5 1.53 6 1 / \setlinear
\arr(5.102,1.42)(5.1,1.421)             \arr(5.102,1.48)(5.1,1.481)
\setquadratic \plot 4 3 5 2.05 6 1 /    \plot 4 3 5 1.95 6 1 / \setlinear
\arr(5.001,2.05)(5,2.051)                    \arr(5.001,1.95)(5,1.951)
\endpicture
\)
\caption{The skeleton of a $2$-graph satisfying condition~\eqref{ass:AT}.}\label{fig:AT 2-graph}
\end{figure}

Figure~\ref{fig:AT 2-graph} illustrates the skeleton of a $2$-graph satisfying condition~\eqref{ass:AT}; in this and all our other diagrams we draw the blue edges as solid lines and the red edges as dashed lines.

\begin{notation} \label{ntn:Ff map}
If $\Lambda$ is a $2$-graph, we write $\Ff$ for the \emph{factorisation map} defined on pairs $(\rho, \tau)$ such that $s(\rho) = r(\tau)$ by 
\[
\Ff(\rho,\tau) := \big((\rho\tau)(0, d(\tau)), (\rho\tau)(d(\tau), d(\rho\tau))\big).
\]
If $\Ff(\rho,\tau) = (\tau',\rho')$, then we write $\Ff_1(\rho,\tau) := \tau'$ and $\Ff_2(\rho,\tau) := \rho'$. So: $\Ff_1(\rho,\tau)$ has the same degree as $\tau$ and the same range as $\rho$; $\Ff_2(\rho,\tau)$ has the same degree as $\rho$ and the same range as $\tau$; and $\Ff_1(\rho,\tau)\Ff_2(\rho,\tau) = \rho\tau$.
\end{notation}

\begin{lemma}\label{lem:cycle factorisations}
Let $(\Lambda,d)$ be a $2$-graph which satisfies condition~\eqref{ass:AT}. Suppose that $\lambda_1$ and $\lambda_2$ are isolated cycles in $\RED\Lambda$. Let $V_1$ and $V_2$ denote the vertices on $\lambda_1$ and $\lambda_2$ respectively.
\begin{itemize}
\item[(1)] For $\mu \in V_1 (\RED\Lambda)$, the map $\Ff_1(\mu,\cdot) : s(\mu) (\BLUE\Lambda) V_2 \to r(\mu)(\BLUE\Lambda) V_2$ obtained from the factorisation map above is a degree-preserving bijection, and
\item[(2)] for $\nu \in V_2 (\RED\Lambda)$, the map $\Ff_2(\cdot,\nu) : V_1 (\BLUE\Lambda) r(\nu) \to V_1 (\BLUE\Lambda) s(\nu)$ obtained from the factorisation map above is also a degree-preserving bijection.
\end{itemize}
\end{lemma}
\begin{proof}
We just prove statement~(1); statement~(2) follows from a very similar argument.

The map is degree-preserving by definition. Since $r(\mu\sigma) = r(\mu)$, the image of $\Ff_1(\mu,\cdot)$ is a subset of $r(\mu)\Lambda$. For $\sigma \in s(\mu) (\BLUE\Lambda) V_2$, we have that $(\nu\sigma) = \Ff_1(\nu,\sigma) \Ff_2(\nu,\sigma)$ where $\Ff_2(\nu,\sigma) \in \Lambda^{d(\nu)} s(\sigma)$. Now $s(\sigma)$ lies on the isolated cycle $\lambda_2$ in $\RED\Lambda$. Since $\Ff_2(\nu,\sigma) \in \RED\Lambda$, it follows that $s(\Ff_1(\nu,\sigma)) = r(\Ff_2(\nu,\sigma)) \in V_2$ as well. Hence $\Ff_1(\nu,\sigma)$ belongs to $r(\nu)(\BLUE\Lambda) V_2$. 

To see that $\Ff_1(\nu,\cdot)$ is bijective, note that for $\sigma' \in r(\nu) \Lambda V_2$, there is a unique path $\nu'(\sigma') \in s(\sigma')\Lambda^{d(\nu)}$ because $\lambda_2$ is isolated. Hence $\sigma' \mapsto \Ff_2(\sigma',\nu'(\sigma'))$ is an inverse for $F_1(\nu,\cdot)$.
\end{proof}

\begin{notation} \label{ntn:mu's and unitaries}
\begin{itemize}
\item[(1)] Let $(\Lambda,d)$ be a $2$-graph which satisfies condition~\eqref{ass:AT}. If $e$ is a red edge and $\mu$ is a red path, then $\langle \mu,e\rangle$ is the number of occurrences of $e$ in $\mu$. 
\item[(2)] Since lower case $e$'s already denote generators of $\NN^k$ and edges of $k$-graphs, we do not use them to denote matrix units. If $I$ is an index set, we write $\{\Theta(i,j) : i,j \in I\}$ for the canonical matrix units in $M_I(\CC)$; we use $\{\theta(i,j) : i,j \in I\}$ to denote a set of matrix units in some other $C^*$-algebra. 
\item[(3)] We write $z$ for the monomial $z \mapsto z \in C(\TT)$, and $z^n$ for $z \mapsto z^n$.
\end{itemize}
\end{notation}

\begin{prop}\label{prp:M_P(C(T))}
Let $(\Lambda,d)$ be a finite $2$-graph which satisfies condition~\eqref{ass:AT}. Suppose that the set $S$ of sources in $\BLUE\Lambda$ are the vertices on a single isolated cycle in $\RED\Lambda$. Let $Y$ denote the collection $(\BLUE\Lambda)S$ of blue paths with source in $S$. Fix an edge $e_\star$ in $S \Lambda^{e_2} S$. Then there is an isomorphism $\pi : C^*(\Lambda) \to M_Y(C(\TT))$ such that for $\alpha,\beta \in Y$ and $\nu \in s(\alpha)(\RED\Lambda)s(\beta)$,
\begin{equation}\label{eq:formula for pi}
\pi(s_\alpha s_\nu s^*_\beta)(z) = z^{\langle \nu,e_\star \rangle} \Theta(\alpha,\beta).
\end{equation}
\end{prop}

The proof of the proposition is long, but not particularly difficult. The strategy is to identify a family of nonzero matrix units $\{\theta(\alpha,\beta) : \alpha,\beta \in Y\}$ and a unitary $U$ in $C^*(\Lambda)$ such that $U$ commutes with each $\theta(\alpha,\beta)$. This gives a homomorphism $\phi$ of $M_Y(\CC) \otimes \TT$ into $C^*(\Lambda)$. Since the matrix units are nonzero, we just have to show that $U$ has full spectrum to see that $\phi$ is injective. We then show that $\phi$ is also surjective, and that taking $\pi = \phi^{-1}$ gives an isomorphism satisfying~\eqref{eq:formula for pi}.

We proceed in a series of lemmas. We begin with a simple technical lemma which we will use frequently.

\begin{lemma}\label{lem:mpctn in f2}
Let $(\Lambda,d)$ be a $2$-graph which satisfies condition~\eqref{ass:AT}, and let $\mu,\nu$ be red paths in $\Lambda$.
\begin{itemize}
\item[(1)] If $r(\mu) \not = r(\nu)$ then $s^*_\mu s_\nu = 0$; otherwise either $\nu = \mu\nu'$ and $s^*_\mu s_\nu = s_{\nu'}$ or $\mu = \nu\mu'$ and $s^*_\mu s_\nu = s^*_{\mu'}$.
\item[(2)] If $s(\mu) \not = s(\nu)$ then $s_\mu s^*_\nu = 0$; otherwise either $\mu = \mu'\nu$ and $s_\mu s^*_\nu = s_{\mu'}$ or $\nu = \nu'\mu$ and $s_\mu s^*_\nu = s^*_{\nu'}$.
\end{itemize}
\end{lemma}
\begin{proof}
Relation~(CK1) shows that $s^*_\mu s_\nu = 0$ when $r(\mu) \not= r(\nu)$ and that $s_\mu s^*_\nu = 0$ when $s(\mu) \not= s(\nu)$. Because $\RED\Lambda$ is a union of isolated cycles, we have either $\mu = \nu\mu'$ or $\nu = \mu\nu'$ when $r(\mu) = r(\nu)$, and either $\mu = \mu'\nu$ or $\nu = \nu'\mu$ when $s(\mu) = s(\nu)$. So~(1) follows from~(CK3), and~(2) follows from~(CK4) because $r(\mu)\Lambda^{\le d(\mu)} = r(\mu)\Lambda^{d(\mu)} = \{\mu\}$ for all $\mu \in \RED\Lambda$.
\end{proof}

\begin{lemma}\label{lem:spanning}
Suppose that $(\Lambda,d)$ is a finite $2$-graph which satisfies condition~\eqref{ass:AT}. Let $S$ denote the collection of sources in $\BLUE\Lambda$, and let $Y := (\BLUE\Lambda)S$. Then 
\[
C^*(\Lambda) = \clsp\{s_\alpha s_\nu s^*_\beta, s_\alpha s^*_\nu s^*_\beta : \alpha,\beta \in Y, \nu \in S (\RED\Lambda) S\}.
\]
\end{lemma}
\begin{proof}
By \cite[Remark~3.8(1)]{RSY1}, 
\[
C^*(\Lambda) = \clsp\{s_\rho s^*_\xi : \rho, \xi \in \Lambda, s(\rho) = s(\xi)\},
\]
so it suffices to show that for all $\rho,\xi$ in $\Lambda$ with $s(\rho) = s(\xi)$, we can write $s_\rho s^*_\xi$ as a sum of elements of the form $s_\alpha s_\nu s^*_\beta$ or $s_\alpha s^*_\nu s^*_\beta$. 

Fix $\rho, \xi \in \Lambda$ with $s(\rho) = s(\xi)$. Since $\BLUE\Lambda$ has no cycles, and $\Lambda^0$ is finite, there exists $N \in \NN$ such that $\Lambda^{n e_1} = \emptyset$ for all $n \ge N$. Thus $\Lambda^{\le N e_1} = (\BLUE\Lambda)S = Y$. Hence~(CK4) gives $s_\rho s^*_\xi = \sum_{\alpha \in s(\rho)Y} s_{\rho\alpha} s^*_{\xi\alpha}$. For fixed $\alpha \in s(\rho)Y$, we factorise $\rho\alpha = \eta\mu$ and $\xi\alpha = \zeta\nu$ where $\eta,\zeta \in \BLUE\Lambda$ and $\mu,\nu \in \RED\Lambda$. Since $\RED\Lambda$ consists of isolated cycles, we must have $s(\eta), s(\zeta) \in S$, so $\eta,\zeta \in Y$, and $\mu,\nu \in S (\RED\Lambda) S$. We have $s(\mu) = s(\eta\mu) = s(\rho\alpha) = s(\alpha) = s(\xi\alpha) = s(\zeta\nu) = s(\nu)$. Hence Lemma~\ref{lem:mpctn in f2} shows that either $s_{\rho\alpha} s^*_{\xi\alpha} = s_\eta s_{\mu'} s^*_\zeta$, or $s_{\rho\alpha} s^*_{\xi\alpha} = s_\eta s^*_{\nu'} s^*_\zeta$, where $\mu', \nu' \in S (\RED\Lambda) S$.
\end{proof}

\begin{lemma}\label{lem:MUs}
Let $(\Lambda,d)$ be a finite $2$-graph which satisfies condition~\eqref{ass:AT}. Suppose that the set $S$ of sources in $\BLUE\Lambda$ is the set of vertices on a single isolated cycle in $\RED\Lambda$. Let $Y := (\BLUE\Lambda)S$. Fix an edge $e_\star$ in $S \Lambda^{e_2} S$, and for $\alpha,\beta \in Y$, let $\nu_0(\alpha,\beta)$ be the unique path connecting $s(\alpha)$ and $s(\beta)$ (in either direction) such that $\langle \nu_0(\alpha,\beta),e_\star \rangle = 0$. Then the elements 
\[
\theta(\alpha,\beta) := \begin{cases}
s_\alpha s_{\nu_0(\alpha,\beta)} s^*_\beta &\text{ if $s(\nu_0(\alpha,\beta)) = s(\beta)$} \\
s_\alpha s^*_{\nu_0(\alpha,\beta)} s^*_\beta &\text{ if $s(\nu_0(\alpha,\beta)) = s(\alpha)$}
\end{cases}
\]
form a collection of nonzero matrix units in $C^*(\Lambda)$.
\end{lemma}
\begin{proof}
The $\theta(\alpha,\beta)$ are nonzero because the partial isometries $\{s_\rho : \rho \in \Lambda\}$ are all nonzero by \cite[Theorem~3.15]{RSY1}. By the definition of the $\theta(\alpha,\beta)$ we have $\theta(\alpha,\beta)^* = \theta(\beta,\alpha)$. So we need only show that $\theta(\alpha,\beta)\theta(\eta,\zeta) = \delta_{\beta,\eta} \theta(\alpha,\zeta)$. Now $\beta,\eta \in \BLUE\Lambda$, and they begin at sources in $\BLUE\Lambda$; it follows that $s_\beta^* s_\eta = \delta_{\beta,\eta} s_{s(\beta)}$. So it suffices to show that $\theta(\alpha,\beta) \theta(\beta,\zeta) = \theta(\alpha,\zeta)$.

Now $\theta(\alpha,\beta) \theta(\beta,\zeta) = s_\alpha t_1 s^*_\beta s_\beta t_2 s^*_\zeta = s_\alpha t_1 t_2 s^*_\zeta$ where $t_1 \in \{s_{\nu_0(\alpha,\beta)}, s^*_{\nu_0(\alpha,\beta)}\}$ and $t_2 \in \{s_{\nu_0(\beta,\zeta)}, s^*_{\nu_0(\beta,\zeta)}\}$. Hence there are four cases to consider. We will deal with two of them; the other two arguments are very similar.

If $t_1 = s_{\nu_0(\alpha,\beta)}$ and $t_2 = s_{\nu_0(\beta,\zeta)}$, then $t_1t_2 = s_{\nu_0(\alpha,\beta)\nu_0(\beta,\zeta)}$. Since neither $\nu_0(\alpha,\beta)$ nor $\nu_0(\beta,\zeta)$ contains an instance of $e_\star$, neither does $\nu = \nu_0(\alpha,\beta)\nu_0(\beta,\zeta)$. Hence $\nu = \nu_0(\alpha,\zeta)$ by definition, so $\theta(\alpha,\beta) \theta(\beta,\zeta) = \theta(\alpha, \zeta)$.

If $t_1 = s_{\nu_0(\alpha,\beta)}$ and $t_2 = s^*_{\nu_0(\beta,\zeta)}$, then Lemma~\ref{lem:mpctn in f2} shows that $t_1 t_2$ has the form $s_{\nu}$ or $s^*_{\nu}$ depending on which of $\nu_0(\alpha,\beta)$ and $\nu_0(\beta,\zeta)$ is longer. In either case, $\nu$ is a path joining $s(\alpha)$ and $s(\zeta)$, and since it is a sub-path of one of $\nu_0(\alpha,\beta)$ and $\nu_0(\beta,\zeta)$, neither of which contains an instance of $e_\star$, we once again have $\nu = \nu_0(\alpha,\zeta)$.
\end{proof}

\begin{lemma} \label{lem:unitary}
Let $(\Lambda,d)$ be a finite $2$-graph which satisfies condition~\eqref{ass:AT}. Suppose that the set $S$ of sources in $\BLUE\Lambda$ is the set of vertices on a single isolated cycle in $\RED\Lambda$. Let $Y := (\BLUE\Lambda)S$, and for $\alpha \in Y$, let $\lambda(\alpha)$ be the unique isolated cycle in $s(\alpha) (\RED\Lambda)$. Let $U := \sum_{\alpha \in Y} s_\alpha s_{\lambda(\alpha)} s^*_\alpha$. Then $U$ is a unitary in $C^*(\Lambda)$ and the spectrum of $U$ is $\TT$.
\end{lemma}
\begin{proof}
For $\alpha,\beta \in Y$, we have $s^*_\alpha s_\beta = \delta_{\alpha,\beta} s_{s(\alpha)}$, so
\[
U U^* 
= \sum_{\alpha,\beta \in Y} s_{\alpha} s_{\lambda(\alpha)} s^*_\alpha s_\beta s^*_{\lambda(\beta)}s^*_\beta 
= \sum_{\alpha \in Y} s_{\alpha} s_{\lambda(\alpha)} s^*_{\lambda(\alpha)} s^*_\alpha 
= \sum_{\alpha \in Y} s_\alpha s^*_\alpha
\] 
by Lemma~\ref{lem:mpctn in f2}. Since $\Lambda$ is finite, there exists $N \in \NN$ such that $\Lambda^{\le Ne_1} = Y$, and so it follows from the calculation above and~(CK4) that $U U^* = \sum_{v \in \Lambda^0} s_v = 1_{C^*(\Lambda)}$.

A similar calculation establishes that $U^* U = 1_{C^*(\Lambda)}$. It remains to show that $U$ has spectrum $\TT$. For this, notice that $d(\lambda(\alpha)) = |S| e_2$ for all $\alpha \in Y$. Hence the gauge action satisfies $\gamma_{1,y}(U) = y^{|S|}U$ for $y \in \TT$. Fix $z \in \sigma(U)$ and $w \in \TT$, and choose $y \in \TT$ such that $y^{|S|} = z\overline{w}$. We have
\begin{align*}
z \in \sigma(U) 
&\implies U - z 1_{C^*(\Lambda)} \not\in C^*(\Lambda)^{-1} \\
&\implies \gamma_{1,y}(U - z 1_{C^*(\Lambda)}) \not\in C^*(\Lambda)^{-1} \\
&\implies y^{|S|} U - z 1_{C^*(\Lambda)} \not\in C^*(\Lambda)^{-1} \\
&\implies z(\overline{w}U - 1_{C^*(\Lambda)}) \not\in C^*(\Lambda)^{-1} \\
&\implies w \in \sigma(U).
\end{align*}
Hence $\sigma(U) = \TT$.
\end{proof}

\begin{lemma} \label{lem:commute}
Let $(\Lambda,d)$ be a finite $2$-graph which satisfies condition~\eqref{ass:AT}, and suppose that the set $S$ of sources in $\BLUE\Lambda$ is the set of vertices on a single isolated cycle in $\RED\Lambda$. Let $Y := (\BLUE\Lambda)S$. Fix an edge $e_\star$ in $S \Lambda^{e_2} S$. Then the matrix units $\theta(\alpha,\beta)$ of Lemma~\ref{lem:MUs} commute with the unitary $U$ of Lemma~\ref{lem:unitary}.
\end{lemma}
\begin{proof}
Fix $\alpha,\beta \in Y$ and $n \in \NN$. Let $\nu_n(\alpha,\beta)$ be the unique path in $s(\alpha)(\RED\Lambda) s(\beta)$ such that $\langle \nu_n(\alpha,\beta), e_\star\rangle = n$. Using (CK2) and Lemma~\ref{lem:mpctn in f2}, it is easy to check that
\begin{equation}\label{eq:U times MU}
U^n \theta(\alpha,\beta) = s_\alpha s_{\nu_n(\alpha,\beta)} s^*_\beta = \theta(\alpha,\beta) U^n.
\end{equation}
Taking $n = 1$ proves the lemma.
\end{proof}

\begin{proof}[Proof of Proposition~\ref{prp:M_P(C(T))}]
Define $\theta(\alpha,\beta)$ and $U$ as in Lemmas \ref{lem:MUs}~and~\ref{lem:unitary}. The universal properties of $M_Y(\CC)$ and of $C(\TT) = C^*(\ZZ)$ ensure that there are homomorphisms $\phi_M : M_Y(\CC) \to C^*(\Lambda)$ and $\phi_\TT : C(\TT) \to C^*(\Lambda)$ such that $\phi_M(\Theta(\alpha, \beta)) = \theta(\alpha,\beta)$ and $\phi_\TT(z) = U$. We know $\phi_M$ is injective because the $\theta(\alpha,\beta)$ are all nonzero by Lemma~\ref{lem:MUs}. We know $\phi_\TT$ is injective because $U$ has full spectrum by Lemma~\ref{lem:unitary}. Since $U$ commutes with the $\theta(\alpha,\beta)$ by Lemma~\ref{lem:commute}, there is a well-defined homomorphism $\phi = \phi_M \otimes \phi_\TT : M_Y(\CC) \otimes C(\TT) \to C^*(\Lambda)$ which satisfies $\phi(\Theta(\alpha,\beta) \otimes z^n) = \theta(\alpha,\beta) U$. Moreover, $\phi$ is injective because both $\phi_M$ and $\phi_\TT$ are injective.

We claim that $\phi$ is also surjective. By Lemma~\ref{lem:spanning}, we need only show that if $\alpha,\beta \in Y$ and $\mu \in s(\alpha) (\RED\Lambda) s(\beta)$, then $s_\alpha s_\mu s^*_\beta$ belongs to the image of $\phi$: taking adjoints then shows that $s_\alpha s^*_\mu s^*_\beta$ also belongs to the image of $\phi$. A straightforward calculation shows that $s_\alpha s_\mu s^*_\beta = U^{\langle\mu, e_\star \rangle} \theta(\alpha, \beta) = \phi(\Theta(\alpha,\beta) \otimes z^{\langle\mu, e_\star\rangle})$.

Let $\pi := \phi^{-1}$. Since $\pi(\theta(\alpha,\beta)) = \Theta(\alpha,\beta)$ and $\pi(U) = z$, equation~\eqref{eq:U times MU} implies that $\pi$ satisfies~\eqref{eq:formula for pi}.
\end{proof}

\begin{prop}\label{prp:oplus}
Let $(\Lambda,d)$ be a finite $2$-graph satisfying condition~\eqref{ass:AT}. Let $S$ be the set of sources of $\BLUE\Lambda$, and write $S = S_1 \sqcup \dots \sqcup S_n$ where each $S_j$ is the set of vertices on one of the isolated cycles in $\RED\Lambda$. For $1 \le j \le n$, let $\Lambda_j := \{\lambda \in \Lambda: s(\lambda) \Lambda S_j \not= \emptyset\}$, and let $d_j$ denote the restriction of the degree map $d$ to $\Lambda_j$.
\begin{itemize}
\item[(1)] Each $(\Lambda_j, d_j)$ is a $2$-graph which satisfies condition~\eqref{ass:AT}, and $S_j = S \cap \Lambda_j$ is precisely the set of sources in $\BLUE\Lambda_j$.
\item[(2)] For $1 \le j \le n$, let $\{s_{j,\rho} : \rho \in \Lambda_j\}$ denote the universal generating Cuntz-Krieger $\Lambda_j$-family. There is an isomorphism of $C^*(\Lambda)$ onto $\bigoplus^n_{j=1} C^*(\Lambda_j)$ which carries $s_\alpha s_\mu s^*_\beta$ to $(0,\dots, 0, s_{j,\alpha} s_{j,\mu} s^*_{j,\beta}, 0,\dots,0)$ for $\alpha,\beta \in Y_j := (\BLUE\Lambda_j)S_j$ and $\mu \in s(\alpha)(\RED\Lambda_j)s(\beta)$.
\end{itemize}
\end{prop}
\begin{proof}
(1) Fix $1 \le j \le n$. To see that $\Lambda_j$ is a category, note that for $\rho \in \Lambda$, we have $\rho \in \Lambda_j$ if and only if $s(\rho) \in \Lambda_j$. Hence $\xi \in \Lambda_j$ implies $\rho\xi \in \Lambda_j$. 

To check the factorisation property, suppose that $\rho \in \Lambda_j$ and $d_j(\rho) = p+q$. The factorisation property for $\Lambda$ ensures that there is a unique factorisation $\rho = \rho'\rho''$ where $d(\rho') = p$ and $d(\rho'') = q$, so we need only show that $\rho', \rho'' \in \Lambda_j$. Since $\rho \in \Lambda_j$, there exists a path $\xi$ in $s(\rho)\Lambda S_j = s(\rho'') \Lambda S_j$, and it follows that $\rho'' \in \Lambda_j$. Moreover, $\rho''\xi \in s(\rho')\Lambda S_j$ giving $\rho' \in \Lambda_j$ as well.

Since $\Lambda_j$ is a subgraph of $\Lambda$ it is row-finite, and $\BLUE\Lambda_j$ is cycle-free. For a given isolated cycle $\lambda$ in $\RED\Lambda$, one of the vertices on $\lambda$ lies in $\Lambda_j$ if and only if they all do, and in this case all the edges in $\lambda$ belong to $\Lambda_j$ as well. Thus each vertex of $\Lambda_j$ lies on an isolated cycle in $\RED\Lambda_j$. So $\Lambda_j$ satisfies condition~\eqref{ass:AT}. 

Let $v \in S_j$ and fix $\rho \in v\Lambda$. Write $\rho = \sigma\mu$ where $\sigma \in \BLUE\Lambda$ and $\mu \in \RED\Lambda$. Since $S_j \subset S$, we must have $d(\sigma) = 0$. But now $\mu \in v (\RED\Lambda) \subset S (\RED\Lambda)$. By assumption, $S (\RED\Lambda) = \bigcup^n_{j=1} S_j (\RED\Lambda) S_j$, so $\rho \in S_j \Lambda S_j$, and $s(\rho) \in S_j$. Since the $S_j$ are disjoint, $S_j \Lambda S_l = \emptyset$ for $j \not= l$, and it follows that $S_j = S \cap \Lambda_j$.

Finally, to see that the sources in $\BLUE\Lambda_j$ are precisely $S \cap \Lambda_j$, note first that elements of $S \cap \Lambda_j$ are clearly sources in $\BLUE\Lambda_j$. For the reverse inclusion, let $v$ be a source in $\BLUE\Lambda_j$, so $v (\BLUE\Lambda_j) = \{v\}$. Since $v \in \Lambda_j$, we have $v \Lambda S_j \not= \emptyset$, say $\rho \in v \Lambda S_j$. Factorise $\rho = \sigma \mu$ where $\sigma \in \BLUE\Lambda$ and $\mu \in \RED\Lambda$. Since $\Lambda$ satisfies condition~\eqref{ass:AT}, we have $S (\RED\Lambda) = (\RED\Lambda) S$. We have $s(\mu) \in S_j$ by choice of $\rho$, and we therefore have $r(\mu) \in S_j$. But $r(\mu) = s(\sigma)$, so $\sigma \in (\BLUE\Lambda)S_j = \BLUE\Lambda_j$. But $v$ was a source in $\BLUE\Lambda_j$, giving $v = r(\sigma) = s(\sigma) \in S_j$.

(2) Define operators $\{t_\rho : \rho \in \Lambda\} \in \bigoplus^n_{j=1} C^*(\Lambda_j)$ by $t_\rho := \bigoplus_{\{j : \rho \in \Lambda_j\}} s_{j,\rho}$. Because each $\{s_{j,\rho} : \rho \in \Lambda_j\}$ is a Cuntz-Krieger family, $\{t_v : v \in \Lambda^0\}$ is a collection of mutually orthogonal projections, so $\{t_\rho : \rho \in \Lambda\}$ satisfies~(CK1). 

If $\rho,\xi \in \Lambda$ satisfy $s(\rho) = r(\xi)$, then
\begin{equation}\label{eq:concat calc}\textstyle
t_\rho t_\xi 
= \Big(\bigoplus_{\{j : \rho \in \Lambda_j\}} s_{j,\rho}\Big) \Big(\bigoplus_{\{l : \xi \in \Lambda_l\}}s_{l,\xi}\Big)
= \bigoplus_{\{j : \rho, \xi \in \Lambda_j\}} s_{j,\rho} s_{j,\xi}
= \bigoplus_{\{j : \rho,\xi \in \Lambda_j\}} s_{j,\rho\xi}.
\end{equation}
Since $\rho\xi \in \Lambda_j$ if and only if $\xi \in \Lambda_j$ and since $\xi \in \Lambda_j \implies \rho \in \Lambda_j$, the right-hand side of~\eqref{eq:concat calc} is equal to $t_{\rho\xi}$. Hence $\{t_\rho : \rho \in \Lambda\}$ satisfies~(CK2).

For $\rho \in \Lambda$, we have
\[\textstyle
t^*_\rho t_\rho 
= \Big(\bigoplus_{\{j : \rho \in \Lambda_j\}} s^*_{j,\rho}\Big) \Big(\bigoplus_{\{l : \rho \in \Lambda_l\}} s_{l,\rho}\Big)
= \bigoplus_{\{j : \rho \in \Lambda_j\}} s^*_{j,\rho} s_{j,\rho} 
= \bigoplus_{\{j : \rho \in \Lambda_j\}} s_{j,s(\rho)}.
\]
Since $\rho \in \Lambda_j$ if and only if $s(\rho) \in \Lambda_j$, $\{t_\rho : \rho \in \Lambda\}$ satisfies~(CK3).

Finally, fix $v \in \Lambda^0$. To establish that the $\{t_\rho : \rho \in \Lambda\}$ satisfy~(CK4), we need only show that $t_v = \sum_{e \in v\Lambda^{e_i}} t_e t^*_e$ for $i = 1,2$. For $i = 2$, this is easy as $v\Lambda^{e_2}$ has precisely one element $f$, and $f \in \Lambda_j$ if and only if $v \in \Lambda_j$. Hence 
\[\textstyle
t_f t^*_f = \bigoplus_{\{j : f \in \Lambda_j\}} s_{j,f} s^*_{j,f} = \bigoplus_{\{j : f \in \Lambda_j\}} s_{j,r(f)} = t_v.
\]
Now consider $i = 1$. Note that for $v \in \Lambda^0$ and $1 \le j \le n$, we have $v \in \Lambda_j$ if and only if there exists a blue edge $e \in v\Lambda^{e_1}$ such that $e \in \Lambda_j$. Using this to reverse the order of summation in the third line below, we calculate:
\begin{align*}
t_v 
&=\textstyle \bigoplus_{\{j : v \in \Lambda_j\}} s_{j,v} \\
&=\textstyle \bigoplus_{\{j : v \in \Lambda_j\}}\Big( \sum_{e \in v\Lambda_j^{e_1}} s_{j,e} s^*_{j,e} \Big)\\
&=\textstyle \sum_{e \in v\Lambda^{e_1}} \Big( \bigoplus_{\{j : e \in \Lambda_j\}} s_{j,e} s^*_{j,e}\Big) \\
&=\textstyle \sum_{e \in v\Lambda^{e_1}} t_e t^*_e,
\end{align*}
and so $\{t_\rho : \rho \in \Lambda\}$ satisfies~(CK4), and hence is a Cuntz-Krieger $\Lambda$-family.

The universal property of $C^*(\Lambda)$ gives a homomorphism $\psi_t : C^*(\Lambda) \to \bigoplus^n_{j=1} C^*(\Lambda_j)$ such that $\psi_t(s_\rho) = t_\rho$ for all $\rho \in \Lambda$. We claim that $\psi_t$ is bijective. To see that $\psi_t$ is injective, let $\gamma$ be the gauge action on $C^*(\Lambda)$, and let $\bigoplus^n_{j=1} \gamma_j$ be the direct sum of the gauge-actions on the $C^*(\Lambda_j)$. Since $d_j(\rho) = d(\rho)$ whenever $\rho \in \Lambda_j$, it is easy to see that $\psi_t \circ \gamma_z = (\bigoplus^n_{j=1} \gamma_j)_z \circ \psi_t$. Moreover, each $\rho \in \Lambda$ belongs to at least one $\Lambda_j$, and hence each $t_\rho$ is nonzero. It now follows from the gauge-invariant uniqueness theorem \cite[Theorem~4.1]{RSY1} that $\psi_t$ is injective.

Finally, we must show that $\psi_t$ is surjective. We just need to show that if $\rho \in \Lambda_j$ then $s_{j,\rho}$ belongs to the image of $\psi_t$. For this, note that if $\rho \in \Lambda_j$ satisfies $s(\rho) \in S_j$, then we must have $s(\rho)\Lambda S_l = \emptyset$ for $j \not= l$. It follows that for such $\rho$, we have $s_{j,\rho} = t_\rho$. Now let $Y_j := (\BLUE\Lambda) S_j = (\BLUE\Lambda_j) S_j$ for $1 \le j \le n$. For $\rho \in \Lambda_j$, 
\[\textstyle
s_{j,\rho} = \sum_{\alpha \in s(\rho) Y_j} s_{j,\rho\alpha} s^*_{j,\alpha} = \sum_{\alpha \in s(\rho) Y_j} t_{\rho\alpha} t^*_\alpha = \psi_t\big(\sum_{\alpha \in s(\rho)Y_j} s_{\rho\alpha} s^*_\alpha\big),
\]
belongs to the image of $\psi_t$.
\end{proof}

\begin{cor} \label{cor:oplus M(C(T))}
Let $(\Lambda,d)$ be a finite $2$-graph which satisfies condition~\eqref{ass:AT}. Let $S$ be the set of sources of $\BLUE\Lambda$, and write $S = S_1 \sqcup \dots \sqcup S_n$ where each $S_j$ is the set of vertices on one of the isolated cycles in $\RED\Lambda$. For each $j$, let $Y_j := (\BLUE\Lambda)S_j$. Then
\[\textstyle
C^*(\Lambda) \cong \bigoplus^n_{j=1} M_{Y_j}(C(\TT)).
\]
\end{cor}
\begin{proof}
Proposition~\ref{prp:oplus}(2) gives $C^*(\Lambda) \cong \bigoplus^n_{i=1} C^*(\Lambda_i)$. By Proposition~\ref{prp:oplus}(1), each $\Lambda_i$ satisfies the hypotheses of Proposition~\ref{prp:M_P(C(T))}. By Proposition~\ref{prp:M_P(C(T))}, we then have $C^*(\Lambda_j) \cong M_{Y_i}(C(\TT))$ for each $i$.
\end{proof}

\begin{proof}[Proof of Theorem~\ref{thm:AT alg}]
By \cite[Proposition~3.2.3]{Ror}, it suffices to show that for every finite collection $a_1, \dots a_n$ of elements of $C^*(\Lambda)$, and every $\varepsilon > 0$ there exists a sub $C^*$-algebra $B \subset C^*(\Lambda)$ and elements $b_1, \dots b_n \in B$ such that $B \cong \bigoplus^n_{i=1} M_{m_i}(C(\TT))$ for some $m_1, \dots m_n \in \NN$, and such that $\|a_i - b_i\| \le \varepsilon$ for $1 \le i \le n$.
 
Since $C^*(\Lambda) = \clsp\{s_\rho s^*_\xi : \rho,\xi \in \Lambda\}$, it therefore suffices to show that for every finite collection $F$ of paths in $\Lambda$, there is a sub $C^*$-algebra $B \subset C^*(\Lambda)$ such that $B \cong \bigoplus^n_{i=1} M_{m_i}(C(\TT))$ for some $m_1, \dots m_n \in \NN$, and such that $B$ contains $\{s_\rho s^*_\xi : \rho,\xi \in F\}$. Our argument is based on the corresponding argument that the $C^*$-algebra of a directed graph with no cycles is AF in \cite[Theorem~2.4]{KPR}.

Fix a finite set $F \subset \Lambda$. Build a set $G \subset \Lambda^{e_1}$ of blue edges as follows. First, let $G_1 = \{\xi(n - e_1, n) : \xi \in F, e_1 \le n \le d(\xi)\}$ be the collection of all blue edges which occur as segments of paths in $F$. Then $G_1$ is finite because $F$ is. Next, obtain $G_2$ by adding to $G_1$ all blue edges $e$ such that $r(e) = s(f)$ for some $f \in G_1$. So $G_2$ has the property that if $e \in G_2$, then either $s(e) \Lambda^{e_1} \subset G_2$ or $s(e)\Lambda^{e_1} \cap G_2 = \emptyset$. Moreover $G_2$ is finite because $\Lambda$ is row-finite. Finally, let $G$ be the collection of all blue edges obtained by applying one of the bijections $\Ff_1(\mu,\cdot)$ of Lemma~\ref{lem:cycle factorisations} to an element of $G_2$; that is $G = \{\Ff_1(\mu, e) : e \in G_2, \mu \in (\RED\Lambda) r(e)\}$. Then $G$ is finite because each isolated cycle has only finitely many vertices and $\Lambda$ is row-finite. 

Let $\Gamma$ be the subset of $\Lambda$ consisting of all vertices which are either the source or the range of an element of $G$ together with all paths of the form $\sigma\mu$ where $\sigma \in \BLUE\Lambda$ is a concatenation of edges from $G$ (that is, $\sigma(n, n+e_1) \in G$ for all $n \le d(\sigma) - e_1$), and $\mu$ is an element of $s(\sigma) (\RED\Lambda)$. By construction, for each isolated cycle $\lambda$ in $\RED\Lambda$, any one of the vertices on $\lambda$ is the source of an edge in $G$ if and only if they all are. It follows that $\Gamma$ is a category because it contains all concatenations of its elements by construction. Moreover, $\Gamma$ satisfies the factorisation property by the construction of $G$ from $G_2$; so $\Gamma$ is a sub $2$-graph of $\Lambda$. 

Since $G$ is finite, $\Gamma^0$ is finite. By construction of $G_2$, for each isolated cycle $\lambda$ in $\RED\Lambda$ whose vertices are the sources of edges in $G$, and for each vertex $v$ on $\Lambda$, either $v\Lambda^{e_1} \subset G$ or $v\Lambda^{e_1} \cap G = \emptyset$. It follows that the Cuntz-Krieger relations for $\Gamma$ are the same as those for $\Lambda$, so there is a homomorphism $\pi : C^*(\Gamma) \to C^*(\Lambda)$ such that $\pi(s_{\Gamma,\rho}) = s_{\Lambda,\rho}$ for all $\rho \in \Gamma$, where $\{s_{\Gamma,\rho} : \rho \in \Gamma\}$, and $\{s_{\Lambda,\rho} : \rho \in \Lambda\}$ are the universal generating Cuntz-Krieger families. The $s_{\Lambda,\rho}$ are automatically nonzero and $\pi$ clearly intertwines the gauge actions on $C^*(\Gamma)$ and $C^*(\Lambda)$, so $\pi$ is an isomorphism of $C^*(\Gamma)$ onto $C^*(\{s_{\Lambda,\rho} : \rho \in \Gamma\}) \subset C^*(\Lambda)$.

But $\Gamma$ satisfies the hypotheses of Proposition~\ref{cor:oplus M(C(T))} by construction, so $C^*(\Gamma) \cong \bigoplus^n_{j=1} M_{Y_j}(C(\TT))$. Hence taking $B := C^*(\{s_{\Lambda,\rho} : \rho \in \Gamma\})$ gives the required circle algebra in $C^*(\Lambda)$.
\end{proof}

\section{\Gbds/ and their $C^*$-algebras} \label{sec:GBDs}

\begin{dfn}\label{dfn:AT 2-graph}
A \emph{\gbd/ of depth $N \in \NN \cup \{\infty\}$} is a row-finite $2$-graph $(\Lambda, d)$ such that $\Lambda^0$ is a disjoint union $\bigsqcup^N_{n=0} V_n$ of nonempty finite sets which satisfy
\begin{itemize}
\item[(1)] for every blue edge $e \in \Lambda^{e_1}$, there exists $n$ such that $r(e) \in V_n$ and $s(e) \in V_{n+1}$;
\item[(2)] all vertices which are sinks in $\BLUE\Lambda$ belong to $V_0$, and all vertices which are sources in $\BLUE\Lambda$ belong to $V_N$ (where this is taken to mean that $\BLUE\Lambda$ has no sources if $N = \infty$); and
\item[(3)] every $v$ in $\Lambda^0$ lies on an isolated cycle in $\RED\Lambda$, and for each red edge $f \in \Lambda^{e_2}$ there exists $n$ such that $r(f), s(f) \in V_n$.
\end{itemize}
\end{dfn}

Conditions (1)~and~(2) say that the blue graph $\BLUE\Lambda$ is the path category of a Bratteli diagram. Condition~(3) says that each $V_n$ is itself a disjoint union $\bigsqcup^{c_n}_{j=1} V_{n,j}$ where each $V_{n,j}$ consists of the vertices on an isolated red cycle.

Every \gbd/ $\Lambda$ satisfies condition~\eqref{ass:AT}, and hence Theorem~\ref{thm:AT alg} implies that $C^*(\Lambda)$ is an A$\TT$ algebra. The inductive structure will give us an inductive limit decomposition which we use to obtain a very detailed description of the internal structure of $C^*(\Lambda)$ including $K$-invariants, ideal structure and real rank. To describe the $K$-invariants, we first need a technical lemma.

\begin{lemma} \label{lem:A,B}
Let $(\Lambda,d)$ be a \gbd/ of depth $N$. Decompose $\Lambda^0 = \bigcup^N_{n = 1} \bigcup^{c_n}_{j=1} V_{n,j}$ as above. For $0 \le n < N$, $1 \le j \le c_n$ and $1 \le i \le c_{n+1}$
\begin{enumerate}
\item the sets $v \Lambda^{e_1} V_{n+1,i}$, $v \in V_{n,j}$ have the same cardinality $A_n(i,j)$;
\item the sets $V_{n,j} \Lambda^{e_1} w$,  $w \in V_{n+1,i}$ have the same cardinality $B_n(i,j)$; and 
\item the integers $A_n(i,j)$ and $B_n(i,j)$ satisfy 
\begin{equation}\label{eq:A,B relation}
A_n(i,j) |V_{n,j}| = |V_{n,j} \Lambda^{e_1} V_{n+1, i}| = |V_{n+1,i}| B_n(i,j).
\end{equation}
\end{enumerate}
The resulting matrices $A_n, B_n \in M_{c_{n+1},c_n}(\ZZ_+)$ have no zero rows or columns. For $0 < n \le N$, let $T_n \in M_{c_n}(\ZZ_+)$ be the diagonal matrix $T_n(j,j) = |V_{n,j}|$. Then $A_n T_n = T_{n+1} B_n$ for $0 \le n < N$.
\end{lemma}
\begin{proof}
For statement~(1), fix two vertices $v,w \in V_{n,j}$. Let $\mu$ be the segment of the isolated cycle round $V_{n,j}$ from $v$ to $w$. Lemma~\ref{lem:cycle factorisations}(1) implies that the factorisation map $\Ff_1(\mu,\cdot)$ restricts to a bijection between $v\Lambda^{e_1} V_{n+1, i}$ and $w\Lambda^{e_1} V_{n+1, i}$. Statement~(2) follows in a similar way from Lemma~\ref{lem:cycle factorisations}(2).

Parts (1)~and~(2) now show that $A_n(i,j) |V_{n,j}|$ and $|V_{n+1,i}| B_n(i,j)$ are each equal to the number $|V_{n,j}\, \Lambda^{e_1}\, V_{n+1, i}|$ of blue edges with source in $V_{n+1,i}$ and range in $V_{n,j}$. This establishes~\eqref{eq:A,B relation}.

Equation~\eqref{eq:A,B relation} shows that $A_n(i,j) = 0$ if and only if $B_n(i,j) = 0$. By (1)~and~(2), the sum of the entries of the $i^{\rm th}$ row of $A_n$ is equal to $|\Lambda^{e_1} w|$ for any $w \in V_{n+1,i}$ and the sum of the entries of the $j^{\rm th}$ column is $|v\Lambda^{e_1}|$ for any $v \in V_{n,j}$. It follows from Definition~\ref{dfn:AT 2-graph}(2) that $A_n$, and hence also $B_n$, has no zero rows or columns.

The last statement follows from~\eqref{eq:A,B relation}.
\end{proof}

Let $\Lambda$ be a \gbd/. We refer to the integers $c_n$ together with the matrices $A_n, B_n$ and $T_n$ arising from Lemma~\ref{lem:A,B} as \emph{the data associated to $\Lambda$}, and say that $\Lambda$ is a \gbd/ \emph{with data $c_n, A_n, B_n, T_n$}.

\begin{theorem} \label{thm:AT and K-th}
Suppose that $(\Lambda,d)$ is a \gbd/ of infinite depth with data $c_n, A_n, B_n, T_n$. Then 
\begin{itemize}
\item[(1)] $C^*(\Lambda)$ is an A$\TT$ algebra;
\item[(2)] $K_0(C^*(\Lambda))$ is order-isomorphic to $\varinjlim (\ZZ^{c_n}, A_n)$ and there is a group isomorphism of $K_1(C^*(\Lambda))$ onto $\varinjlim(\ZZ^{c_n}, B_n)$;
\item[(3)] $K_1(C^*(\Lambda))$ is isomorphic to a subgroup $H$ of $K_0(C^*(\Lambda))$ such that the quotient group $K_0(C^*(\Lambda))/H$ has rank zero as an abelian group;
\item[(4)] the map $I \mapsto \langle [s_v] \vert s_v \in I\rangle \subset K_0(C^*(\Lambda))$ is an isomorphism of the lattice of gauge-invariant ideals of $C^*(\Lambda)$ onto the lattice of order-ideals of $K_0(C^*(\Lambda))$.
\end{itemize}
\end{theorem}

Theorem~\ref{thm:AT and K-th}(1) follows from Theorem~\ref{thm:AT alg}. We will show next that statement~(3) of Theorem~\ref{thm:AT and K-th} follows from statement~(2).

\begin{proof}[Proof of Theorem~\ref{thm:AT and K-th}(3)]
Suppose for now that Theorem~\ref{thm:AT and K-th}(2) holds. By Lemma~\ref{lem:A,B}, we have the following commuting diagram.
\begin{equation}\label{eq:comm diagram}
\beginpicture
\setcoordinatesystem units <2.5em, 2.5em>
\put{$\ZZ^{c_n}$} at 0 .8
\put{$\ZZ^{c_n}$} at 0 -.8
\put{$\ZZ^{c_{n+1}}$} at 3 .8
\put{$\ZZ^{c_{n+1}}$} at 3 -.8
\put{$\varinjlim(\ZZ^{c_n}, A_n) \cong K_0(C^*(\Lambda))$}[l] at 6.4 .8
\put{$\varinjlim(\ZZ^{c_n}, B_n) \cong K_1(C^*(\Lambda))$}[l] at 6.4 -.8
\arrow <0.15cm> [0.25,0.75] from 0  -.5 to 0 .5
\put{$T_n$} at -.3 0
\arrow <0.15cm> [0.25,0.75] from 3  -.5 to 3 .5
\put{$T_{n+1}$} at 2.5 0
\setdashes <0.417em>
\arrow <0.15cm> [0.25,0.75] from 7.5  -.5 to 7.5 .5
\put{$T_{\infty}$} at 7.1 0
\setquadratic
\plot 0 1.1 3.75 1.7 7.5 1.1 /
\arrow <0.15cm> [0.25,0.75] from 7.4 1.14 to 7.5 1.1
\put{$A_{\infty,n}$}[t] at 3.75 1.65
\plot 0 -1.1 3.75 -1.7 7.5 -1.1 /
\arrow <0.15cm> [0.25,0.75] from 7.4 -1.14 to 7.5 -1.1
\put{$B_{\infty,n}$}[b] at 3.75 -1.65
\setlinear
\setsolid
\arrow <0.15cm> [0.25,0.75] from -1.5  -.8 to -.5 -.8
\arrow <0.15cm> [0.25,0.75] from -1.5  .8 to -.5 .8
\arrow <0.15cm> [0.25,0.75] from .5  -.8 to 2.3 -.8
\put{$A_n$}[b] at 1.5 .85
\arrow <0.15cm> [0.25,0.75] from .5  .8 to 2.3 .8
\put{$B_n$}[t] at 1.5 -.85
\arrow <0.15cm> [0.25,0.75] from 3.5  -.8 to 4.4 -.8
\arrow <0.15cm> [0.25,0.75] from 3.5  .8 to 4.4 .8
\arrow <0.15cm> [0.25,0.75] from 5.3  -.8 to 6.2 -.8
\arrow <0.15cm> [0.25,0.75] from 5.3  .8 to 6.2 .8
\put{$\cdots$} at -1.85 -.8 
\put{$\cdots$} at -1.85 .8 
\put{$\cdots$} at 4.95 -.8 
\put{$\cdots$} at 4.95 .8 
\endpicture
\end{equation}
The induced map $T_{\infty}$ is injective because each $T_n$ is. 

Let $H = T_\infty(K_1(C^*(\Lambda)))$. We must show that every element of $K_0(C^*(\Lambda))/H$ has finite order. Fix $g \in K_0(C^*(\Lambda))$, and fix $n \in \NN$ and $m \in \ZZ^{c_n}$ for which $g = A_{\infty, n}(m)$. Let $t \in \NN$ be the least common multiple of the diagonal entries of $T_n$. For each $j$, $(t m)_j$ is divisible by $T_n(j,j)$, and it follows that $tm \in T_n \ZZ^{c_n}$, say $tm = T_n p$. Now $tg = t A_{\infty, n}(m) = A_{\infty, n}(t m) = A_{\infty, n}\circ T_n(p) = T_\infty \circ B_{\infty,n}(p) \in H$. It follows that the order of $[g]$ in $K_0(C^*(\Lambda))/H$ is at most $t$. Since $g \in K_0(C^*(\Lambda))$ was arbitrary, it follows that every element of $K_0(C^*(\Lambda))/H$ has finite order as required.
\end{proof}

Next we prove Theorem~\ref{thm:AT and K-th}(2). To do so, we need some technical results.

\begin{lemma}\label{lem:P}
Let $(\Lambda,d)$ be a \gbd/ of depth $N$. Let $X := V_0 \BLUE\Lambda V_N$. Then $X = \Lambda^{N e_1} = V_0 \Lambda^{\le N e_1}$. The projection $P := \sum_{v \in V_0} s_v$ is full in $C^*(\Lambda)$, and is equal to $\sum_{\alpha \in X} s_\alpha s^*_\alpha$.
\end{lemma}
\begin{proof}
That $X$ is equal to $\Lambda^{N e_1}$ follows from property~(1) of \gbds/, and that this is also equal to $V_0 \Lambda^{\le N e_1}$ follows from property~(2). To see that $P$ is full, fix $\xi \in \Lambda$. By Definition~\ref{dfn:AT 2-graph}(2), there exists $\alpha \in V_0 \Lambda r(\xi)$. Hence $s_\xi = s^*_\alpha s_\alpha s_\xi = s^*_\alpha P s_\alpha s_\xi \in C^*(\Lambda) P C^*(\Lambda)$. It follows that the generators of $C^*(\Lambda)$ belong to the ideal generated by $P$, so $P$ is full. For the final statement, fix $v \in V_0$. The first statement of the Lemma gives $vX = v\Lambda^{\le N e_1}$. Hence $s_v = \sum_{\alpha \in vX} s_\alpha s^*_\alpha$ by~(CK4). It follows that 
\[
P = \sum_{v \in V_0} s_v = \sum_{v \in V_0} \sum_{\alpha \in vX} s_\alpha s^*_\alpha = \sum_{\alpha \in X} s_\alpha s^*_\alpha.
\]
because $X = \bigsqcup_{v \in V_0} vX$.
\end{proof}

\begin{lemma}\label{lem:corner alg}
Let $(\Lambda,d)$ be a \gbd/ of depth $N$ such that the sources in $\BLUE\Lambda$ all lie on a single isolated cycle in $\RED\Lambda$. Let $P$ be the projection $P = \sum_{v \in V_0} s_v$, and let $X := V_0 (\BLUE\Lambda) V_N$. For each edge $e_\star \in V_N \Lambda^{e_2}$, there is an isomorphism $\pi : P C^*(\Lambda) P \to M_X(C(\TT))$ such that
\begin{equation}\label{eq:pi formula}
\pi(s_\alpha s_\mu s^*_\beta)(z) = z^{\langle\mu,e_\star\rangle} \Theta(\alpha,\beta) \quad\text{for all $\alpha,\beta \in X$ and $\mu \in s(\alpha) (\RED\Lambda) s(\beta)$}
\end{equation}
where $\langle \mu,e_\star \rangle$ is the number of occurrences of $e_\star$ in $\mu$ as in Notation~\ref{ntn:mu's and unitaries}. 
\end{lemma}
\begin{proof}
Let $Y$ be the set $(\BLUE\Lambda)V_N$ of paths appearing in Proposition~\ref{prp:M_P(C(T))}, which then gives an isomorphism $\pi : C^*(\Lambda) \to M_{Y}(C(\TT))$. By Lemma~\ref{lem:P}, $P$ is full. Restricting $\pi$ to $P C^*(\Lambda) P$ gives an injection, also called $\pi$ which satisfies~\eqref{eq:pi formula} and has range $M_X(C(\TT)) \subset M_Y(C(\TT))$.
\end{proof}

\begin{cor}\label{cor:corner K-th}
Let $(\Lambda,d)$, $P$ and $X$ be as in Lemma~\ref{lem:corner alg}.
\begin{itemize}
\item[(1)] There is an isomorphism $\phi_0 : K_0(P C^*(\Lambda) P) \to \ZZ$ such that $\phi_0([s_\alpha s^*_\alpha]) = 1$ for every $\alpha \in X$; and
\item[(2)] For $\alpha \in X$, let $\lambda(\alpha)$ be the isolated cycle in $\RED\Lambda$ whose range and source are equal to $s(\alpha)$. Then there is an isomorphism $\phi_1 : K_1(P C^*(\Lambda) P) \to \ZZ$ such that 
\[\textstyle
\phi_1\big(\big[ s_\alpha s_\lambda(\alpha) s^*_\alpha + \sum_{\beta \in X \setminus\{\alpha\}} s_\beta s^*_\beta\big]\big) = 1
\]
for every $\alpha \in X$.
\end{itemize}
\end{cor}
\begin{proof}
The rank-1 projections $\Theta(\alpha,\alpha) \in M_X(C(\TT))$ all represent the same class in $K_0(M_X(C(\TT)))$, and this class is the identity of $K_0$. Likewise the unitaries $z \mapsto z\Theta(\alpha,\alpha) + \sum_{\beta \not= \alpha} \Theta(\beta,\beta)$ all have the same class in $K_1(M_X(C(\TT)))$ and this class is the identity of $K_1$. The result therefore follows from Lemma~\ref{lem:corner alg}.
\end{proof}

\begin{lemma} \label{lem:direct sum K-th}
Let $(\Lambda,d)$ be a \gbd/ of depth $N$ and write $V_N = \bigsqcup^{c_n}_{j=1} V_{N,j}$ as before. Let $Y = (\BLUE\Lambda) V_N$ and $X = V_0 (\BLUE\Lambda) V_N$ as before, and let $X_j := V_0 (\BLUE\Lambda)V_{N,j}$ for $1 \le j \le c_N$. Let $P = \sum_{v \in V_0} s_v$.
\begin{itemize}
\item[(1)] There is an isomorphism $\phi_0 : K_0(P C^*(\Lambda) P) \to \ZZ^{c_n}$ such that $\phi_0([s_\alpha s^*_\alpha]) = e_j$ for every $\alpha \in X_j$.
\item[(2)] For each $\alpha \in X$, let $\lambda(\alpha)$ be the isolated cycle in $s(\alpha) (\RED\Lambda) s(\alpha)$. There is an isomorphism $\phi_1 : K_1(P C^*(\Lambda) P) \to \ZZ^{c_n}$ such that
\begin{equation}\label{eq:K1 image}\textstyle
\phi_1\big(\big[ s_{\alpha} s_{\lambda(\alpha)} s^*_\alpha + \sum_{\beta \in X_j \setminus \{\alpha\}} s_\beta s^*_\beta + \sum_{\beta \in X_k, k \not= j} s_\beta s^*_\beta \big]\big) = e_j
\end{equation}
for every $\alpha \in X_j$.
\end{itemize}
\end{lemma}
\begin{proof}
By Lemma~\ref{lem:P}, we have $P = \sum_{\alpha \in X} s_\alpha s^*_\alpha = \sum^{c_n}_{j=1} \sum_{\alpha \in X_j} s_\alpha s^*_\alpha$. Let $P_j := \sum_{\alpha \in X_j} s_\alpha s^*_\alpha$ for each $j$. The $P_j$ are mutually orthogonal, and hence $P C^*(\Lambda) P = \bigoplus^{c_n}_{j=1} P_j C^*(\Lambda) P_j$. We saw in Proposition~\ref{prp:oplus} that if $\Lambda_j := \{\rho \in \Lambda : s(\rho)\Lambda V_{N,j} \not= \emptyset\}$, then $\Lambda_j$ is a \gbd/ in which the sources in $\BLUE\Lambda_j$ lie on a single isolated cycle in $\RED\Lambda_j$, and there is an injection $\iota_j $ of $C^*(\Lambda_j) = C^*(\{t_\lambda\})$ into $C^*(\Lambda)$ which carries $t_\alpha t_\mu t^*_\beta$ to $s_\alpha s_\mu s^*_\beta$ for $\alpha,\beta \in Y_j := YV_{N,j}$. For each $j$, let $P(\Lambda_j) \in C^*(\Lambda_j)$ be the projection obtained by applying Lemma~\ref{lem:P}. Then each injection $\iota_j$ restricts to an isomorphism of $P(\Lambda_j) C^*(\Lambda_j) P(\Lambda_j)$ onto $P_j C^*(\Lambda) P_j$.

Let $\phi^j_i : K_i(P_j C^*(\Lambda) P_j) \to \ZZ$ be the isomorphisms obtained from Corollary~\ref{cor:corner K-th}, so for $\alpha \in X_j$, we have $\phi^j_0([s_\alpha s^*_\alpha]) = 1$ and 
\[\textstyle
\phi^j_1\big(\big[s_\alpha s_{\lambda(\alpha)} s^*_\alpha + \sum_{\beta \in X_j \setminus\{\alpha\}} s_\beta s^*_\beta\big]\big) = 1.
\]

These injections combine to give the desired isomorphisms $\phi_* := \bigoplus^{c_n}_{k=1} \phi^k_*$ from
\[\textstyle
K_*(P C^*(\Lambda) P) = \bigoplus^{c_n}_{k=1} K_*(P_k C^*(\Lambda) P_k)\quad\text{onto}\quad \bigoplus^{c_n}_{k=1} \ZZ = \ZZ^{c_n}.
\]
To see that this satisfies~\eqref{eq:K1 image} for each $j$, fix $1 \le j \le c_N$. For each $k$, the class $\big[\sum_{\beta \in X_k} s_\beta s^*_\beta\big]_1$ is the identity element of the direct summand $K_1(P_k C^*(\Lambda) P_k$. Hence the class of the unitary
\[
s_{\alpha} s_{\lambda(\alpha)} s^*_\alpha + \sum_{\beta \in X_j \setminus \{\alpha\}} s_\beta s^*_\beta + \sum_{\beta \in X_k, k \not= j} s_\beta s^*_\beta
\]
represents the generator of the $j^{\rm th}$ summand of $\bigoplus^{c_N}_{k = 1} K_1(P_k C^*(\Lambda) P_k)$ and the zero element in the other summands.
\end{proof}

We now consider a \gbd/ $\Lambda$ of infinite depth. For each $N \in \NN$, we consider the sub $2$-graph $\Lambda_N := \big(\bigcup^N_{n = 0} V_n\big) \Lambda \big(\bigcup^N_{n = 0} V_n\big)$ consisting of paths connecting vertices in the first $N$ levels of $\Lambda^0$ only.

\begin{lemma} \label{lem:direct limits}
Let $\Lambda$ be a \gbd/ of infinite depth. For each $N \in \NN$, the subalgebra $C^*(\{s_\rho : \rho \in \Lambda_N\})$ of $C^*(\Lambda)$ is canonically isomorphic to $C^*(\Lambda_N)$. If we use these isomorphisms to identify the $C^*(\Lambda_N)$ with the subalgebras of $C^*(\Lambda)$, then
\begin{equation}\label{eq:direct limit}
\textstyle C^*(\Lambda) = \overline{\bigcup^\infty_{N=1} C^*(\Lambda_N)}.
\end{equation}
Moreover, $P = \sum_{v \in V_0} s_v$ is a full projection in $C^*(\Lambda)$ and in each $C^*(\Lambda_N)$, and
\begin{equation}\label{eq:corner direct limit}
\textstyle P C^*(\Lambda) P = \overline{\bigcup^\infty_{N = 1} P C^*(\Lambda_N) P}.
\end{equation}
\end{lemma}
\begin{proof}
For each $N \in \NN$, the $2$-graph $\Lambda_N$ is locally convex in the sense of~\cite{RSY1}, and the elements $\{s_\rho : \rho \in \Lambda_N\}$ form a Cuntz-Krieger $\Lambda_N$-family in $C^*(\Lambda)$. Thus an application of the gauge-invariant uniqueness theorem \cite[Theorem~4.1]{RSY1} gives the required isomorphism of $C^*(\Lambda_N)$ onto $C^*(\{s_\rho : \rho \in \Lambda_N\})$. Since each generator of $C^*(\Lambda)$ lies in some $C^*(\Lambda_N)$, we have~\eqref{eq:direct limit}.

The projection $P$ is full in each $C^*(\Lambda_N)$ by Lemma~\ref{lem:P}. Hence the ideal generated by $P$ in $C^*(\Lambda)$ contains the dense subalgebra $\bigcup^\infty_{N=1} C^*(\Lambda_N)$ and it follows that $P$ is full in $C^*(\Lambda)$. Since compression by $P$ is continuous, \eqref{eq:corner direct limit} follows from~\eqref{eq:direct limit}.
\end{proof}

\begin{proof}[Proof of Theorem~\ref{thm:AT and K-th}(2)] \label{pg:pf of main thm(2)}
Lemma~\ref{lem:direct limits} implies that the projection $P = \sum_ {v \in V_0} s_v$ is full. Hence the inclusion $P C^*(\Lambda) P \subset C^*(\Lambda)$ induces isomorphisms on $K$-theory. The direct limit decomposition of $P C^*(\Lambda) P$ in~\eqref{eq:corner direct limit} and the continuity of $K$-theory imply that
\[
K_*(C^*(\Lambda)) = K_*(P C^*(\Lambda) P)) = \varinjlim K_*(P C^*(\Lambda_N) P).
\]
It therefore suffices to show that the inclusions $i_N : C^*(\Lambda_N) \hookrightarrow C^*(\Lambda_{N+1})$ and the isomorphisms $\phi^N_*$ of $K_*(P C^*(\Lambda_N) P)$ with $\ZZ^{c_N}$ provided by Lemma~\ref{lem:direct sum K-th} fit into commutative diagrams
\begin{equation}\label{eq:commuting diagrams}
\beginpicture
\put{\hbox{\(
\beginpicture
\setcoordinatesystem units <2.5em,2.5em>
\put{$K_0(P C^*(\Lambda_N) P)$}[r] at 0 2          \arr(0.1,2)(1.4,2)            \put{$\ZZ^{c_N}$}[l] at 1.5 2
    \arr(-1.5,1.7)(-1.5,0.3)                                                          \arr(1.7,1.7)(1.7,0.3)
\put{$K_0(P C^*(\Lambda_{N+1}) P)$}[r] at 0 0        \arr(0.1,0)(1.4,0)            \put{$\ZZ^{c_{N+1}}$}[l] at 1.5 0
                        \put{$\phi^N_0$} [b] at 0.75 2.1
\put{$(i_N)_*$} [r] at  -1.6 1              \put{$A_N$} [r] at 1.6 1
                        \put{$\phi^{N+1}_0$} [t] at 0.75 -0.1
\setcoordinatesystem units <2.5em,2.5em> point at -7 0
\put{$K_1(P C^*(\Lambda_N) P)$}[r] at 0 2          \arr(0.1,2)(1.4,2)            \put{$\ZZ^{c_N}$}[l] at 1.5 2
    \arr(-1.5,1.7)(-1.5,0.3)                                                          \arr(1.7,1.7)(1.7,0.3)
\put{$K_1(P C^*(\Lambda_{N+1}) P)$}[r] at 0 0        \arr(0.1,0)(1.4,0)            \put{$\ZZ^{c_{N+1}}$}[l] at 1.5 0
                        \put{$\phi^N_1$} [b] at 0.75 2.1
\put{$(i_N)_*$} [r] at  -1.6 1              \put{$B_N$} [r] at 1.6 1
                        \put{$\phi^{N+1}_1$} [t] at 0.75 -0.1
\endpicture
\)}} at 0 0
\endpicture
\end{equation}

Fix $N \in \NN$ and $j \le c_N$.

As in Lemma~\ref{lem:direct sum K-th} we write $X^N := V_0 (\BLUE\Lambda) V_N$ and decompose $X^N = \bigsqcup^{c_N}_{j = 1} X^N_j$. The inclusion $i_N$ of $C^*(\Lambda_N)$ in $C^*(\Lambda_{N+1})$ carries $s_\alpha s^*_\alpha$ into $\sum_{e \in s(\alpha) \Lambda^{e_1}} s_{\alpha e} s^*_{\alpha e}$. For each $i \le c_{N+1}$, exactly $A_N(i,j)$ of the paths $\alpha e$ lie in $X_i^{N+1}$, and hence the class of $s_\alpha s^*_\alpha$ in $K_0(P C^*(\Lambda_{N+1}) P)$ is given by
\[\textstyle
[s_\alpha s^*_\alpha] 
= \Big[ \sum_{e \in s(\alpha)\Lambda^{e_1}} s_{\alpha e} s^*_{\alpha e} \Big] 
= \sum_{e \in s(\alpha)\Lambda^{e_1}}\big[ s_{\alpha e} s^*_{\alpha e} \big]
= \sum^{c_{N+1}}_{i = 1} A_N(i,j) e_i.
\]
This establishes the commuting diagram on the left of~\eqref{eq:commuting diagrams}.

Now for the diagram on the right of~\eqref{eq:commuting diagrams}. For $1 \le j \le c_N$, let $e_j$ denote the generator of the $j^{\rm th}$ copy of $\ZZ$ in $K_1(P C^*(\Lambda_N) P)$. We set $M := \operatorname{lcm}\{A_N(i,j) : 1 \le i \le c_{N+1}, A_N(i,j) \not= 0\}$, and compute the image of $M e_j$ in $K_1(P C^*(\Lambda_{N+1})P)$. Let $\alpha \in X^N_j$. By Lemma~\ref{lem:direct sum K-th},
\begin{equation}\label{eq:M times generator}\textstyle
M e_j = \Big[ s_\alpha s_{\lambda(\alpha)^M} s^*_\alpha  + \sum_{\beta \in X^N\setminus\{\alpha\}} s_\beta s^*_\beta \Big].
\end{equation}
The effect of multiplying by $M$ is that if $e \in s(\alpha) \Lambda^{e_1} V_{N+1, i}$, then the path $\lambda(\alpha)^M e$ factors as $\sigma_i(e) \lambda(\alpha  e)^{M_i}$ where the integer $M_i$ is related to $M$ by 
\begin{equation}\label{eq:Ms related}
M |V_{N,j}| = M_i |V_{N+1, i}|,
\end{equation}
and where $\sigma_i$ is a permutation of $s(\alpha) \Lambda^{e_1} V_{N+1, i}$ which preserves the source map. The inclusion $i_N$ of $C^*(\Lambda_N)$ in $C^*(\Lambda_{N+1})$ carries the right-hand side of~\eqref{eq:M times generator} to the class of 
\begin{align*}
S &:= \sum_{e \in s(\alpha)\Lambda^{e_1}} s_\alpha s_{\lambda(\alpha)^M} s_e s^*_{\alpha e} + \sum_{\substack{\beta \in X^N\setminus\{\alpha\}\\ f \in s(\beta)\Lambda^{e_1}}} s_{\beta f} s^*_{\beta f} \\
&=  \Big(\sum^{c_{N+1}}_{i=1} \sum_{e \in s(\alpha)\Lambda^{e_1} V_{N+1, i}} s_{\alpha \sigma(e)} s_{\lambda(\alpha e)^{M_i}} s^*_{\alpha e}\Big) + \Big(\sum_{\substack{\beta \in X^N\setminus\{\alpha\}\\ f \in s(\beta)\Lambda^{e_1}}} s_{\beta f} s^*_{\beta f}\Big).
\end{align*}
To express $S$ in terms of our generators for $K_1(P C^*(\Lambda_{N+1}) P)$, let
\[
U := \sum_{\substack{1 \le i \le c_{N+1} \\ e \in s(\alpha)\Lambda^{e_1} V_{N+1, i}}} s_{\alpha e} s^*_{\alpha \sigma(e)} + \sum_{\substack{\beta \in X^N\setminus\{\alpha\}\\ f \in s(\beta)\Lambda^{e_1}}} s_{\beta f} s^*_{\beta f}.
\]
Then $U$ is unitary because the $\sigma_i$ are permutations. Moreover,
\[
US = \Big(\sum^{c_{N+1}}_{i=1} \sum_{e \in s(\alpha)\Lambda^{e_1} V_{N+1, i}} s_{\alpha e} s_{\lambda(\alpha e)^{M_i}} s^*_{\alpha e}\Big) + \Big(\sum_{\substack{\beta \in X^N\setminus\{\alpha\}\\ f \in s(\beta)\Lambda^{e_1}}} s_{\beta f} s^*_{\beta f}\Big).
\]
For any choice of distinguished edges in the isolated cycles at level $N+1$, the associated isomorphism of $P C^*(\Lambda_{N+1}) P$ onto $\bigoplus^{c_{N+1}}_{i=1} M_{X^{N+1}_i}(C(\TT))$ obtained from Propositions \ref{prp:M_P(C(T))}~and~\ref{prp:oplus} takes $U$ to a constant function, and hence $[U]$ is the identity element of $K_1(P C^*(\Lambda) P)$. Thus $(i_n)_*(M e_j) = [S] = [U] + [S] = [US]$. Moreover,
\begin{align*}
[US]
&= \Bigg[\Big(\sum^{c_{N+1}}_{i=1} \sum_{e \in s(\alpha)\Lambda^{e_1} V_{N+1, i}} s_{\alpha e} s_{\lambda(\alpha e)^{M_i}} s^*_{\alpha e}\Big) + \Big(\sum_{\substack{\beta \in X^N\setminus\{\alpha\}\\ f \in s(\beta)\Lambda^{e_1}}} s_{\beta f} s^*_{\beta f}\Big) \Bigg]\\
&= \Bigg[\prod_{i=1}^{c_{N+1}} \prod_{e \in s(\alpha)\Lambda^{e_1} V_{N+1, i}}\Big(s_{\alpha e} s_{\lambda(\alpha e)^{M_i}} s^*_{\alpha e} + \sum_{\beta f \in X^{N+1}\setminus\{\alpha e\}} s_{\beta f} s^*_{\beta f}\Big)\Bigg] \\
&= \sum_{i=1}^{c_{N+1}} \sum_{e \in s(\alpha)\Lambda^{e_1} V_{N+1, i}} M_i e_i \quad\text{by~\eqref{eq:K1 image}}\\
&= \sum_{i=1}^{c_{N+1}} |s(\alpha) \Lambda^{e_1} V_{N+1, i}|M_i e_i.
\end{align*}
Hence 
\begin{equation}\label{eq:K_1 expression}
(i_N)_*(e_j) = \frac{1}{M} \sum_{i=1}^{c_{N+1}} |s(\alpha) \Lambda^{e_1} V_{N+1, i}|M_i e_i = \sum_{i=1}^{c_{N+1}} A_N(i,j) \frac{M_i}{M} e_i.
\end{equation}
Recall from~\eqref{eq:Ms related} that for each $i$, the quantities $M$ and $M_i$ satisfy $M |V_{N,j}| = M_i |V_{N+1,i}|$. By Lemma~\ref{lem:A,B}(3), we also have $A_N(i,j) |V_{N,j}| = B_N(i,j) |V_{N+1, i}|$, so  
\[
\frac{M_i}{M} = \frac{B_N(i,j)}{A_N(i,j)}.
\]
Substituting this into~\eqref{eq:K_1 expression} gives $(i_N)_*(e_j) = \sum_{i=1}^{c_{N+1}} B_N(i,j) e_i$, and this establishes the commuting diagram on the right of~\eqref{eq:commuting diagrams}.
\end{proof}

To conclude the proof of Theorem~\ref{thm:AT and K-th}, it remains to prove assertion~(4). The idea is as follows. We construct a $1$-graph $B$ such that $C^*(B)$ is AF and $K_0(C^*(B))$ is canonically isomorphic to $K_0(C^*(\Lambda))$. We use the classifications of ideals in the $C^*$-algebras of graphs which satisfy condition~(K) \cite[Theorem~6.6]{KPRR} and the classification of gauge-invariant ideals in the $C^*$-algebras of $k$-graphs \cite[Theorem~5.2]{RSY1} to establish a lattice isomorphism between the ideals of $C^*(B)$ and the gauge-invariant ideals of $C^*(\Lambda)$. Finally, we use \cite[Theorem~1.5.3]{Ror} to obtain an isomorphism from the lattice of ideals of $C^*(B)$ to the lattice of order-ideals of $K_0(C^*(B))$.

The next result amounts to a restatement of results of Bratteli \cite{Bra} and Elliott \cite{Ell1} in the language of $1$-graph algebras. We give the result and the proof here for two reasons: firstly, the language and notation of $1$-graph algebras is more convenient to our later arguments than the traditional notation of Bratteli diagrams; and secondly, we want to establish explicit formulas linking the $K_0$-group and ideal structure of $C^*(\Lambda)$ for a \gbd/ $\Lambda$ to the $K_0$-group and ideal structure of an associated AF graph algebra.

\begin{prop}\label{prp:associated BD}
Let $\Lambda$ be a \gbd/ of infinite depth, and let $c_n$, $A_n$, $B_n$, $T_n$ be the data associated to $\Lambda$. Let $B$ be the $1$-graph with vertices $B^0 = \bigsqcup^\infty_{n=0} W_n$ where $W_n = \{w_{n,1}, \dots, w_{n, c_n}\}$ and with $A_n(i,j)$ edges from $w_{n+1,i}$ to $w_{n,j}$ for all $n, i, j$. Let $\{t_\beta : \beta \in B\}$ be the universal generating Cuntz-Krieger $B$-family in $C^*(B)$, and let $Q := \sum_{w \in W_0} t_w$. Then 
\begin{itemize}
\item[(1)] $Q$ is a full projection in $C^*(B)$; 
\item[(2)] for $n \in \NN$, the set $F_n = \lsp\{s_\alpha s^*_\beta : \alpha,\beta \in W_0 B W_n\}$ is a subalgebra of $Q C^*(B) Q$ and is canonically isomorphic to $\bigoplus_{j=1}^{c_n} M_{W_0 B w_{n,j}}(\CC)$;
\item[(3)] $F_n \subset F_{n+1}$ for all $n$, and $Q C^*(B) Q$ is equal to $\overline{\bigcup^\infty_{n=1} F_n}$ and hence is a unital AF algebra;
\item[(4)] there is an isomorphism $\phi : K_0(P C^*(\Lambda) P) \to K_0(Q C^*(B) Q)$ which satisfies $\phi([s_\eta s^*_\eta]) = [t_\beta t^*_\beta]$ for all $\eta \in V_0 (\BLUE\Lambda) V_{n,j}$ and $\beta \in W_0 B w_{n,j}$; and
\item[(5)] there is a lattice isomorphism between the ideals of $Q C^*(B) Q$ and the gauge-invariant ideals of $P C^*(\Lambda) P$ which takes $J \lhd Q C^*(B) Q$ to the ideal generated by $\{s_\eta s^*_\eta : \eta \in V_0 (\BLUE\Lambda) V_{n,j}, s_\beta s^*_\beta \in J \text{ for } \beta \in W_0 B w_{n,j}\}$ in $P C^*(\Lambda) P$.
\end{itemize}
\end{prop}
\begin{proof}
An argument more or less identical to the proof of \cite[Proposition~2.12]{CBMSbk} establishes claims~(1), (2)~and~(3)

(4) For $\beta \in W_0 B w_{n,j}$, $t_\beta t^*_\beta$ is a minimal projection in the $j^{\rm th}$ summand of $F_n$ and hence its class in $K_0(F_n)$ is the $j^{\rm th}$ generator $e_j$ of $\ZZ^{c_n}$. The inclusion map $\iota : F_n \to F_{n+1}$ takes a minimal projection $t_\beta t^*_\beta$ in $F_n$ to $\sum_{e \in s(\beta)B^1} t_{\beta e} t^*_{\beta e} \in F_{n+1}$. Hence $t_\beta t^*_\beta$ embeds in the $i^{\rm th}$ summand of $F_{n+1}$ as a projection of rank $|s(\beta) B^1 w_{n+1, i}| = A_n(i,j)$. It follows that 
\begin{equation}\label{eq:K-th of C*(B)}
K_0(Q C^*(B) Q) = \varinjlim(K_0(F_n), K_0(\iota)) = \varinjlim(\ZZ^{c_n}, A_n) 
\end{equation}
and this is equal to $K_0(P C^*(\Lambda) P)$ by Theorem~\ref{thm:AT and K-th}(2).
 
Equation~\eqref{eq:K-th of C*(B)} shows that for $\beta \in W_0 B w_{n,j}$, the class of $t_\beta t^*_\beta$ in $K_0(Q C^*(B) Q)$ is $A_{\infty, n}(e_j)$. But this is precisely the class of $s_\eta s^*_\eta$ in $K_0(P C^*(\Lambda) P)$ for any $\eta \in V_0 (\BLUE\Lambda) V_{n,j}$ by Lemma~\ref{lem:direct sum K-th}(2) and the left-hand commuting diagram of equation~\eqref{eq:commuting diagrams}.

(5) Since $B$ has no cycles, it satisfies condition~(K) of \cite{KPRR}. Hence \cite[Theorem~6.6]{KPRR} implies that the lattice of ideals of $C^*(B)$ is isomorphic to the lattice of saturated hereditary subsets of $B^0$ via $I \mapsto H_I := \{v : s_v \in I\}$. We have $I = \clsp\{t_\alpha t^*_\beta : s(\alpha) = s(\beta) \in H_I\}$. Since~(1) shows that $Q$ is full, the map $I \mapsto Q I Q$ is a lattice-isomorphism between ideals of $C^*(B)$ and ideals of $Q C^*(B) Q$. Hence if $J$ is an ideal in $Q C^*(B) Q$, we can sensibly define $H_J := H_I$ where $J = Q I Q$, and we have $J = \clsp\{t_\alpha t^*_\beta : \alpha,\beta \in W_0 B H_J\}$.

A similar analysis, using \cite[Theorem~5.2]{RSY1} instead of \cite[Theorem~6.6]{KPRR} shows that the lattice of gauge-invariant ideals of $P C^*(\Lambda) P$ is isomorphic to the lattice of saturated hereditary subsets of $\Lambda^0$ via $J \mapsto \{v : s_\tau s^*_\tau \in J\text{ for }\tau \in V_0 \Lambda v\}$. For $\tau \in \Lambda$ we can decompose $\tau = \eta\mu$ where $\eta \in \BLUE\Lambda$ and $\mu \in \RED\Lambda$, and we have $s_\tau s^*_\tau = s_\eta s^*_\eta$. If $H$ is hereditary, then $s(\tau) \in H$ if and only if $s(\eta) \in H$ because $\Lambda$ satisfies condition~\eqref{ass:AT}. Hence $J \mapsto \{v : s_\eta s^*_\eta \in J \text{ for } \eta \in V_0 (\BLUE\Lambda) v\}$ is a lattice-isomorphism between gauge-invariant ideals of $PC^*(\Lambda)P$ and saturated hereditary subsets of $\Lambda^0$.

Now if $H \subset \Lambda^0$ is hereditary, then for $n \in \NN$ and $1 \le j \le c_n$, either $V_{n,j} \subset H$ or $V_{n,j} \cap H = \emptyset$. It is easy to check using this that there is a bijection between the saturated hereditary subsets of $B^0$ and those of $\Lambda^0$ characterised by $H \subset B^0 \mapsto \bigcup\{V_{n,j} : w_{n,j} \in H\}$, and this completes the proof.
\end{proof}

\begin{proof}[Proof of Theorem~\ref{thm:AT and K-th}(4)]
Let $B$ be the 1-graph of Proposition~\ref{prp:associated BD}. Proposition~\ref{prp:associated BD}(3) shows that $Q C^*(B) Q$ is AF, so it is stably finite and has real-rank zero \cite[p~23]{Ror}. By \cite[Theorem~1.5.3]{Ror}, the map $J \mapsto K_0(J)$ is therefore an isomorphism from the ideal lattice of $Q C^*(B) Q$ to the order-ideal lattice of $K_0(Q C^*(B) Q)$.

The image of an ideal $J$ in $K_0(Q C^*(B) Q)$ is equal to $\varinjlim(K_0(J \cap F_n), A_n|_{K_0(J \cap F_n)})$. Since the ideals of $F_n$ are precisely its direct summands, $J \cap F_n$ is a direct sum of some subset of the direct summands of $F_n$, and so $K_0(J \cap F_n) = \langle [t_\beta t^*_\beta] : t_\beta t^*_\beta \in J \cap F_n\rangle$. Hence $K_0(J) = \langle [t_\beta t^*_\beta] : t_\beta t^*_\beta \in J\rangle \subset K_0(Q C^*(B) Q)$. By Proposition~\ref{prp:associated BD}(4), it follows that the image $\phi^{-1}(K_0(J))$ of $K_0(J)$ in $K_0(P C^*(\Lambda)P)$ is equal to 
\[
\langle [s_\eta s^*_\eta] : \eta \in V_0 (\BLUE\Lambda) V_{n,j}, s_\beta s^*_\beta \in J \text{ for } \beta \in W_0 B w_{n,j} \rangle.
\]
Proposition~\ref{prp:associated BD}(5) now establishes that there is an isomorphism $\theta$ from the lattice of gauge-invariant ideals of $P C^*(\Lambda) P$ to the lattice of order-ideals of $K_0(P C^*(\Lambda) P)$ which takes $J$ to $\langle [s_\eta s^*_\eta] : s_\eta s^*_\eta \in J\rangle$. 

Since $P$ is full, compression by $P$ induces an isomorphism $\phi_P$ of $K_0(C^*(\Lambda))$ onto $K_0(P C^*(\Lambda) P)$. For $\eta \in \BLUE\Lambda$, we have $s_\eta s^*_\eta \sim s_{s(\eta)}$, so $[s_{s(\eta)}] = [s_\eta s^*_\eta] \in K_0(C^*(\Lambda))$. Since each $s_\eta s^*_\eta \le P$ it follows that $\phi([s_{s(\eta)}]) = [s_\eta s^*_\eta] \in K_0(P C^*(\Lambda) P)$. Hence 
\[
\theta(J) = \phi_P(\langle [s_v] : v = s(\eta)\text{ for some }\eta\text{ with }s_\eta s^*_\eta \in J\rangle)
\]
for each ideal $J \in P C^*(\Lambda) P$. For an ideal $I$ of $C^*(\Lambda)$, $s_\eta s^*_\eta \in P I P$ if and only if $s_v \in I$. Since $P$ is full and gauge-invariant, $I \mapsto P I P$ is an isomorphism between the lattice of gauge-invariant ideals of $C^*(\Lambda)$ and that of $P C^*(\Lambda) P$. Thus $I \mapsto \phi_P^{-1}(\theta(P I P))$ is the desired lattice-isomorphism.
\end{proof}

\subsection*{Order units and dimension range}

Given a $C^*$-algebra $A$, we write $\mathcal{D}_0(A)$ for the \emph{dimension range}
\[
\mathcal{D}_0(A) = \{[p]_0 : p \in A\text{ is a projection}\} \subset K_0(A).
\]
Elliott's classification theorem implies that if $A$ is a simple A$\TT$ algebra with real-rank zero, then $A$ is determined up to isomorphism by the data $(K_0(A), K_1(A), [1_A]_0)$ if $A$ is unital, and by the data $(K_0(A), K_1(A), \mathcal{D}_0(A))$ if $A$ is non-unital (see \cite[Theorem~3.2.6]{Ror} and the subsequent discussion). In Section~\ref{sec:LPF} we identify conditions on a \gbd/ $\Lambda$ which ensure that $C^*(\Lambda)$ (and hence $P C^*(\Lambda) P$) is simple and has real-rank zero, so it is worth identifying the class $[P] \in K_0(P C^*(\Lambda) P)$ and the dimension range $\mathcal{D}_0(C^*(\Lambda)) \subset K_0(C^*(\Lambda))$.

\begin{lemma} \label{cor:position of unit}
Let $\Lambda$ be a \gbd/ of infinite depth. 
\begin{itemize}
\item[(1)] There is an order-isomorphism of $K_0(P C^*(\Lambda) P)$ onto $\varinjlim(\ZZ^{c_n}, A_n)$ which takes $[P]$ to the image of $(|V_{0,1}|, \dots, |V_{0, c_0}|) \in \ZZ^{c_0}$, and an isomorphism of $K_1(P C^*(\Lambda) P)$ onto $\varinjlim(\ZZ^{c_n}, B_n)$.
\item[(2)]  For $n \in \NN$ and $1 \le j \le c_n$ let $Y_{n,j} := (\BLUE\Lambda)V_{n,j}$, and let $D_n$ denote the subset $\{m \in \ZZ^{c_n} : 0 \le m_j \le |Y_{n,j}|\text{ for each }1 \le j \le c_n\}$. The isomorphism of $K_0(C^*(\Lambda))$ onto $\varinjlim(\ZZ^{c_n}, A_n)$ described in Theorem~\ref{thm:AT and K-th}(2) takes $\mathcal{D}_0(C^*(\Lambda))$ to the subset $\bigcup_{n = 0}^\infty A_{n,\infty}(D_n)$.
\end{itemize}
\end{lemma} 
\begin{proof}
The first statement follows from Lemmas \ref{lem:direct sum K-th}~and~\ref{lem:direct limits}. The second statement follows from a similar argument using Corollary~\ref{cor:oplus M(C(T))} to see that the isomorphism of $C^*(\Lambda_n)$ onto $\bigoplus^{c_n}_{j=1} M_{Y_{n,j}}(C(\TT))$ takes the class of the identity in $K_0(C^*(\Lambda_n))$ to $(|Y_{n,1}|, \dots, |Y_{n, c_n}|) \in \ZZ^{c_n}$ for each $n$.
\end{proof}

\begin{rmk}
For any choice of vertices $v_j \in V_{0,j}$, it is easy to check that the projection $p = \sum_{j=1}^{c_0} s_{v_j}$ is full in $C^*(\Lambda)$, and we can use Theorem~\ref{thm:AT and K-th} and Lemmas \ref{lem:direct sum K-th}~and~\ref{lem:direct limits} to see that there is an order-isomorphism of $K_0(p C^*(\Lambda) p)$ onto $\varinjlim(\ZZ^{c_n}, A_n)$ which takes $[p]$ to the usual order-unit $A_{\infty, 0}(1, \dots, 1)$.
\end{rmk}

\section{Large-permutation factorisations, simplicity, and real rank zero} \label{sec:LPF}

In this section we characterise the \gbds/ $\Lambda$ whose $C^*$-algebras are simple, and describe a condition on $\Lambda$ which ensures that $C^*(\Lambda)$ also has real-rank zero. Elliott's classification theorem for A$\TT$ algebras (see Theorem~3.2.6 and the discussion that follows it in \cite{Ror}) then implies that $C^*(\Lambda)$ is determined up to isomorphism by its $K$-theory.

In a \gbd/ the factorisation property induces a permutation $\Ff$ of the set $\BLUE\Lambda$ of blue paths: for $\alpha \in \BLUE\Lambda$, let $f$ be the unique red edge with $s(f) = r(\alpha)$, and define $\Ff(\alpha)$ to be the unique blue path such that $f \alpha = \Ff(\alpha) f'$ for some red edge $f'$. (In the notation of \S3, $\Ff(\alpha) = \Ff_1(f,\alpha)$.) For $\alpha \in \BLUE\Lambda$, the \emph{order} $o(\alpha)$ of $\alpha$ is the smallest $k > 0$ such that $\Ff^k(\alpha) = \alpha$. If $r(\alpha) \in V_{n,j}$ and $\mu$ is the unique red path with $s(\mu) = r(\alpha)$ and $|\mu| = o(\alpha)$, then $\mu\alpha$ has the form $\alpha\mu'$, and $\mu = \lambda(r(\alpha))^m$ for $m = o(\alpha)/|V_{n,j}|$.

Recall from \cite[Definition~4.7]{KP} that a $k$-graph $\Lambda$ is \emph{cofinal} if for every vertex $v$ and every infinite path $x$ there exists $n \in \NN^2$ such that $v \Lambda x(n)$ is nonempty.

\begin{theorem}\label{critsimple}
Let $\Lambda$ be a \gbd/. Then $C^*(\Lambda)$ is simple if and only if $\Lambda$ is cofinal and $\{o(e) : e \in \Lambda^{e_1}\}$ is unbounded.
\end{theorem}

To prove the theorem, we first need to establish some properties of the order function~$o$.

\begin{lemma}\label{proposofo}
Let $\Lambda$ be a \gbd/.
\begin{enumerate}
\item Suppose that $\alpha=\mu g\nu$ and $\alpha=\beta h\gamma$ are two factorisations of $\alpha\in \Lambda$ in which $d(\mu)$ and $d(\beta)$ have the same $1^{\rm st}$ coordinates and $g,h \in \Lambda^{e_1}$. Then $o(g)=o(h)$.
\item For every blue path $\beta$ of length $n$, $o(\beta)=\lcm(o(\beta_1),\dots,o(\beta_n))$.
\end{enumerate}
\end{lemma}
\begin{proof}
For (1), write $d(\mu)=(n,k)$ and $d(\beta)=(n,l)$, and without loss of generaility suppose $l\geq k$.  Then $\mu(0,(n,k))=\alpha(0,(n,k))=\delta$, say, and the factorisation property implies that $\beta=\delta\beta'$. Since $d(\beta')=(l-k)e_2$, $\beta'$ is the unique red path of length $l-k$ from $r(h)$ to $s(\delta)=r(g)$. Thus $g$ is the image $\Ff^{l-k}(h)$ of $h$ under the $(l-k)^{\rm th}$ iteration of the permutation $\Ff$, and hence has the same order as $h$. (This is a general property of permutations of sets.) 

For (2), notice that the uniqueness of factorisations implies that 
\[
\Ff^k(\beta)=\Ff^k(\beta_1)\Ff^k(\beta_2)\cdots\Ff^k(\beta_n)
\]
is equal to $\beta$ if and only if $\Ff^k(\beta_i)=\beta_i$ for all $i$.
\end{proof}

\begin{cor}\label{cor:preserves order}
Let $\Lambda$ be a \gbd/, and suppose that $\mu \alpha = \alpha' \mu'$ where $\mu, \mu'$ are red and $\alpha, \alpha'$ are blue. Then $o(\alpha) = o(\alpha')$.
\end{cor}

We aim to prove simplicity of $C^*(\Lambda)$ by verifying that $\Lambda$ satisfies the aperiodicity Condition (A) of \cite{KP}, so we begin by recalling some definitions from \cite{KP}. We denote by $\Omega_k$ the $k$-graph with vertices $\Omega_k^0 := \NN^k$, paths $\Omega_k^m = \{(n, n+m) : n \in \NN^k\}$ for $m \in \NN^k$, $r((n, n + m)) = n$ and $s((n, n+m)) = n+m$. The \emph{infinite paths} in a $k$-graph $\Lambda$ with no sources are the degree preserving functors $x : \Omega_k \to \Lambda$. The collection of all infinite paths of $\Lambda$ is denoted $\Lambda^\infty$, and the range of $x$ is the vertex $x(0)$. For $p \in \NN^k$ and $x \in \Lambda^\infty$, $\sigma^p(x) \in \Lambda^\infty$ is defined by $\sigma^p(x)(n) := x(n+p)$ \cite[Definitions~2.1]{KP}. A path $x\in\Lambda^\infty$ is \emph{aperiodic} if $\sigma^p(x)=\sigma^q(x)$ implies $p=q$.

The next lemma will help us recognise aperiodic paths.

\begin{lemma}\label{critaper}
Suppose that $x$ is an infinite path in a \gbd/ $\Lambda$ such that $o(x(0,ne_1))\to \infty$ as $n\to \infty$. Then $x$ is aperiodic.
\end{lemma} 

\begin{proof}
Suppose that $p,q \in \NN^2$ satisfy $\sigma^p(x) = \sigma^q(x)$. We must show that $p = q$. If $n$ is the integer such that $r(x)\in V_n$, then $r(\sigma^p(x)) \in V_{n+p_1}$ and $r(\sigma^q(x)) \in V_{n+ q_1}$. Thus $\sigma^p(x) = \sigma^q(x)$ implies that $p_1 = q_1$. Without loss of generality, we may suppose $q_2 \ge p_2$. Now the infinite path $y := \sigma^{(p_1, p_2)}(x) = \sigma^{(q_1, p_2)}(x)$ satisfies $\sigma^{le_2}(y) = y$ where $l := q_2 - p_2$, so $y=y(0,le_2)y$. Since $o(x(0,p_1e_1))$ is finite, Corollary~\ref{cor:preserves order} implies that the path $y$ also satisfies $o(y(0, ne_1)) \to \infty$ as $n \to \infty$. But for every $n$ we have
\[
y(0,ne_1)=(y(0,le_2)y)(0,ne_1)=\Ff^l(y(0,ne_1)),
\]
and hence we must have $l=0$, $p_2=q_2$ and $p=q$.
\end{proof} 

\begin{proof}[Proof of sufficiency in Theorem~\ref{critsimple}]
Suppose that $\Lambda$ is cofinal and $\{o(e):e\in\Lambda^{e_1}\}$ is unbounded. To show that $C^*(\Lambda)$ is simple, it suffices by \cite[Proposition~4.8]{KP} to show that for each $w\in \Lambda^0$ there is an aperiodic path $x$ with $r(x)=w$. 

We first claim that for every $v\in \Lambda^0$ there exists $N$ such that $v\Lambda V_{N,i}$ is nonempty for every $i\leq c_N$. To prove this, we suppose to the contrary that there exists $v \in \Lambda^0$, say $v \in V_n$, and a sequence $\{i_m : m > n\}$ such that $v \Lambda V_{m,i_m} = \emptyset$ for all $m$. By assumption the sinks in $\Lambda$ belong to $V_0$, so for each $m > n$ there exists a path $\xi_m \in V_n \Lambda V_{m, i_m}$. Let $p_0 := \{\xi_m : m > n\}$. Since $\Lambda$ is row-finite, there exists $g_1 \in V_n\Lambda^{e_1}$ such that $p_1 := \{\eta \in p_0 : \eta(0, e_1) = g_1\}$ is infinite. For the same reason, there then exists $g_2 \in s(g_1)\Lambda^{e_1}$ such that $p_2 := \{\eta \in p_1 : \eta(e_1, 2e_1) = g_2\}$ is infinite. Continuing in this way, we obtain a sequence $g_i$ of blue edges such that for each $i$ there are infinitely many $m$ with $\xi_m(0, ie_1) = g_1 \dots g_i$. By choice of the $\xi_m$, we then have $v \Lambda s(g_i) = \emptyset$ for all $i$. For each $i$, let $x_i := g_1\lambda(s(g_1))g_2\lambda(s(g_2))\dots g_i\lambda(s(g_i))$. By \cite[Remark~2.2]{KP}, there is a unique infinite path $x \in \Lambda^\infty$ such that $x(0, d(x_i)) = x_i$ for all $i$. By construction, we have $v \Lambda x(n) = \emptyset$ for all $n \in \NN^2$. This contradicts the cofinality of $\Lambda$, and we have justified the claim.

We now fix $w\in \Lambda^0$, and construct an aperiodic path with range $v$. By the claim there exists $N \in \NN$ such that $v\Lambda V_{N, j} \not= \emptyset$ for all $j \leq c_{N}$. Since the sinks in $\Lambda^0$ belong to $V_0$, we then have $w\Lambda V_{M, i} \not= \emptyset$ for all $M\ge N$ and $i \le c_M$. Since $\sup\{o(e) : e \in \Lambda^{e_1}\} = \infty$, and since $\bigcup^{N - 1}_{n = 0} V_n \Lambda^{e_1}$ is finite, there exists $M \ge N$ and $g \in V_{M} \Lambda^{e_1}$ such that $o(g) \ge 2$. By choice of $g$ there exists a path $\alpha_2 \in v \Lambda g$, and we may assume that $d(\alpha_2) \ge (1,1)$. Repeating this procedure at the vertex $v=s(\alpha_2)$ gives a path $\alpha_3$ in $s(\alpha_2)\Lambda$ such that $o(\alpha_3(d(\alpha_3) - e_1, d(\alpha_3))) \ge 3$ and $d(\alpha_3) \ge (1,1)$. By continuing this way we can inductively construct a sequence of paths $\alpha_i$ with $s(\alpha_i) = r(\alpha_{i+1})$, $d(\alpha_i) \ge (1,1)$ and $o(\alpha_i(d(\alpha_i) - e_1, d(\alpha_i))) \ge i$. By \cite[Remark~2.2]{KP}, there is a unique infinite path $x$ such that $x(0, d(\alpha_2) + \cdots + d(\alpha_i)) = \alpha_2 \dots \alpha_i$ for all $i$. Part (1) of Lemma~\ref{proposofo} implies that $d(x(0,ne_1))\geq o(\alpha_i(d(\alpha_i) - e_1, d(\alpha_i)))\geq i$ for sufficiently large $n$, and hence $d(x(0,ne_1))\to \infty$ as $n\to \infty$. Thus it follows from Lemma~\ref{critaper} that $x$ is aperiodic, and since $r(x)=r(\alpha_2)=w$, this completes the proof.
\end{proof}

For the other direction in Theorem~\ref{critsimple}, we show that if $\{o(e)\}$ is bounded then the graph is periodic, and apply the following general result.

\begin{prop}\label{prp:non-gi ideal}
Let $\Lambda$ be a row-finite $k$-graph. Suppose that there is a vertex $v \in \Lambda^0$ and an element $p \in \NN^k$ such that every $x \in v\Lambda^\infty$ satisfies $x = \sigma^p(x)$. Then $C^*(\Lambda)$ is not simple.
\end{prop}
\begin{proof}
Let $\{S_\lambda : \lambda \in \Lambda\}$ be the Cuntz-Krieger family on $\ell^2(\Lambda^\infty)$ given by $S_\lambda e_x = \delta_{s(\lambda), r(x)} e_{\lambda x}$ \cite[Proposition~2.11]{KP}. Then $S^*_\lambda e_x = \delta_{\lambda, x(0,d(\lambda))} e_{\sigma^{d(\lambda)}(x)}$ for all $\lambda \in \Lambda$ and $x \in \Lambda^{\infty}$. Let $\pi_S$ be the corresponding representation of $C^*(\Lambda)$.

Fix $\mu \in v\Lambda^p$. By assumption on $v\Lambda^{\infty}$, 
\[
S^*_\mu e_x = \delta_{\mu, x(0,p)} e_{\sigma^p(x)} = \delta_{\mu, x(0,p)} e_x = S_\mu S^*_\mu e_x
\]
for all $x \in \Lambda^\infty$. Hence $S^*_\mu = S_\mu S^*_\mu$. Let $\gamma$ be the gauge action on $C^*(\Lambda)$. Fix $z \in \TT^k$ such that $\overline{z}^p = -1$. Then $\gamma_z(s^*_\mu) = -s^*_\mu$ and $\gamma_z(s_\mu s^*_\mu) = s_\mu s^*_\mu$, so $s_\mu s^*_\mu - s^*_\mu \not= 0$. However $\pi_S(s_\mu s^*_\mu - s^*_\mu) = S^*_\mu S_\mu - S^*_\mu = 0$, so the kernel of $\pi_S$ is a nontrivial ideal in $C^*(\Lambda)$.
\end{proof}

\begin{proof}[Proof of necessity in Theorem~\ref{critsimple}]
Suppose that $C^*(\Lambda)$ is simple. Let $B$ be the $1$-graph associated to $\Lambda$ as in Proposition~\ref{prp:associated BD}. By Proposition~\ref{prp:associated BD}(5), $Q C^*(B) Q$ is simple, so Proposition~\ref{prp:associated BD}(1) shows that $C^*(B)$ is simple. It now follows from \cite[Proposition~5.1]{BPRS} that $B$ is cofinal, and from the definition of $B$ that $\Lambda$ is also cofinal.

We now argue by contradiction that $\{o(e) : e \in \Lambda^{e_1}\}$ is unbounded. Suppose to the contrary that $o(e) \le l$ for all $e \in \Lambda^{e_1}$. Then $o(e)$ divides $l!$ for all $e \in \Lambda^{e_1}$. Let $p = l! e_2$. We claim that $\sigma^p(x) = x$ for every $x \in \Lambda^\infty$. To compute $\sigma^p(x)$, we first factor $x$ as
\begin{equation}\label{eq:factor x}
x = \mu g_1 \lambda(s(g_1)) g_2 \lambda(s(g_2)) \dots
\end{equation}
where $d(\mu) = p$ and $d(g_i) = e_1$. Since $o(g_1)$ divides $l!$ and $\mu$ is the unique path of length $l!$ starting at $r(g_1)$, $\mu g_1$ has the form $g_1 \mu_1$ where $d(\mu_1) = d(\mu) = p$. Since the cycle $\lambda(s(g_1))$ is isolated, we have $\mu_1\lambda(s(g_1)) g_2 = \lambda(s(g_1))\mu_1 g_2$. Now $\mu_1$ is the unique red path of length $l!$ starting at $r(g_2)$. Since $o(g_2)$ divides $l!$, $\lambda(s(g_1))\mu_1 g_2$ has the form $\lambda(s(g_1)) g_2 \mu_2$ where $d(\mu_2) = d(\mu_1) = p$. Continuing this way shows that we can also factor $x$ as
\[
x = g_1 \lambda(s(g_1)) g_2 \lambda(s(g_2)) \dots,
\]
and then~\eqref{eq:factor x} implies that $\sigma^p(x) = x$, establishing the claim. Proposition~\ref{prp:non-gi ideal} now implies that $C^*(\Lambda)$ is not simple, which is a contradiction.
\end{proof}

We now turn our attention to the problem of deciding when $C^*(\Lambda)$ has real-rank zero.

\begin{dfn}\label{dfn:LPF}
We say that a \gbd/ $\Lambda$ has \emph{large-permutation factorisations} if for each $v \in \Lambda^0$ and each integer $l > 0$ there exists $N \in \NN$ such that 
\begin{equation}\label{eq:alternate LPF}
\text{$o(\alpha) > l$ for all $\alpha \in v\Lambda^{N e_1}$.}
\end{equation}
\end{dfn}

Since \gbds/ are row-finite, a \gbd/ with large-permutation factorisations must have infinite depth, and Lemma~\ref{critaper} implies that every infinite path in $\Lambda$ is aperiodic.

There are several ways to ensure that a \gbd/ has large-permutation factorisations. For example, this is automatically the case if the red cycles get larger as $n$ grows, or more precisely if $\min_j |V_{n,j}| \to \infty$ as $n \to \infty$. Alternatively, we can keep $|V_{n,j}|$ small but add lots of blue edges entering each $V_{n,j}$ and define the factorisation property to ensure that $\min\{o(e):r(e)\in V_n\} \to \infty$ as $n\to \infty$.

\begin{theorem} \label{thm:LPF -> RR0}
Let $\Lambda$ be a \gbd/ with large-permutation factorisations.
\begin{itemize}
\item[(1)] Every ideal of $C^*(\Lambda)$ is gauge-invariant, and the lattice of ideals of $C^*(\Lambda)$ is isomorphic to the lattice of order-ideals of $K_0(C^*(\Lambda))$ via the map described in Theorem~\ref{thm:AT and K-th}(4).
\item[(2)] If $\Lambda$ is cofinal, then $C^*(\Lambda)$ is simple and $C^*(\Lambda)$ has real-rank zero.
\end{itemize}
\end{theorem}

To prove Theorem~\ref{thm:LPF -> RR0}(2), we show that the projections in $C^*(\Lambda)$ separate the tracial states with a view to applying \cite[Theorem~1.3]{BBEK}.

Recall that for $\alpha \in \Lambda$, $\lambda(\alpha)$ denotes the isolated cycle with range and source $s(\alpha)$. For the next result, we adopt the convention that for a negative integer $m$, $s_{\lambda(\alpha)^m} := s^*_{\lambda(\alpha)^{-m}}$.

\begin{lemma}\label{lem:traces diagonal}
Let $\Lambda$ be a \gbd/ and let $\tau$ be a trace on $C^*(\Lambda)$. Let $\alpha,\beta \in \BLUE\Lambda$ and $\mu \in \RED\Lambda$. Suppose that $\tau(s_\alpha s_\mu s^*_\beta) \not= 0$ or that $\tau(s_\alpha s^*_\mu s^*_\beta) \not= 0$. Then $\alpha = \beta$ and $\mu = \lambda(\alpha)^m$ for some $m \in \ZZ$. 
\end{lemma}
\begin{proof}
We argue the case where $\tau(s_\alpha s_\mu s^*_\beta) \not= 0$; the other case is similar. Since $\tau$ is a trace, we have $\tau(s_\beta^* s_\alpha s_\mu) = \tau(s_\alpha s_\mu s^*_\beta) \not= 0$. Since $\mu \in \RED\Lambda$, both $s(\alpha)$ and $s(\beta)$ belong to the same level of the \gbd/ $\Lambda$; say $s(\alpha), s(\beta) \in V_n$. Furthermore, since $s^*_\beta s_\alpha s_\mu \not= 0$, the ranges of $\alpha$ and $\beta$ must coincide, and in particular belong to the same $V_m$. Since $\Lambda$ is a \gbd/ it follows that $d(\alpha) = d(\beta) = (n-m)e_1$. This forces $\alpha = \beta$. But now $r(\mu) = s(\beta) = s(\alpha) = s(\mu)$, and it follows that $\mu = \lambda(\alpha)^m$ for some $m \in \NN$.
\end{proof}

\begin{lemma}\label{lem:LPF -> traces kill unitaries}
Let $\Lambda$ be a \gbd/. Suppose that $\Lambda$ has large-permutation factorisations, and let $\tau$ be a trace on $C^*(\Lambda)$. If $\alpha \in \BLUE\Lambda$ and $\tau(s_\alpha s_{\lambda(s(\alpha))^m} s^*_\alpha) \not = 0$, then $m = 0$.
\end{lemma}
\begin{proof}
We show that $m > 0$ implies that $\tau(s_\alpha s_{\lambda(s(\alpha))^m} s^*_\alpha) = 0$; a similar argument shows that $m<0$ is also impossible, so that $m$ must be $0$.

So suppose that $m > 0$. Taking $v = s(\alpha)$ and $l = m|\lambda(\alpha)|$ in Definition~\ref{dfn:LPF} we obtain an integer $n$ such that for all $\beta \in v \Lambda^{n e_1}$, $\Ff^{m|\lambda(v)|}(\beta) \not= \beta$. But now
\begin{align*}
\tau(s_\alpha s_{\lambda(v)^m} s^*_\alpha)
&= \tau\Big(\sum_{\beta \in v\Lambda^{n e_1}} s_\alpha s_{\lambda(v)^m} s_\beta s^*_\beta s^*_\alpha\Big) \quad\text{by~(CK4)} \\
&= \sum_{\beta \in v\Lambda^{n e_1}} \tau(s_{\alpha \Ff^{m|\lambda(v)|}(\beta)} s_{\mu(\beta)} s^*_{\alpha\beta})
\end{align*}
where for each $\beta$, $\mu(\beta)$ is the unique element of $s(\beta) \Lambda$ of degree $d(\lambda(v)^m)$. By choice of $n$, $\alpha \Ff^{m|\lambda(v)|}(\beta) \not = \alpha\beta$ for each term in the sum, and it follows from Lemma~\ref{lem:traces diagonal} that $\tau(s_\alpha s_{\lambda(v)^m} s^*_\alpha) = 0$.
\end{proof}

\begin{cor}\label{cor:projections separate traces}
Let $\Lambda$ be a \gbd/ with large-permutation factorisations. Then the projections in $C^*(\Lambda)$ separate traces on $C^*(\Lambda)$.
\end{cor}
\begin{proof}
Let $\tau_1$ and $\tau_2$ be traces on $C^*(\Lambda)$ which agree on all the projections in $C^*(\Lambda)$. Then $\tau_1(s_\alpha s^*_\alpha) = \tau_2(s_\alpha s^*_\alpha)$ for all $\alpha \in \BLUE\Lambda$. By Lemmas \ref{lem:traces diagonal}~and~\ref{lem:LPF -> traces kill unitaries}, it follows that $\tau_1$ and $\tau_2$ agree on $\clsp\{s_\alpha s_\mu s^*_\beta, s_\alpha s^*_\mu s_\beta : \alpha,\beta \in \BLUE\Lambda, \mu \in \RED\Lambda\}$, which by Lemma~\ref{lem:spanning} and the first assertion of Lemma~\ref{lem:direct limits} is all of $C^*(\Lambda)$. Thus $\tau_1 = \tau_2$.
\end{proof}

\begin{proof}[Proof of Theorem~\ref{thm:LPF -> RR0}]
If $H$ is a saturated hereditary subset of $\Lambda^0$, then $\Gamma_H := \Lambda \setminus \Lambda H = \{\eta \in \Lambda : s(\eta) \not\in H\}$ is a $k$-graph by \cite[Theorem~5.2(b)]{RSY1}. Theorem~5.3 of \cite{RSY1} implies that if $\Gamma_H$ satisfies \cite[Condition~(B)]{RSY1} for every saturated hereditary $H \subset \Lambda^0$, then every ideal of $C^*(\Lambda)$ is gauge-invariant. Remark~(4.4) of \cite{RSY1} shows that if $\Gamma_H$ has no sources and satisfies the aperiodicity condition \cite[Condition~(A)]{KP}, then $\Gamma_H$ satisfies \cite[Condition~(B)]{RSY1}. It therefore suffices to show that if $H \subset \Lambda^0$ is saturated and hereditary, then $\Gamma_H$ has no sources and satisfies the aperiodicity condition.

Fix a saturated hereditary subset $H$ of $\Lambda^0$. If $v \in \Gamma_H^0$, then $v\Lambda^{e_i}$ is nonempty. If $v \Lambda^{e_i} \subset \Lambda H$, then $v \in H$ because $H$ is saturated, contradicting $v \in \Gamma_H^0$. Thus there exists $e \in v\Lambda^{e_i} \setminus \Lambda H = v\Gamma_H^{e_i}$, and $\Gamma_H$ has no sources. Each infinite path of $\Gamma_H$ is also an infinite path of $\Lambda$, and hence is aperiodic by Lemma~\ref{critaper}. Thus $\Gamma_H$ satisfies the aperiodicity condition of \cite{KP}.

The rest of (1) now follows from Theorem~\ref{thm:AT and K-th}(4).

For (2), we first deduce from Lemma~\ref{proposofo} that $\{o(e) : e \in \Lambda^{e_1}\}$ has to be unbounded, and hence Theorem~\ref{critsimple} implies that $C^*(\Lambda)$ is simple. Corollary~\ref{cor:projections separate traces} implies that the projections in $C^*(\Lambda)$ separate the tracial states. Theorem~\ref{thm:AT and K-th} guarantees that $C^*(\Lambda)$ is an A$\TT$ algebra. By \cite[Theorem~1.3]{BBEK}, a simple A$\TT$ algebra $A$ has real-rank zero if and only if the projections of $A$ separate the tracial states, and this proves the result.
\end{proof}

\section{Achievability of classifiable algebras} \label{sec:Achievability}

In this section we characterise the $K$-group pairs which can arise as those of $P C^*(\Lambda) P$ when $\Lambda$ is a \gbd/ with large-permutation factorisations. We have already shown that the data associated to a \gbd/ $\Lambda$ consists of integers $c_n$, matrices $A_n, B_n \in M_{c_{n+1}, c_n}(\ZZ_+)$ with no zero rows or columns, and diagonal matrices $T_n \in M_{c_n}(\ZZ_+)$ with positive diagonal entries such that $K_0(C^*(\Lambda)) = \varinjlim(\ZZ^{c_n}, A_n)$, $K_1(C^*(\Lambda)) = \varinjlim(\ZZ^{c_n},B_n)$, and $A_n T_n = T_{n+1} B_n$ for all $n$. Here we establish a converse and characterise the $K$-groups that can arise when $PC^*(\Lambda)P$ is simple with real-rank zero.

\begin{dfn}
We say that an integer matrix $M$ is \emph{proper} if all entries of $M$ are nonnegative, and each row and each column of $M$ contains at least one nonzero entry (cf \cite[\S A4]{Eff}). Note that a diagonal matrix is proper if and only if all diagonal entries are nonzero. 
\end{dfn}

The data associated to a \gbd/ always consists of proper matrices.

\begin{theorem}\label{thm:classifiable achievability}
\begin{itemize}
\item[(1)] Let $\Lambda$ be a \gbd/ of infinite depth and suppose that $C^*(\Lambda)$ is simple. Then $K_0(C^*(\Lambda))$ is a simple dimension group which is not isomorphic to $\ZZ$.
\item[(2)] Let $\{c_n : n \in \NN\}$ be positive integers. For each $n$, let $A_n, B_n \in M_{c_{n+1}, c_n}(\ZZ_+)$ be proper matrices, and let $T_{c_n} \in M_n(\ZZ_+)$ be a proper diagonal matrix. Suppose additionally that $A_n T_n = T_{n+1} B_n$ for all $n$. Then there exists a \gbd/ $\Lambda$ such that $K_0(C^*(\Lambda)) \cong \varinjlim (\ZZ^{c_n}, A_n)$ and $K_1(C^*(\Lambda)) \cong \varinjlim (\ZZ^{c_n}, B_n)$. If $\varinjlim(\ZZ^{c_n}, A_n)$ is simple dimension group which is not isomorphic to $\ZZ$, then $\Lambda$ can be chosen so that $C^*(\Lambda)$ is simple with real-rank zero.
\end{itemize}
\end{theorem}

\begin{rmk}
In Theorem~\ref{thm:classifiable achievability}(2), we do not claim that there is a \gbd/ $\Lambda$ with data $c_n$, $A_n$, $B_n$, $T_n$. We \emph{can} always build a \gbd/ $\Lambda$ with the specified data (see Proposition~\ref{prp:build the gbd}), and if $\varinjlim(\ZZ^{c_n}, A_n)$ is a simple dimension group, then $\Lambda$ will be cofinal. However, to ensure that $C^*(\Lambda)$ is simple and has real-rank zero, we construct a \gbd/ with large-permutation inclusions, and to do this, we have to choose a subsequence of $\NN$ and adjust the data $c_n, A_n, B_n, T_n$ accordingly.
\end{rmk}

\begin{proof}[Proof of Theorem~\ref{thm:classifiable achievability}(1)]
We will show that if $K_0(P C^*(\Lambda) P)$ is isomorphic to $\ZZ$ or is not simple as a dimension group, then $C^*(\Lambda)$ is not simple. Let $B$ be the $1$-graph of Proposition~\ref{prp:associated BD}. Then parts~(3)~and~(4) of Proposition~\ref{prp:associated BD} imply that $Q C^*(B) Q$ is a unital AF algebra with $K_0(Q C^*(B) Q) = \varinjlim(\ZZ^{c_n}, A_n)$ isomorphic as a dimension grpoup to $K_0(P C^*(\Lambda) P)$. 

First suppose that $K_0(P C^*(\Lambda) P)$ is isomorphic to $\ZZ$. Then $K_0(Q C^*(B) Q)$ is order-isomorphic to $\ZZ$, and hence $Q C^*(B) Q \cong M_n(\CC)$ where $n = [1_A]$ is the class of the unit \cite[Theorem~7.3.4]{RLL}. Since $Q C^*(B) Q$ is finite-dimensional, the approximating subalgebras $F_n$ of Proposition~\ref{prp:associated BD}(2) must equal $Q C^*(B) Q$ for large $n$. Thus $F_n$ eventually has just one direct summand, so $c_n = |W_n| = 1$ for large $n$ by Proposition~\ref{prp:associated BD}(2). Moreover, since $F_n=F_{n+1}$ for large $n$, we must have $|W_0 B W_n| = |W_0 B W_{n+1}|$ for large $n$ by Proposition~\ref{prp:associated BD}(2), so $A_n(1,1) = |W_n B^1 W_{n+1}| = 1$ for large $n$, say for $n \ge M$.

So for $n \ge M$, $V_n \subset \Lambda^0$ consists of the vertices on a single red cycle, and each vertex in $V_n$ receives exactly one blue edge from $V_{n+1}$. Since $A_n(1,1) = 1$ for all $n \ge M$, Equation~\eqref{eq:A,B relation} implies that $(|V_n|)_{n \ge M}$ is a decreasing sequence of positive integers, and hence is eventually constant; say $|V_n| = c$ for $n \ge N \ge M$. The set $H = \bigcup^\infty_{n = N} V_n$ is hereditary in the sense of \cite[Section~5]{RSY1}, and its saturation $\overline{H}$ is all of $\Lambda^0$. Thus $C^*(\Lambda)$ is Morita equivalent to $C^*(H\Lambda)$ by \cite[Theorem~5.2]{RSY1}. Now $H\Lambda$ is isomorphic to the product graph $\Omega_1 \times C_c$, where $C_c$ is the red cycle with $c$ vertices, and $\Omega_1$ is the one-sided infinite path of blue edges (see \cite[Examples~1.7(ii)]{KP}). Corollary~3.5(iv) of \cite{KP} shows that $C^*(\Omega_1 \times C_c) \cong C^*(\Omega_1) \otimes C^*(C_c) \cong \Kk \otimes M_c(C(\TT))$, which is not simple.

Now suppose that $K_0(P C^*(\Lambda) P)$ is not simple as a dimension group. Then the AF algebra $Q C^*(B) Q$ is not simple either \cite[Corollary~1.5.4]{Ror}, and  Proposition~\ref{prp:associated BD}(5) implies that $C^*(\Lambda)$ is not simple.
\end{proof}

For the second claim of the theorem, we need to know how to build a \gbd/ from the data $c_n, A_n, B_n$ and $T_n$. Recall from Definition~\ref{dfn:LPF} that for $e \in \Lambda^{e_1}$, the order $o(e)$ of $e$ is the length of the shortest nontrivial path $\mu \in \RED\Lambda$ such that $(\mu e)(0, e_1) = e$.

\begin{prop}\label{prp:build the gbd}
Let $c_n$, $A_n$, $B_n$ and $T_n$ be as in Theorem~\ref{thm:classifiable achievability}(2). There is a \gbd/ $\Lambda$ with this data which has the following property: for each blue edge $e \in \Lambda^1$, say $r(e) \in V_{n,j}$ and $s(e) \in V_{n+1, i}$, we have $o(e) = A_n(i,j)|V_{n,j}|$.
\end{prop}
\begin{proof}
Since the data of~\eqref{eq:comm diagram} is all contained in the $1$-skeleton of the \gbd/, \cite[Section~6]{KP} shows that we need only construct a $1$-skeleton with the right number of edges, and an allowable collection of commuting squares so that the order of each blue edge in $V_{n,j} \Lambda^{e_1} V_{n+1, i}$ is $A_n(i,j) |V_{n,j}|$. By~\eqref{eq:A,B relation}, this is equivalent to showing that the order of each blue edge in $V_{n,j} \Lambda^{e_1} V_{n+1, i}$ is maximal.

For each $n$, the matrix $T_n$ defines a collection of $c_n$ isolated cycles $\lambda_{n,j}$ ($1 \le j \le c_n$) where $\lambda_{n,j}$ has $T_n(j,j)$ vertices. The collection of all paths in these cycles is $\RED\Lambda$, and the vertices on each $\lambda_{n,j}$ are the elements of $V_{n,j}$.

We want to show that for each $j \le c_n$ and $i \le c_{n+1}$, we can: 
\begin{itemize}
\item[(1)] add blue edges from vertices in $\lambda_{n+1,i}$ to vertices in $\lambda_{n,j}$ so that the number of blue edges to each vertex on $\lambda_{n,j}$ from $\lambda_{n+1,i}$ is $a := A_n(i,j)$ and the number of blue edges to $\lambda_{n,j}$ from each vertex on $\lambda_{n+1,i}$ is $b := B_n(i,j)$; and 
\item[(2)] specify an allowed collection of commuting squares so that the resulting permutation of the blue edges in $V_{n,j} \Lambda^{e_1} V_{n+1,i}$ is maximal. 
\end{itemize}
Let $v := T_n(j,j)$ be the number of vertices on $\lambda_{n,j}$, and let $w := T_{n+1}(i,i)$ be the number of vertices on $\lambda_{n+1,i}$. The commutativity of~\eqref{eq:comm diagram} says that $av = wb$.

We first demonstrate that it suffices to show how to add the desired blue edges when $a$ and $b$ have no common divisors; that is, when $(a,b) = 1$. To see this, suppose that $a = da'$ and $b = db'$, and that we can add the desired edges to obtain the data $A_n(i,j) = a'$ and $B_n(i,j) = b'$. Then we take the resulting diagram, and add $d-1$ blue edges $e(1), \dots, e(d)$ parallel to each blue edge $e$, so that we now have $A_n(i,j) = a$ and $B_n(i,j) = b$, and define the factorisation property by lifting the old factorisation cycle $(e,\Ff(e),\Ff^2(e),\dots \Ff^{a'v}(e)=e)$ to
\[
(e(1),\Ff(e)(1),\dots \Ff^{a'v-1}(e)(1),e(2),\Ff(e)(2),\dots, \Ff^{a'v-1}(e)(d),e(1)).
\]

Next we demonstrate that it suffices to show how to add the desired blue edges when $(v,w) = 1$. To see this, suppose that $v = dv'$ and $w = dw'$, and that we can add the desired edges in the diagram corresponding to $T_n(i,i) = v'$ and $T_{n+1}(j,j) = w'$. Then we may take the resulting diagram, add $d-1$ vertices between pairs of consecutive vertices on $\lambda_{n,j}$ and on $\lambda_{n+1,i}$ and augment each commuting square between vertices on $\lambda_{n,j}$ and $\lambda_{n+1,i}$ to a sequence of $d-1$ commuting squares as shown in Figure~\ref{fig:augment} for $d = 3$.
\begin{figure}[ht]
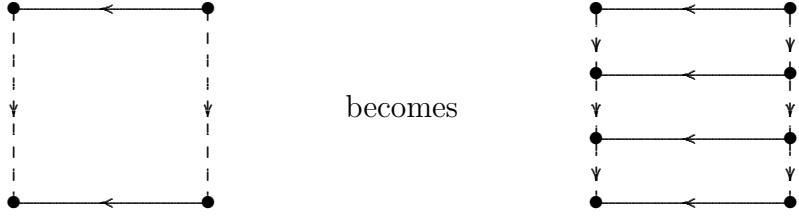

\(
\setcoordinatesystem units <6.25em, 6.25em>
\beginpicture
\put{$\bullet$} at 0 0   \put{$\bullet$} at 1 0   \put{$\bullet$} at 0 1   \put{$\bullet$} at 1 1
\setdashesnear <0.375em> for <6.25em>
\plot 0 0 0 1 /   \plot 1 0 1 1 /
\arr(0,0.451)(0,0.45)  \arr(1,0.451)(1,0.45)
\setsolid
\plot 0 0 1 0 /   \plot 0 1 1 1 /
\arr(0.451,0)(0.45,0)  \arr(0.451,1)(0.45,1)
\put{\text{becomes}} at 2 0.5
\put{$\bullet$} at 3 0       \put{$\bullet$} at 4 0
\put{$\bullet$} at 3 0.333   \put{$\bullet$} at 4 0.333
\put{$\bullet$} at 3 0.666   \put{$\bullet$} at 4 0.666
\put{$\bullet$} at 3 1       \put{$\bullet$} at 4 1
\setdashesnear <0.417em> for <2.083cm>
\plot 3 0 3 0.333 /   \plot 3 0.333 3 0.666 /   \plot 3 0.666 3 1 /
\plot 4 0 4 0.333 /   \plot 4 0.333 4 0.666 /   \plot 4 0.666 4 1 /
\setsolid
\arr(3,0.117)(3,0.116)  \arr(4,0.117)(4,0.116)
\arr(3,0.451)(3,0.45)   \arr(4,0.451)(4,0.45)
\arr(3,0.784)(3,0.783)  \arr(4,0.784)(4,0.783)
\plot 3 0 4 0 /   \plot 3 0.33 4 0.33 /   \plot 3 0.66 4 0.66 /   \plot 3 1 4 1 / 
\arr(3.451,0)(3.45,0)  \arr(3.451,0.33)(3.45,0.33)   \arr(3.451,0.66)(3.45,0.66)   \arr(3.451,1)(3.45,1)  
\endpicture
\)
\caption{Augmenting the commuting squares when $(v,w) \not= 1$.}\label{fig:augment}
\end{figure}

The factorisation property is uniquely determined in each of these augmented paths, and we obtain a diagram with the desired data $T_n(j,j) = v$ and $T_{n+1}(i,i) = w$. Notice that we have multiplied both the order of the factorisation permutation and the number of edges in the picture by the same number $d$.

Finally, we demonstrate how to add the desired blue edges when $(a,b) = 1$ and $(v,w) = 1$. To do this, note that the conditions $(a,b) = 1$ and $(v,w) = 1$ together force $a = w$ and $v = b$. Now adding the complete bipartite graph from the vertices on $\lambda_{n,j}$ to the vertices on $\lambda_{n+1,i}$ gives a 1-skeleton with the desired data and a unique factorisation property. We have $r(\Ff_1(\mu,e)) = r(e)$ if and only if $d(\mu) = k v e_2$ for some $k \in\NN$ and $s(\Ff_1(\mu,e)) = s(e)$ if and only if $d(\mu) = l w e_2$ for some $l$. Since $(v,w) = 1$, it follows that $\Ff_1(\mu,e) \not= e$ for $0 < d(\mu) < vw = va$. 
\end{proof}

\begin{proof}[Proof of Theorem~\ref{thm:classifiable achievability}(2)]
The first claim follows immediately from Proposition~\ref{prp:build the gbd} and Theorem~\ref{thm:AT and K-th}(2). Now suppose that $\varinjlim(\ZZ^{c_n}, A_n)$ is simple and is not isomorphic to $\ZZ$. For $n \ge m$ let $A_{n,m} := A_{n-1} A_{n-2} \dots A_m$. We begin by showing that there is a subsequence $(\ell(n))^\infty_{n=1}$ of $\NN$ for which all entries of $A_{\ell(n+1),\ell(n)}$ are at least $n$.

Since the matrices $A_n$ are proper, the second paragraph of the proof of \cite[Lemma~A4.3]{Eff} shows that we can find a subsequence $k(n)$ of $\NN$ such that all the entries in the matrices $A_{k(n+1),k(n)}$ are positive and nonzero. Since $G \not\cong \ZZ$, \cite[Lemma~A4.4]{Eff} implies that for every $m \in \ZZ^{c_{k(n)}}$, 
\[
\min\{(A_{k(l),k(n)} m)_i : 1 \le i \le c_{k(l)}\}\to \infty\ \mbox{ as $l\to \infty$.}
\] 
It follows that for each $N \in \NN$ there exists $l \ge n$ such that every entry of $A_{k(l),k(n)}$ is greater than $N$. Thus there is a subsequence $\ell(n)$ of $k(n)$ such that every entry of $A_{\ell(n+1),\ell(n)}$ is at least $n$. 

Now $\{\ell(n)\}$ is cofinal in $\NN$ and each $A_{\ell(n+1),\ell(n)}$ is proper by choice. So if we let $c'_n := c_{\ell(n)}$, $A'_n := A_{\ell(n+1), \ell(n)}$, $B'_n := B_{\ell(n+1), \ell(n)}$ and $T'_n := T_{\ell(n)}$ for all $n$, we obtain a commuting diagram of the form~\eqref{eq:comm diagram} in which every entry of $A'_n$ is at least $n$.

Let $\Lambda$ be the \gbd/ obtained by applying Proposition~\ref{prp:build the gbd} to the data $c'_n, A'_n, B'_n, T'_n$. For a blue edge $e$ with range in $V_n$, the order of $e$ is bounded below by the smallest entry of $A'_n$, hence is at least $n$. Thus condition~\eqref{eq:alternate LPF} holds  for $N = 1$. Hence $\Lambda$ has large-permutation factorisations.

Now $P C^*(\Lambda) P$ has the desired $K$-theory by Theorem~\ref{thm:AT and K-th}(2), and it is simple with real-rank zero because it is a full corner in $C^*(\Lambda)$, which is such an algebra by Theorem~\ref{thm:LPF -> RR0}(2) (see \cite[Theorem~3.1.8]{Lin}).
\end{proof}

\begin{example}[The irrational rotation algebras] \label{ex:irration rotation algs}
Fix an irrational number $\theta \in (0,1)$. The irrational rotation algebra $A_\theta$ is the universal $C^*$-algebra generated by unitaries $U, V$ satisfying $UV = e^{2\pi\imath\theta} VU$. Elliott and Evans have proved that that $A_\theta$ is a simple unital A$\TT$ algebra with real-rank zero \cite{EE}, and work of Rieffel and Pimsner-Voiculescu combines to show that $K_0(A_\theta)$ is order-isomorphic to $\ZZ + \theta\ZZ$, and $K_1(A_\theta)$ is isomorphic to $\ZZ^2$ (see \cite{rieff, pv}).

Let $[a_1, a_2, a_3, \dots]$ be the unique simple continued fraction expansion for $\theta$ \cite[Theorem~169]{HW}, and define 
\[
A_n  := \Big(\begin{tabular}{cc} $a_n$ & $1$ \\ $1$ & $0$ \end{tabular}\Big).
\]
Theorem~3.2 of \cite{ES} says that $\ZZ + \theta\ZZ$ is order-isomorphic to $\varinjlim(\ZZ^2, A_n)$. Let $T_n := \id_{c_n}$ and let $B_n = A_n$. Since $\ZZ + \theta \ZZ$ is group-isomorphic to $\ZZ^2$, we obtain a commuting diagram of the form~\eqref{eq:comm diagram} with $\varinjlim(\ZZ^{c_n}, A_n) = \ZZ + \theta \ZZ$ and $\varinjlim(\ZZ^{c_n}, B_n) = \ZZ^2$.

Since $\ZZ + \theta \ZZ$ is a simple dimension group, it follows from Theorem~\ref{thm:classifiable achievability}(2) that there is a \gbd/ $\Lambda_\theta$ such that $C^*(\Lambda_\theta)$ is a simple A$\TT$ algebra with real-rank zero with $K_0(C^*(\Lambda_\theta))$ order-isomorphic to $\ZZ + \theta\ZZ$, and with  $K_1(C^*(\Lambda_\theta))$ isomorphic to $\ZZ^2$. Corollary~\ref{cor:position of unit} implies that $P C^*(\Lambda) P$ has the same $K$-theory with the usual order-unit for $K_0$. Now Elliott's classification theorem for A$\TT$ algebras (as in \cite[Theorem~3.2.6]{Ror}) implies that $P C^*(\Lambda_\theta) P$ is isomorphic to $A_\theta$.

To draw such a \gbd/ $\Lambda_\theta$, take $A_n$ as above. If $\ell(n):=n(n+1)/2$ is the sequence of triangular numbers, then every entry  of $A_{\ell(n+1),\ell(n)}$ is greater than or equal to $n$. Let $\phi_n := A_{\ell(n+1), \ell(n)}$ for all $n$. Then the skeleton of $\Lambda_\theta$ is illustrated by Figure~\ref{fig:Irrational Rotation}, where the label $n$ on a solid edge indicates the presence of $n$ parallel blue edges. The factorisation rules are specified by $\lambda(v)e = \sigma(e)\lambda(v)$ for maximal permutations $\sigma$ of parallel blue edges.
\begin{figure}[ht]
\(
\beginpicture \setcoordinatesystem units <1.65em,1.4em>
\setquadratic
\setdashes <0.417em>
\put{$\bullet$} at 0 0
\put{\plot 0 0 -0.34 .96 0 1.33 /
\plot 0 0 0.34  .96 0 1.33 /
\arrow <0.15cm> [0.25,0.75] from 0.1  1.314 to 0.08 1.328
}[b] at 0 0
\put{$\bullet$} at 0 5
\put{\plot 0 0 -0.34 .96 0 1.33 /
\plot 0 0 0.34  .96 0 1.33 /
\arrow <0.15cm> [0.25,0.75] from 0.1  1.314 to 0.08 1.328
}[b] at  0 5
\put{$\bullet$} at 5 0
\put{\plot 0 0 -0.34 .96 0 1.33 /
\plot 0 0 0.34  .96 0 1.33 /
\arrow <0.15cm> [0.25,0.75] from 0.1  1.314 to 0.08 1.328
}[b] at 5 0
\put{$\bullet$} at 5 5
\put{\plot 0 0 -0.34 .96 0 1.33 /
\plot 0 0 0.34  .96 0 1.33 /
\arrow <0.15cm> [0.25,0.75] from 0.1  1.314 to 0.08 1.328
}[b] at 5 5
\put{$\bullet$} at 10 0
\put{\plot 0 0 -0.34 .96 0 1.33 /
\plot 0 0 0.34  .96 0 1.33 /
\arrow <0.15cm> [0.25,0.75] from 0.1  1.314 to 0.08 1.328
}[b] at 10 0
\put{$\bullet$} at 10 5
\put{\plot 0 0 -0.34 .96 0 1.33 /
\plot 0 0 0.34  .96 0 1.33 /
\arrow <0.15cm> [0.25,0.75] from 0.1  1.314 to 0.08 1.328
}[b] at 10 5
\put{$\bullet$} at 15 0
\put{\plot 0 0 -0.34 .96 0 1.33 /
\plot 0 0 0.34  .96 0 1.33 /
\arrow <0.15cm> [0.25,0.75] from 0.1  1.314 to 0.08 1.328
}[b] at 15 0
\put{$\bullet$} at 15 5
\put{\plot 0 0 -0.34 .96 0 1.33 /
\plot 0 0 0.34  .96 0 1.33 /
\arrow <0.15cm> [0.25,0.75] from 0.1  1.314 to 0.08 1.328
}[b] at 15 5
\put{$\dots$} at 16.5 0
\put{$\dots$} at 16.5 5
\setlinear
\setsolid
\plot 0 0 15 0 /
\arrow <0.5em> [0.25,0.75] from 2.51 0 to 2.49 0
\put{$\scriptstyle\phi_1(2,2)$} at 2.51 -.5
\arrow <0.5em> [0.25,0.75] from 7.51 0 to 7.49 0
\put{$\scriptstyle\phi_2(2,2)$} at 7.51 -.5
\arrow <0.5em> [0.25,0.75] from 12.51 0 to 12.49 0
\put{$\scriptstyle\phi_3(2,2)$} at 12.51 -.5
\plot 0 5 15 5 /
\arrow <0.5em> [0.25,0.75] from 2.51 5 to 2.49 5
\put{$\scriptstyle\phi_1(1,1)$} at 2.51 5.5
\arrow <0.5em> [0.25,0.75] from 7.51 5 to 7.49 5
\put{$\scriptstyle\phi_2(1,1)$} at 7.49 5.5
\arrow <0.5em> [0.25,0.75] from 12.51 5 to 12.49 5
\put{$\scriptstyle\phi_3(1,1)$} at 12.49 5.5
\plot 0 0 5 5 /
\arrow <0.5em> [0.25,0.75] from 1.26 1.26 to 1.24 1.24
\put{$\scriptstyle\phi_1(2,1)$} at 2.2 4
\plot 0 5 5 0 /
\arrow <0.5em> [0.25,0.75] from 1.26 3.74 to 1.24 3.76
\put{$\scriptstyle\phi_1(1,2)$} at 2.2 1
\plot 5 0 10 5 /
\arrow <0.5em> [0.25,0.75] from 6.26 1.26 to 6.24 1.24
\put{$\scriptstyle\phi_2(2,1)$} at 7.2 4
\plot 5 5 10 0 /
\arrow <0.5em> [0.25,0.75] from 6.26 3.74 to 6.24 3.76
\put{$\scriptstyle\phi_2(1,2)$} at 7.2 1
\plot 10 0 15 5 /
\arrow <0.5em> [0.25,0.75] from 11.26 1.26 to 11.24 1.24
\put{$\scriptstyle\phi_3(2,1)$} at 12.2 4
\plot 10 5 15 0 /
\arrow <0.5em> [0.25,0.75] from 11.26 3.74 to 11.24 3.76
\put{$\scriptstyle\phi_3(1,2)$} at 12.2 1
\endpicture
\)
\caption{A \gbd/ for the irrational rotation algebra $A_\theta$.}\label{fig:Irrational Rotation}
\end{figure}
\end{example}

\begin{example}\label{eg:other dim groups}
More generally, let $G$ be a simple dimension group other than $\ZZ$. Write $G = \varinjlim(\ZZ^{c_n}, A_n)$, let $T_n = \id_{c_n}$ and let $B_n = A_n$. As above, we obtain a \gbd/ $\Lambda(G)$ such that $P C^*(\Lambda(G)) P$ is a simple unital $C^*$-algebra with real-rank zero, $K_0(P C^*(\Lambda(G)) P)$ is order-isomorphic to $G$ with the usual order unit and $K_1(P C^*(\Lambda(G)) P)$ is group-isomorphic to $G$. Elliott's classification theorem for A$\TT$ algebras then implies that $P C^*(\Lambda(G)) P$ is the unique A$\TT$ algebra with these properties.
\end{example}

\begin{example}[The Bunce-Deddens algebras]\label{eg:Bunce-Deddens}
As in \cite[Section~7.4]{RLL}, a \emph{supernatural number} is a sequence $\mathbf{m} = (m_i)_{i=1}^\infty$ where each $m_i \in \{0, 1, 2, \dots, \infty\}$. We think of $\mathbf{m}$ as the formal product $\prod^\infty_{i=1} p_i^{m_i}$ where $p_i$ is the $i^{\rm th}$ prime number. We say $\mathbf{m}$ is \emph{infinite} if $\prod^\infty_{i=1} p_i^{m_i} = \infty$, or equivalently if $\sum^\infty_{j=1} m_j = \infty$.

For each supernatural number $\mathbf{m}$, $Q(\mathbf{m})$ denotes the subgroup of $\QQ$ consisting of the fractions of the form $x(\prod_{j = 1}^N p_j^{-q_j})$ with $0 \le q_j \le m_j$ for all $j$. Each $Q(\mathbf{m})$ is a simple dimension group. If $\mathbf{m}$ is finite, then $Q(\mathbf{m}) \cong \ZZ$, so there is no simple A$\TT$ algebra with real-rank zero and $K_0$-group $Q(\mathbf{m})$. If $\mathbf{m}$ is infinite, then $Q(\mathbf{m})_+$ contains no minimal elements, and so $Q(\mathbf{m})$ is not isomorphic to $\ZZ$. Elliott's classification theorem says there is a unique simple unital A$\TT$-algebra $A$ with real-rank zero and $(K_0(A), K_1(A)) = (Q(\mathbf{m}),\ZZ)$. This $C^*$-algebra is known as the Bunce-Deddens algebra of type $\mathbf{m}$; there are several concrete realisations of these algebras, for example as the $C^*$-algebras generated by families of weighted shifts (see \cite{BD} or \cite[\S V.3]{david}), or as crossed products by odometer actions (see \cite[\S VIII.4]{david}). We will demonstrate that for each infinite supernatural number $\mathbf{m}$ there is a \gbd/ $\Lambda(\mathbf{m})$ such that $P C^*(\Lambda(\mathbf{m})) P$ is isomorphic to the Bunce-Deddens algebra of type $\mathbf{m}$.

Fix an infinite supernatural number $\mathbf{m}$. Let $\{a_j\}^\infty_{j=1}$ be any sequence of primes in which each prime $p_i$ occurs with cardinality $m_i$. Then $\varinjlim(\ZZ, \times a_j) \cong Q(\mathbf{m})$ by \cite[Lemma~7.4.4]{RLL}. For $n \in \NN$, let $c_n := 1$, let $A_n := [a_n]$, let $B_n = [1]$ and let $T_n := [a_n]$. This data gives a diagram of the form~\eqref{eq:comm diagram} in which $\varinjlim(\ZZ^{c_n}, A_n) = Q(\mathbf{m})$ and $\varinjlim(\ZZ^{c_n}, B_n) = \ZZ$. It follows from Theorem~\ref{thm:classifiable achievability}(2) that there is a \gbd/ $\Lambda(\mathbf{m})$ such that $C^*(\Lambda(\mathbf{m}))$ is simple and has real-rank zero and $K$-groups $Q(\mathbf{m}), \ZZ$. Corollary~\ref{cor:position of unit} implies that $P C^*(\Lambda) P$ has the same $K$-theory with the usual order-unit for $K_0$. Elliott's classification theorem then implies that $P C^*(\Lambda(\mathbf{m})) P$ is isomorphic to the Bunce-Deddens algebra of type $\mathbf{m}$.

For example, the skeleton of $\Lambda(2^\infty)$ is given in Figure~\ref{fig:Bunce-Deddens}; the factorisation rules are uniquely determined by the skeleton.
\begin{figure}[ht]
\[
\beginpicture
\setcoordinatesystem units <2.5em,1.875em>
\put{$\bullet$} at -2 0
%\put{$\scriptstyle(0,0)$} at -2 -.4
%
\put{$\bullet$} at 0 1
%\put{$\scriptstyle(1,1)$} at 0 1.4
\put{$\bullet$} at 0 -1
%\put{$\scriptstyle(1,0)$} at 0 -1.4
%
\put{$\bullet$} at 4 3
%\put{$\scriptstyle(2,3)$} at 4.45 2.875
\put{$\bullet$} at 4 1
%\put{$\scriptstyle(2,2)$} at 4.45 .95
\put{$\bullet$} at 4 -1
%\put{$\scriptstyle(2,1)$} at 4.45 -.95
\put{$\bullet$} at 4 -3
%\put{$\scriptstyle(2,0)$} at 4.45 -2.875
%
\put{$\bullet$} at 6 4
%\put{$\scriptstyle(3,7)$} at 6.4 4
\put{$\bullet$} at 6 2.857
%\put{$\scriptstyle(3,6)$} at 6.4 2.857
\put{$\bullet$} at 6 1.744
%\put{$\scriptstyle(3,5)$} at 6.4 1.744
\put{$\bullet$} at 6 .602
%\put{$\scriptstyle(3,4)$} at 6.4 .602
\put{$\bullet$} at 6 -.602
%\put{$\scriptstyle(3,3)$} at 6.4 -.602
\put{$\bullet$} at 6 -1.744
%\put{$\scriptstyle(3,2)$} at 6.4 -1.744
\put{$\bullet$} at 6 -2.857
%\put{$\scriptstyle(3,1)$} at 6.4 -2.857
\put{$\bullet$} at 6 -4
%\put{$\scriptstyle(3,0)$} at 6.4 -4
%
\put{$.$} at 6.8 0
\put{$.$} at 7.0 0
\put{$.$} at 7.2 0
\put{$.$} at 6.8 4.4
\put{$.$} at 7.0 4.5
\put{$.$} at 7.2 4.6
\put{$.$} at 6.8 -4.4
\put{$.$} at 7.0 -4.5
\put{$.$} at 7.2 -4.6
\plot -2 0 0 1 /
\arrow <0.15cm> [0.25,0.75] from -0.98  0.51 to  -1.02 0.49
\plot -2 0 0 -1 /
\arrow <0.15cm> [0.25,0.75] from -0.98  -0.51 to  -1.02 -0.49
\plot 0 1 4 3 /
\arrow <0.15cm> [0.25,0.75] from 1.6  1.8 to  1.58 1.79
\plot 0 1 4 -1 /
\arrow <0.15cm> [0.25,0.75] from 1.6  .2 to  1.58 .21
\plot 0 -1 4 1 /
\arrow <0.15cm> [0.25,0.75] from 1.6  -.2 to  1.58 -.21
\plot 0 -1 4 -3 /
\arrow <0.15cm> [0.25,0.75] from 1.6  -1.8 to  1.58 -1.79
\plot 4 3 6 4 /
\arrow <0.15cm> [0.25,0.75] from 5  3.5 to  4.98 3.49
\plot 4 3 6 -.602 /
\arrow <0.15cm> [0.25,0.75] from 5  1.199 to  4.98 1.235
\plot 4 1 6 2.857 /
\arrow <0.15cm> [0.25,0.75] from 5  1.9285 to  4.98 1.9099
\plot 4 1 6 -1.744 /
\arrow <0.15cm> [0.25,0.75] from 5  -.372 to  4.98 -.3445
\plot 4 -1 6 1.744 /
\arrow <0.15cm> [0.25,0.75] from 5  .372 to  4.98 .3445
\plot 4 -1 6 -2.857 /
\arrow <0.15cm> [0.25,0.75] from 5  -1.9285 to  4.98 -1.9099
\plot 4 -3 6 .602 /
\arrow <0.15cm> [0.25,0.75] from 5  -1.199 to  4.98 -1.235
\plot 4 -3 6 -4 /
\arrow <0.15cm> [0.25,0.75] from 5  -3.5 to  4.98 -3.49
\setdashes <0.417em>
\plot 0 -1.05 0 1 /
\arrow <0.15cm> [0.25,0.75] from 0.002  .075 to  0.002 .08
\plot 4 -2.92 4 3 /
\arrow <0.15cm> [0.25,0.75] from 4  1.78 to  4 1.775
\arrow <0.15cm> [0.25,0.75] from 4  -.095 to  4 -.1
\arrow <0.15cm> [0.25,0.75] from 4  -1.99 to  4 -1.995
\plot 6 4.025 6 2.857 /
\arrow <0.15cm> [0.25,0.75] from 6  3.3335 to  6 3.3285
\plot 6 2.89 6 1.744 /
\arrow <0.15cm> [0.25,0.75] from 6  2.2205 to  6 2.2155
\plot 6 1.76 6 .602 /
\arrow <0.15cm> [0.25,0.75] from 6  1.0785 to  6 1.0735
\plot 6 .602 6 -.602 /
\arrow <0.15cm> [0.25,0.75] from 6  -.095 to  6 -.1
\plot 6 -.585 6 -1.744 /
\arrow <0.15cm> [0.25,0.75] from 6  -1.2685 to  6 -1.2745
\plot 6 -1.715 6 -2.857 /
\arrow <0.15cm> [0.25,0.75] from 6  -2.3905 to  6 -2.3955
\plot 6 -2.825 6 -4 /
\arrow <0.15cm> [0.25,0.75] from 6  -3.5035 to  6 -3.5085
\setquadratic
\plot -2 0 -2.226 .72 -2 1 /
\plot -2 0 -1.774 .72 -2 1 /
\arrow <0.15cm> [0.25,0.75] from -1.875  0.971 to -1.89 0.99
%\arrow <0.15cm> [0.25,0.75] from -1.9  1.314 to -1.918 1.322
%
\plot -0.02 .97 -.3 0 0 -1 /
\arrow <0.15cm> [0.25,0.75] from -.3  .075 to  -.3 .08
\plot 3.95 2.85 3.5 0 4 -3 /
\arrow <0.15cm> [0.25,0.75] from 3.498  .075 to  3.498 .08
\plot 5.89 3.72 5.2 0 6 -4 /
\arrow <0.15cm> [0.25,0.75] from 5.1988  .075 to  5.1988 .08
\endpicture
\]
  \caption{A \gbd/ for the Bunce-Deddens algebra of type $2^\infty$.}\label{fig:Bunce-Deddens}
\end{figure}
\end{example}

\section{\Gbds/ with length-1 cycles} \label{sec:trace simplex}

\newcommand{\conv}{\mathrm{conv}}

\newcommand{\tr}{{\rm tr}}
\newcommand{\cM}{{\mathcal{M}}}
\newcommand{\ep}{{\varepsilon}}
\newcommand{\cE}{{\mathcal{E}}}
\newcommand{\Z}{{\mathbb Z}}
\newcommand{\N}{{\mathbb N}}
\newcommand{\T}{{\mathbb T}}
\newcommand{\R}{{\mathbb R}}
\newcommand{\Cs}{{$C^*$-al\-ge\-bra}}
\newcommand{\AT}{{A$\T$}}
\newcommand{\sh}{{$^*$-ho\-mo\-mor\-phism}}

In this section we restrict attention to \gbds/ in which all the isolated cycles have length~1. We show that in this situation, the associated inclusions of circle algebras are \emph{standard permutation mappings}, and that every directed system of direct sums of circle algebras under standard permutation mappings arises from a \gbd/ in which all the isolated cycles have length~1. In the next section we use this to investigate simplicity, real-rank and the trace simplex of the associated $C^*$-algebra in greater detail than the results obtained for general \gbds/ in the previous sections. To state the main theorem of this section, we first need to set up some notation.

Let $\cM(\T)$ denote the set of (positive) probability measures on $\T$. Each $\mu \in \cM(\T)$ induces a functional on $C(\T)$, again 
denoted by $\mu$, given by integration 
\[
\mu(f) = \int_\T f \, d\mu, \qquad f \in  C(\T).
\]
Thus $\cM(\T)$ is a subset of $C(\T)^*$ whence it inherits the weak-$^*$ topology. 

A \emph{Markov operator} on $C(\T)$ is a positive linear mapping $E \colon C(\T) \to C(\T)$ that maps the constant function $1$ to itself. Each Markov operator $E$ induces an affine mapping $\widehat{E} \colon \cM(\T) \to \cM(\T)$ given by $\widehat{E}(\mu)(f) = \mu(E(f))$ which is continuous both wrt.\ the weak-$^*$ and the norm topology on $\cM(\T)$. 

For each $\theta \in \R$ let $\rho_\theta$ be the Markov operator on $C(\T)$ given by rotation by angle $\theta$, that is, $\rho_\theta(f)(z) = f(e^{i\theta}z)$. For each $k \in \N$, let $R_k$ be the Markov operator \label{pg:R_k}
\begin{equation} \label{eq:R_k}
R_k = \frac{1}{k} \sum_{j=0}^{k-1} \rho_{2\pi j/k}.
\end{equation}
Observe that 
\begin{equation} \label{eq:R-comp}
R_k \circ R_\ell = R_{\, \lcm(k,\ell)}.
\end{equation}
It follows from \eqref{eq:R-comp} that if $N$ is a natural number, and if $E_1, E_2 \in \conv\{R_k : k \mid N\}$ and $F_1, F_2 \in \conv\{R_k: k \nmid N\}$, then 
\begin{equation} \label{eq:comp}
E_1E_2 \in \conv\{R_k : k \mid N\}, \qquad F_1E_2,E_1F_2,F_1F_2 \in \conv\{R_k : k \nmid N\}.
\end{equation}
In agreement with the convention mentioned above, $\widehat{\rho}_\theta$ and $\widehat{R}_k$ will be the corresponding affine mappings on $\cM(\T)$.

Given a unital $C^*$-algebra $A$, we write $T(A)$ for the Choquet simplex of tracial states on $A$ endowed with the weak-$^*$ topology.

\begin{dfn} \label{dfn:c}
Let $\sigma$ be a permutation on a finite set $S$. For $s \in S$, $o(s)$ denotes the order $\min\{n > 0 : \sigma^n(s) = s\}$ of $s$ under $\sigma$. We write $\kappa(\sigma)$ for $\max\{o(s) : s \in S\}$. For $\ell \le |S|$, we write $c_\ell(\sigma)$ for the number of orbits under $\sigma$ of size $\ell$. Note that $\sum_\ell \ell c_\ell(\sigma) = |S|$. For $N \in \NN$, we define
\[
\beta_N(\sigma) := \frac{1}{|S|} \sum_{\ell \mid N} \ell c_\ell(\sigma) = \frac{1}{|S|} \big|\{s \in S : o(s) \text{ divides } N\}\big|.
\]
\end{dfn}

The goal of the section is to prove the following theorem. Recall that for a vertex $v$ of a \gbd/, the path $\lambda(v)$ is the isolated cycle in $v (\RED\Lambda) v$. Recall also that $\Ff$ denotes the permutation of $\BLUE\Lambda$ of Section~\ref{sec:LPF}.

\begin{theorem}\label{thm:trace simplex}
Let $\Lambda$ be a \gbd/ of infinite depth in which all red cycles have length $1$. For $n \in \NN$ and $1 \le j \le c_n$, let $v_{n,j}$ be the unique element of $V_{n,j}$. For $n \in \NN$, $1 \le j \le c_n$ and $1 \le i \le c_{n+1}$, let $\Ff_n^{i,j}$ be restriction of $\Ff$ to $v_{n,j} \Lambda^{e_1} v_{n+1, i}$. For $N \in \NN$, let 
\[
\underline{\alpha}_N := \sum_{n \in \NN} (1 - \max_{A_n(i,j) \not= 0}\beta_N(\Ff_n^{i,j}))
\quad\text{and}\quad
\overline{\alpha}_N := \sum_{n \in \NN} (1 - \min_{A_n(i,j) \not= 0}\beta_N(\Ff_n^{i,j})).
\]
Let $P = \sum_{v \in V_0} s_v$. Identify $C^*(\BLUE\Lambda)$ with the subalgebra $C^*(\{s_\xi : \xi \in \BLUE\Lambda\})$ of $C^*(\Lambda)$. Then $C^*(\BLUE\Lambda)$ is AF and is simple if and only if $\Lambda$ is cofinal. Moreover each trace $\tau$ on $P C^*(\Lambda) P$ restricts to a trace on $P C^*(\BLUE\Lambda) P$. The A$\TT$ algebra $C^*(\Lambda)$ is simple if and only if $\Lambda$ is cofinal and 
\[
\sup\{\kappa(\Ff_n^{i,j}) : n \in \NN, 1 \le j \le c_N, 1 \le i \le c_{N+1}, A_n(i,j) \not= 0\} = \infty.
\]

Suppose that $\Lambda$ is cofinal.
\begin{itemize}
\item[(1)] If $\underline{\alpha}_N = \infty$ for all $N$, then $C^*(\Lambda)$ is simple, $P C^*(\Lambda) P$ has real-rank zero, and $\tau \mapsto \tau|_{P C^*(\BLUE\Lambda) P}$ determines an isomorphism between the traces on $P C^*(\Lambda) P$ and the traces on $P C^*(\BLUE\Lambda) P$.
\item[(2)] If $\overline{\alpha}_N < \infty$ for some $N$, then $P C^*(\Lambda) P$ has real-rank one and there is an injective mapping
\begin{equation}\label{eq:trace mapping}
T(P C^*(\BLUE\Lambda) P) \times \{\mu \in \Mm(\TT) : \widehat{R}_N(\mu) = \mu\} \to T(P C^*(\Lambda) P) \qquad (\tau,\mu) \mapsto \bar\tau_\mu
\end{equation}
such that $\bar\tau_\mu|_{P C^*(\BLUE\Lambda) P} = \tau$ for all $\tau \in T(P C^*(\BLUE\Lambda) P)$.
\end{itemize}
Suppose that $|V_n| = c_n = 1$ for all $n$. Then $P C^*(\BLUE\Lambda) P$ has unique trace and $\overline{\alpha}_N = \underline{\alpha}_N$ for all $N$. If $\overline{\alpha}_N = \underline{\alpha}_N < \infty$ for some $N$ then the injection~\eqref{eq:trace mapping} is continuous and affine. If $|V_n| = c_n = 1$ for all $n$ and $\overline{\alpha}_1  = \underline{\alpha}_1 < \infty$ then the injection~\eqref{eq:trace mapping} is a homeomorphism of $\Mm(\TT)$ onto $T(P C^*(\Lambda) P)$.
\end{theorem}

We will prove this theorem on page~\pageref{pg:pf of thm}, after an analysis of the trace simplices of certain A$\TT$ algebras. To apply the results of this analysis, we need to show that the partial inclusions of circle subalgebras of $C^*(\Lambda)$ in the setting of Theorem~\ref{thm:trace simplex} are of a standard form. 

Let $m$ be a natural number, and let $S_m$ denote the group of permutations on $m$ letters. For $\sigma \in S_m$, let $\psi_\sigma \colon C(\T) \to M_m(C(\T))$ \label{pg:psi_sigma} be the $^*$-homomorphism which sends the canonical generator $z$ for $C(\T)$ into the unitary element $\sum_{j=1}^m z e_{j, \sigma(j)}$ in $M_m(C(\T))$, where $\{e_{i,j}\}$ is the set of canonical matrix units for $M_m$. In the special case where $\sigma$ is the $m$-cycle $(1 \; 2 \; 3 \; \cdots \; m)$ the associated \sh{} $\psi_\sigma$ will also be denoted by $\psi_m$, and it is given by 
\begin{equation} \label{eq:psi_k}
\psi_m(z) = \begin{pmatrix}
0 & z & 0 & \cdots & 0 \\ 0 & 0 & z & \cdots & 0 \\
\vdots & \vdots & & \ddots & \vdots \\
0 & 0 & 0 & \cdots & z \\
z & 0 & 0 & \cdots & 0
\end{pmatrix},
\end{equation}
where again $z$ is the canonical generator of $C(\T)$.

If $n,m$ are natural numbers and $\sigma \in S_m$, then we shall also let $\psi_\sigma$ denote the amplified \sh{} $M_n(C(\T)) \to M_{mn}(C(\T))$, or more generally, the not necessarily unital amplified \sh{} $M_n(C(\T)) \to M_k(C(\T))$, where $k$ is any natural number greater than or equal to $mn$, obtained by viewing $M_{mn}(C(\T))$ as a (non-unital) sub-\Cs{} of $M_k(C(\T))$.  We refer to a \sh{} of this form as a \emph{standard permutation mapping}.

\begin{prop}\label{prp:std partial inclusions}
Let $\Lambda$ be a \gbd/ of infinite depth in which each red cycle has length~1. Let $\pi_n : P C^*(\Lambda_n) P \to M_{X_n}(\CC) \otimes C(\TT)$ be the isomorphisms obtained from Proposition~\ref{prp:M_P(C(T))} and Proposition~\ref{prp:oplus}. For each $n \in \NN$, $1 \le j \le c_n$ and $1 \le i \le c_{n+1}$, let $\iota_n^{i,j}$ denote the partial inclusion of the $j^{\rm th}$ summand of $P C^*(\Lambda_N) P$ into the $i^{\rm th}$ summand of $P C^*(\Lambda_{N+1}) P$. Let $\Ff_n^{i,j}$ be as in Theorem~\ref{thm:trace simplex}. Then $\pi_{n+1} \circ \iota_n^{i,j} \circ \pi^{-1}_{n}$ is the standard permutation mapping $\psi_{(\Ff_n^{i,j})^{-1}}$.
\end{prop}

\begin{rmk}
When the red cycles have length~$1$ we need not distinguish a red edge $e_\star$ in each red cycle to obtain the isomorphisms $\pi_n$. Thus the constant matrix units in Lemma~\ref{lem:corner alg} are precisely the $s_\alpha s^*_\beta$ for $\alpha, \beta \in \BLUE\Lambda$. If we demand only that all red cycles have the same length, a result similar to Proposition~\ref{prp:std partial inclusions} holds, but we have to work much harder to show that $\pi_{n+1} \circ \iota_n^{i,j} \circ \pi^{-1}_{n}$ is unitarily equivalent to the standard permutation mapping $\psi_{(\Ff_n^{i,j})^{-1}}$. 

Corollary~\ref{cor:all std perm mappings} shows that any direct system of standard permutation mappings can be realised with the simpler construction where each red cycle has length one, so we omit the more complicated analysis for longer cycles.
\end{rmk}

\begin{proof}[Proof of Proposition~\ref{prp:std partial inclusions}]
For $\alpha, \beta \in X_{n,j}$ and $a,b \in v_{n,j} \Lambda^{e_1} v_{n+1,i}$, let 
\[\textstyle
\theta(\alpha,\beta) := s_\alpha s^*_\beta
\qquad
\theta(\alpha a, \beta b) := s_{\alpha a} s^*_{\beta b}\quad\text{ and }\quad 
\theta(a,b) := \sum_{\eta \in X_{n,j}} s_{\eta a} s^*_{\eta b}. 
\]
Relation~(CK4) shows that for $\alpha,\beta \in X_{n,j}$, the image of $\theta(\alpha,\beta)$ in $P C^*(\Lambda_{n+1, i}) P$ is equal to 
\[
\sum_{a \in v_{n,j} \Lambda^{e_1} v_{n+1, i}} s_{\alpha a} s^*_{\beta a} = \sum_{a \in v_{n,j} \Lambda^{e_1} v_{n+1, i}} \theta(\alpha a, \beta a).
\] 

Using the Cuntz-Krieger relations, we therefore have $\theta(\alpha,\beta) \theta(a,b) = \theta(\alpha a, \beta b) = \theta(a,b) \theta(\alpha, \beta)$ for all $\alpha, \beta \in X_{n,j}$ and $a,b \in v_{n,j} \Lambda^{e_1} v_{n+1, i}$.

Since $P C^*(\Lambda_{n+1, i}) P$ is generated by the matrix units $\theta(\alpha a, \beta b)$, $\alpha a, \beta b \in X_{n+1,i}$ and the unitary $U_{n+1, i} := \sum_{\alpha a \in X_{n+1, i}} s_\alpha s_{\lambda(v_{n+1, i})} s^*_\alpha$, we now have an isomorphism 
\begin{equation}\label{eq:tensor decomp}
M_{X_{n,j}}(\CC) \otimes M_{v_{n,j} \Lambda^{e_1} v_{n+1, i}}(\CC) \otimes C(\TT) \cong P C^*(\Lambda_{n+1, i}) P
\end{equation}
which takes $z \mapsto \Theta(\alpha,\beta) \otimes \Theta(a,b) \otimes z$ to $\theta(\alpha,\beta)\theta(a,b)U_{n+1, i}$. Under this identification, $\pi_{n+1} \circ \iota_n^{i,j} \circ \pi^{-1}_{n}$ takes $\Theta(\alpha,\beta)$ to $\sum_{a \in v_{n, j} \Lambda^{e_1} v_{n+1, i}} \Theta(\alpha,\beta) \otimes \Theta(a,b) \otimes 1$.

Let $U_{n,j} := \sum_{\alpha\in X_{n,j}} s_\alpha s_{\lambda(n,j)} s^*_\alpha$. If we identify $M_{X_{n,j}}(C(\TT))$ with $M_{X_{n,j}}(\CC) \otimes C(\TT)$ in the usual way, then $\pi_n$ takes $\theta(\alpha,\beta) U$ to $\Theta(\alpha,\beta) \otimes z$. The partial inclusion of $P C^*(\Lambda_{n,j}) P$ into $P C^*(\Lambda_{n+1, i}) P$ takes $\theta(\alpha,\beta) U = s_\alpha s_{\lambda(v_{n,j})} s^*_\beta$ to 
\[\textstyle
\sum_{a \in v_{n,j} \Lambda^{e_1} v_{n+1, i}} s_\alpha s_{\lambda(v_{n,j})} s_a s^*_a s^*_\beta = \sum_{a \in v_{n,j} \Lambda^{e_1} v_{n+1, i}} s_{\alpha \Ff_n^{i,j}(a)} s_{\lambda(v_{n+1,i})^m} s^*_{\beta a}
\] 
by~(CK4). Hence under the identification~\eqref{eq:tensor decomp},
\[\textstyle
\pi_{n+1} \circ \iota_n^{i,j} \circ \pi^{-1}_{n}(z) = 
1_{X_{n,j}} \otimes \Big(\sum_{a \in v_{n,j} \Lambda^{e_1} v_{n+1, i}} \Theta(\Ff_n^{i,j}(a), a)\Big) \otimes z.
\]

This completes the proof.
\end{proof}

\begin{cor}\label{cor:all std perm mappings}
Fix integers $c_n, X_{n,j} \in \NN$ for $n \in \NN$ and $1 \le j \le c_n$. Suppose that for each $n$, $\psi_n : \bigoplus^{c_n}_{j=1} M_{X_{n,j}}(C(\TT)) \to \bigoplus^{c_{n+1}}_{i=1} M_{X_{n+1,i}}(C(\TT))$ is a unital inclusion in which all nonzero partial inclusions $\psi_n^{i,j} : M_{X_{n,j}}(C(\TT)) \hookrightarrow M_{X_{n+1, i}}(C(\TT))$ are standard permutation mappings. Then there is a \gbd/ $\Lambda$ in which all red cycles have length~$1$ such that $P C^*(\Lambda) P \cong \varinjlim(\bigoplus^{c_n}_{j=1} M_{X_{n,j}}(C(\TT)), \psi_n)$.
\end{cor}
\begin{proof}
We may assume without loss of generality that for each $n \in \NN$ and each $1 \le j \le c_n$ there exists $1 \le i \le c_{n+1}$ so that $\psi_n^{i,j} \not = 0$. For each $n \in \NN$ and $1 \le i \le c_{n+1}$ there exists $1 \le j \le c_n$ such that $\psi_n^{i,j} \not= 0$ because each $\psi_n$ is unital. 

For each $n,i,j$ such that $\psi_n^{i,j} \not= 0$, let $\sigma^n_{i,j}$ be the permutation such that $\psi_n^{i,j} = \psi_{\sigma^n_{i,j}}$. When $\psi_n^{i,j} \not= 0$, we define $A_n(i,j)$ to be the size of the set of letters acted upon by $\sigma^n_{i,j}$ and regard $\sigma_n^{i,j}$ as a permutation of $\{1, 2, \dots, A_n(i,j)\}$. If $\psi_n^{i,j} = 0$, we define $A_n(i,j) = 0$. The previous paragraph shows that the matrices $A_n$ obtained in this way are all proper. 

We construct $\Lambda$ as follows. Each $V_n$ contains $c_n$ vertices $\{v_{n,1}, \dots, v_{n,c_n}\}$. Each vertex $v_{n,j}$ hosts a single red loop $\lambda_{n,j}$. Insert blue edges $\{e(l) : 1 \le l \le A_n(i,j)\}$ from $v_{n+1, i}$ to $v_{n,j}$ for each $n,i,j$. Specify the commuting squares by $\lambda_{n,j} a(\sigma_n^{i,j}(l)) = a(l) \lambda_{n+1,i}$. This data specifies a unique \gbd/ $\Lambda$ by \cite[page~101]{RSY1}. Proposition~\ref{prp:std partial inclusions} implies that $P C^*(\Lambda) P \cong \varinjlim(\bigoplus^{c_n}_{j=1} M_{X_{n,j}}(C(\TT)), \psi_n)$.
\end{proof}

\section{Real-rank and the trace simplex}\label{sec:mikael}

The results of this section are inspired by Goodearl's paper \cite{Goo}. In this section we continue to use the notation established in Section~\ref{sec:trace simplex}. Let $\sigma$ be the cyclic permutation on $m$ letters. Note that the associated standard permutation mapping $\psi_m = \psi_\sigma$ satisfies
\begin{equation} \label{eq:psi_k-f}
\psi_m(f)(z) \sim_u \begin{pmatrix}
f(z) & 0 &  \cdots & 0 \\ 0 & f(\omega z)  & \cdots & 0 \\
\vdots & \vdots & \ddots &  \vdots \\
0 & 0  & \cdots & f(\omega^{m-1}z)
\end{pmatrix}, \qquad f \in C(\T), \, \, z \in \T,
\end{equation}
where $\omega = \exp(2\pi/m)$. A general permutation $\sigma$ on $m$ letters is the product of disjoint cycles $\sigma_1 \sigma_2 \cdots \sigma_r$ (where we include all 1-cycles). Let $\ell_j$ denote the order of $\sigma_j$ (or, equivalently, the length of the cycle $\sigma_j$). Then $\psi_\sigma$ is unitarily equivalent (by a permutation unitary) to the direct sum $\bigoplus_{j=1}^r \psi_{\ell_j}$. Moreover $\psi_\sigma$ is unitarily equivalent (again with a permutation unitary) to the direct sum $\bigoplus_\ell 1_{c_\ell(\sigma)} \otimes \psi_{\ell}$, where $1_c \otimes \psi$ denotes the $c$-fold direct sum of copies of $\psi$. Notice that $m = \sum_\ell \ell c_\ell(\sigma)$ for all $\sigma \in S_m$.

There is a norm on the linear span of $\cM(\T)$ which on differences of elements from $\cM(\T)$ is the total variation: $\|\mu - \nu\| = |\mu-\nu|(\T)$, and is equal to the operator norm of $\mu-\nu$ when viewed as a functional on $C(\T)$. 

Recall the definitions of $\widehat{R}_k$ and $\widehat{\rho}_k$ from page~\pageref{pg:R_k}.
 
\begin{lemma} \label{lm:m}
\begin{equation} \label{eq:m}
\bigcap_{n=1}^\infty \overline{\conv \{\widehat{R}_k(\mu) : k \ge n, \, \mu \in \cM(\T)\}} = \{m\},
\end{equation}
where the closure is with respect to the norm-topology, and where $m$ denotes the Lebesgue measure (or the normalized Haar measure) on $\T$.
\end{lemma}

\begin{proof} The Lebesgue measure $m$ is the unique rotation invariant measure in $\cM(\T)$, i.e., the only measure that satisfies $\widehat{\rho}_\theta (m) = m$ for all $\theta \in \R$. In particular, $\widehat{R}_k(m) =m$ for all $k$, so $m$ belongs to the left-hand side of \eqref{eq:m}. Suppose, conversely, that $\nu$ is any element belonging to the left-hand side of \eqref{eq:m}. We show that $\widehat{\rho}_\theta (\nu) = \nu$ for all $\theta \in \R$. This will entail that $\nu = m$ and will complete the proof. 

Let $f \in C(\T)$ and let $\ep >0$ be given. Find $\delta > 0$ such that $\|\rho_\theta(f) - f\|_\infty \le \ep$ whenever $|\theta| \le
\delta$.  Note that $\rho_\theta R_k = \rho_{\theta'} R_k$ whenever $k \in \N$ and $\theta, \theta' \in \R$ satisfy $\theta -\theta' \in 2\pi k^{-1} \Z$. For any $k \ge \pi\delta^{-1}$ and for any $\theta \in \R$ we can choose $\theta' \in \R$ such that $|\theta'| \le \delta$ and $\theta -\theta' \in 2\pi k^{-1} \Z$. Then, for any $\mu \in \cM(\T)$, 
\begin{eqnarray*}
|(\widehat{\rho}_\theta \widehat{R}_k \mu)(f) - (\widehat{R}_k \mu)(f)| 
&=& |\mu\big(\rho_\theta R_k(f)-R_k(f)\big)| \; = \;|\mu\big(\rho_{\theta'} R_k(f)-R_k(f)\big)| \\
&=& |\mu\big(R_k(\rho_{\theta'}(f) -f )\big)| \; \le \; \|\rho_{\theta'}(f) - f\| \le \ep.
\end{eqnarray*}
Thus $|(\widehat{\rho}_\theta \mu')(f) - \mu'(f)| \le \ep$ for all 
\[
\mu' \in \overline{\conv \{\widehat{R}_k(\mu) : k \ge \pi\delta^{-1}, \, \mu \in \cM(\T)\}}.
\]
In particular, $|(\widehat{\rho}_\theta \nu)(f) - \nu(f)| \le \ep$. As $f \in C(\T)$ and $\ep>0$ were arbitrary it follows that $\widehat{\rho}_\theta(\nu) = \nu$ for all $\theta \in \R$, as desired.
\end{proof}

Let $\tr_n$ denote the normalized trace on $M_n$, and for $\mu \in \cM(\T)$ and $n\in \N$, let $\tr_{n,\mu}$ denote the normalized trace on $M_n(C(\T))$ given by 
\begin{equation} \label{eq:tr1}
\tr_{n,\mu}(f) = \int_\T \tr_n(f(z)) \, d\mu(z), \qquad f \in M_n(C(\T)).
\end{equation}
Every tracial state on $M_n(C(\T))$ is of the form $\tr_{n,\mu}$ for some $\mu \in \cM(\T)$. 

Consider again the unital \sh{} $\psi_\sigma \colon M_n(C(\T)) \to M_{mn}(C(\T))$ associated to a permutation $\sigma \in S_m$. The induced mapping $T(\psi_\sigma) \colon T(M_{mn}(C(\T))) \to T(M_n(C(\T)))$ is by~\eqref{eq:R_k} and~\eqref{eq:psi_k-f} given as follows: 
\begin{equation} \label{eq:tr2}
\forall \mu_1,\mu_2 \in \cM(\T): \, \,\tr_{nm,\mu_2} \circ \psi_\sigma = \tr_{n,\mu_1} \iff
\mu_1 = \sum_\ell \frac{\ell c_\ell(\sigma)}{m} \widehat{R}_\ell(\mu_2).
\end{equation}

We shall often use the next identity, that holds for any $n \in \N$ and any $\mu,\nu \in \cM(\T)$:
\begin{equation} \label{eq:tr3}
  \|\tr_{n,\mu} - \tr_{n,\nu}\| = \|\mu-\nu\|.
\end{equation}

Finally recall that if $\{s_j\}_{j=1}^\infty$ is a sequence in $(0,1]$, then
\begin{equation} \label{eq:sum-prod}
\prod_{j=1}^\infty s_j > 0 \iff \sum_{j=1}^\infty (1-s_j) < \infty.
\end{equation}

We remind the reader again that for a permutation $\sigma$, the quantities $\kappa(\sigma)$, $c_\ell(\sigma)$ and $\beta_N(\sigma)$ are defined in Definition~\ref{dfn:c}, and that the Markov operators $R_k$ and $\widehat{R}_k$ and the standard permuation mapping $\psi_\sigma$ are defined on pp.~ \pageref{pg:R_k}--\pageref{pg:psi_sigma}. 

\begin{theorem} \label{thm:2} Consider a direct limit of \Cs s
\begin{equation} \label{eq:A2}
\beginpicture
\setcoordinatesystem units <2.5em, 2.5em>
\put{$\bigoplus_{i=1}^{r_1} M_{n_{1,i}}(C(\T))$}[r] at 0 0 
\arrow <0.4em> [0.25, 0.75] from 0.1 0 to 0.9 0
\put{$\scriptstyle \varphi_1$}[b] at 0.5 0.1
\put{$\bigoplus_{i=1}^{r_2}  M_{n_{2,i}}(C(\T))$}[l] at 1 0
\arrow <0.4em> [0.25, 0.75] from 4.1 0 to 4.8 0
\put{$\scriptstyle\varphi_2$}[b] at 4.4 0.1
\put{$\bigoplus_{i=1}^{r_3}  M_{n_{3,i}}(C(\T))$}[l] at 4.9 0 
\arrow <0.4em> [0.25,0.75] from 8 0 to 8.4 0
\put{$\cdots$} at 8.9 0
\arrow <0.4em> [0.25,0.75] from 9.4 0 to 9.9 0 
\put{$A,$}[l] at 10 0 
\endpicture
\end{equation}
with unital connecting maps $\varphi_j$. Let $A_j$ denote the $j^{th}$ algebra in the sequence, so that $A_j = \bigoplus_{i=1}^{r_j} A_{j,i}$, where $A_{j,i} = M_{n_{j,i}}(C(\T))$. Suppose that each of the  partial mappings $\varphi_{j}^{s,t} \colon A_{j,s} \to A_{j+1,t}$ induced by $\varphi_j$ either is zero or is a standard permutation mapping, say of the form $\psi_{\sigma_{j}^{s,t}}$, where $\sigma_{j}^{s,t}$ is a permutation on  $m_j^{s,t}$ letters. Set 
\[
X_j = \{(s,t) : 1 \le t \le r_j, 1 \le s \le r_{j+1}, \varphi_j^{s,t} \ne 0\}\text{ and } X_j(t) = \{s : (s,t) \in X_j\}.
\]

Let $B$ be the AF-algebra associated with the inductive limit in~\eqref{eq:A2}, defined as follows. Let $B_j \subseteq A_j$ be the sub-\Cs{} consisting of all constant functions (so that $B_j = \bigoplus_{i=1}^{r_j} M_{n_{j,i}}$), and observe that $\varphi_j(B_j) \subseteq B_{j+1}$. Set $B = \overline{\bigcup_{j=1}^\infty \varphi_{\infty,j}(B_j)} \subseteq A$, where $\varphi_{\infty,j} \colon A_j \to A$ is the inductive limit map; or equivalently, $B$ is the inductive limit of the sequence $B_1 \to B_2 \to B_3 \to \cdots$. 

Suppose that $B$ is simple. Then:
\begin{itemize}
\item[{\rm(i)}] $A$ is simple if and only if $\sup\{\kappa(\sigma_{j}^{s,t}) : j \in \N, \, (s,t) \in X_j\} = \infty$.
\end{itemize} 
For each natural number $N$, set
\begin{alignat*}{2}
 \overline{\beta}(N,j) &=  \max \{\beta_N(\sigma_{j}^{s,t}) : (s,t) \in X_j\}, &\qquad \underline{\beta}(N,j) 
 &= \min \{\beta_N(\sigma_{j}^{s,t}) : (s,t) \in X_j\},\\
\underline{\alpha}_N  &= \sum_{j=1}^\infty (1-\overline{\beta}(N,j)), &\qquad  \overline{\alpha}_N &= \sum_{j=1}^\infty (1-\underline{\beta}(N,j)).
\end{alignat*}
\begin{itemize}
\item[{\rm{(ii)}}] If $\underline{\alpha}_N=\infty$ for all natural numbers $N$, then $A$ is simple with real-rank zero, and the inclusion mapping $B \to A$ induces an isomorphism $T(A) \to T(B)$ at the level of traces. 
\item[{\rm{(iii)}}] If $\overline{\alpha}_N < \infty$ for some $N$, then $A$ has real-rank one, and there is an injective mapping 
\[
T(B) \times \{\mu \in \cM(\T) : \widehat{R}_N(\mu) = \mu\} \to T(A), \quad (\tau,\mu) \mapsto \overline{\tau}_\mu,
\]
such that each $\overline{\tau}_\mu$ extends $\tau$. 
\item[{\rm{(iv)}}] Suppose that each $r_j = 1$ so each $A_j$ has just a single direct summand. Then $B$ is a UHF algebra and the quantities $\overline{\alpha}_N$ and $\underline{\alpha}_N$ coincide for all $N$. If $\underline{\alpha}_N < \infty$ then the injection of part~(iii) is continuous and affine. If $\underline{\alpha}_1 < \infty$, then the injection of part~(iii) is a homeomorphism of $\Mm(\TT)$ onto $T(A)$.
\end{itemize}
\end{theorem}
\begin{proof}
We can and will assume that the restriction of each $\varphi_j$ to each summand $A_{j,s}$ is non-zero. This will ensure that the connecting maps $\varphi_j$ are injective. 

The connecting mapping $A_j \to A_i$, for $j < i$, is denoted by $\varphi_{i,j}$, and the corresponding partial mapping $A_{j,t} \to A_{i,s}$ is denoted by $\varphi_{i,j}^{s,t}$. As already mentioned, we let $\varphi_{\infty,j}$ denote the inductive limit mapping $A_j \to A$. We identify each $A_{j,s}$ with a sub-\Cs{} of $A_j$, let $\pi_{j,s} \colon A_j \to A_{j,s}$ be the natural conditinal expectation, and denote the unit of $A_{j,s}$ by $e_{j,s}$. Note that $fe_{j,s} = \pi_{j,s}(f)$ for $f \in A_j$. Note also that the projections $\{\varphi_{\infty,j}(e_{j,s})\}$ separate traces on $B$. 

(i). It suffices to show that $\varphi_{\infty,j}(f)$ is full in $A$ for $j \in \NN$ and $f \in A_j \setminus\{0\}$. Given a non-zero $f$ in $A_j$, then $\pi_{j,t_0}(f) \ne 0$ for some $t_0 = 1, \dots, r_j$. Take any non-zero element $b_0$ in $B_j \cap A_{j,t_0}$. Because $B$ is simple and the connecting maps are unital and injective there is $j' > j$ such that $\varphi_{j',j}(b_0)$ is full in $B_{j'}$ (and hence in $A_{j'}$). 

For each $i \ge j$ and for each $s=1, \dots, r_i$ put 
\[
U_{i,s} = \{z \in \T : (\pi_{i,s} \circ \varphi_{i,j})(f)(z) \ne 0\} \subseteq \T.
\]
Suppose that $\varphi_i^{s,t} \ne 0$. Use~\eqref{eq:psi_k-f} and the fact that $\varphi_i^{s,t} = \psi_{\sigma_i^{s,t}}$ is unitarily equivalent to $\bigoplus_\ell 1_{c_\ell(\sigma_i^{s,t})} \otimes \psi_\ell$, to conclude that $U_{i,t} \subseteq U_{i+1,s}$, and that $U_{i+1,s} = \T$ if $U_{i,t}$ contains a closed connected arc of length at least $2\pi \kappa(\sigma_i^{s,t})^{-1}$.   

The set $U := U_{j,t_0}$ is non-empty because $\pi_{j,t_0}(f) \ne 0$. The partial mapping $\varphi_{i,j}^{s,t_0}$ which takes $A_{j,t_0}$ to $A_{i,s}$ is non-zero for all $i \ge j'$ and for all $s=1, \dots, r_i$ because $\varphi_{i,j}(b_0)$ is full in $A_i$ when $i \ge j'$. The argument above therefore shows that $U \subseteq U_{i,s}$ for all $i \ge j'$ and for all $s$. The assumption that $\{\kappa(\sigma_{i}^{s,t})\}$ is unbounded implies that there is $i \ge j'$ and $(s,t) \in X_i$ such that $U$ contains a closed connected arc of length at least $2 \pi \kappa(\sigma_{i}^{s,t})^{-1}$. Thus $U_{i+1,s} = \T$, or, in other words, $(\pi_{i+1,s} \circ \varphi_{i+1,j})(f)$ is full in $A_{i+1,s}$. This shows that the ideal in $A$ generated by $\varphi_{\infty,j}(f)$ contains $\varphi_{\infty,i+1}(A_{i+1,s})$, and hence has non-zero intersection with $B$, so is equal to $A$. 

Suppose now that $\{\kappa(\sigma_{i}^{s,t})\}$ is bounded. Then there is a natural number $N$ such that $\ell \mid N$ for all $\ell$ for which $c_\ell(\sigma_j^{s,t}) \ne 0$ for some $j$ and some $(s,t) \in X_j$. Let $g_{N,j} \in A_j = C(\T,B_j)$ be given by $g_N(z) = z^N 1_{B_j}$. It follows from~\eqref{eq:psi_k} that $\varphi_j(g_{N,j}) = g_{N,j+1}$ for all $j$. As $g_{N,j}$ is central in $A_j$ for all $j$, $\varphi_{\infty,1}(g_{N,1})$ belongs to the centre of $A$. Hence $A$ has non-trivial centre, so $A$ is non-simple. 

(ii) and (iii). Each tracial state $\tau_j$ on $A_j$ is of the form 
\begin{equation} \label{eq:tau_j}
\tau_j(f) = \sum_{t=1}^{r_j} a_{j,t} \tr_{n_{j,t},\mu_{j,t}}(\pi_{j,t}(f)),
\end{equation} 
for some $\mu_{j,1}, \dots, \mu_{j,r_j} \in \cM(\T)$ where each $a_{j,t} \ge 0$ is the value $\tau_j(e_{j,t})$ of $\tau$ at the unit $e_{j,t}$ for $A_{j,t}$. We show first that if $\tau_j$ and $\tau_{j+1}$ are traces on $A_j$ and $A_{j+1}$, respectively, given as in~\eqref{eq:tau_j}, if $\tau_{j+1} \circ \varphi_j = \tau_j$, and if $a_{j,t} \ne 0$, then
\begin{eqnarray} 
\mu_{j,t} 
&=& \sum_{s \in X_j(t)} \frac{a_{j+1,s}}{a_{j,t}} \frac{m_j^{s,t} n_{j,t}}{n_{j+1,s}} \sum_\ell \frac{\ell c_\ell(\sigma_j^{s,t})}{m_j^{s,t}} \, \widehat{R}_\ell(\mu_{j+1,s}), \label{eq:a} \\
1 &= &\sum_{s\in X_j(t)} \frac{a_{j+1,s}}{a_{j,t}} \frac{m_j^{s,t} n_{j,t}}{n_{j+1,s}} \, =  \, \sum_\ell \frac{\ell c_\ell(\sigma_j^{s,t})}{m_j^{s,t}}. \label{eq:b}
\end{eqnarray} 
The second identity in \eqref{eq:b} follows by the definition of the coefficients $c_\ell$. The first identity in \eqref{eq:b} follows from the 
calculation:
\begin{eqnarray*}
a_{j,t} 
&=& \tau_j(e_{j,t}) = \tau_{j+1}(\varphi_j(e_{j,t})) \\
&=& \sum_{s \in X_j(t)} a_{j+1,s} \tr_{n_{{j+1},s},\mu_{j+1,s}}(\varphi_j^{s,t}(e_{j,t})) \\
&=& \sum_{s\in X_j(t)} a_{j+1,s} \frac{m_j^{s,t} n_{j,t}}{n_{j+1,s}},
\end{eqnarray*}
where we have used that $\dim(e_{j,t}) = n_{j,t}$ and that the multiplicity of $\varphi_j^{s,t}$ is $m_j^{s,t}$. We proceed to prove \eqref{eq:a}. Two applications of~\eqref{eq:tau_j} yield
\begin{equation}\label{eq:a_j,t relation}
a_{j,t} \tr_{n_{j,t},\mu_{j,t}} = \tau_j \circ \pi_{j,t} = \tau_{j+1} \circ \varphi_j \circ \pi_{j,t} = \sum_{s\in X_j(t)} a_{j+1,s} \tr_{n_{{j+1},s},\mu_{{j+1,s}}} \circ \psi_{\sigma_j^{s,t}}.
\end{equation}
We wish to apply~\eqref{eq:tr2} to right-hand side of~\eqref{eq:a_j,t relation}, but we must take into account that the \sh{} $\varphi_j^{s,t} \colon A_{j,t} \to A_{j+1,s}$ is not unital. This is done by adjusting the right-hand side of~\eqref{eq:tr2} by the factor $\dim(\varphi_j^{s,t}(e_{j,t})) / n_{j+1,s} = n_{j,t}m_j^{s,t}/n_{j+1,s}$. Now, \eqref{eq:a} follows from~\eqref{eq:a_j,t relation} and from the modified~\eqref{eq:tr2}.

For any natural number $N$ we rewrite \eqref{eq:a} as 
\begin{align*}
\mu_{j,t} 
&= \sum_{s \in X_j(t)} \frac{a_{j+1,s}}{a_{j,t}} \frac{m_j^{s,t} n_{j,t}}{n_{j+1,s}} \Big(\sum_{\ell \, \mid N} \frac{\ell c_\ell(\sigma_j^{s,t})}{m_j^{s,t}} \widehat{R}_\ell(\mu_{j+1,s}) + \sum_{\ell \, \nmid N} \frac{\ell c_\ell(\sigma_j^{s,t})}{m_j^{s,t}} \widehat{R}_\ell(\mu_{j+1,s})\Big) \\
&= \sum_{s \in X_j(t)} \gamma_{N,j}^{s,t} \, \widehat{E}_{N,j}^{s,t}(\mu_{j+1,s}) + \sum_{s \in X_j(t)} \eta_{N,j}^{s,t} \, \widehat{F}_{N,j}^{s,t}(\mu_{j+1,s}),
\end{align*}
where $E_{N,j}^{s,t} \in \conv\{R_k : k \mid N\}$, $F_{N,j}^{s,t} \in \conv\{R_k : k \nmid N\}$, and for $(s,t) \in X_j$, 
\[
\gamma_{N,j}^{s,t} = \frac{a_{j+1,s}}{a_{j,t}} \frac{m_j^{s,t} n_{j,t}}{n_{j+1,s}} \beta_N(\sigma_j^{s,t}), \qquad
\eta_{N,j}^{s,t} = \frac{a_{j+1,s}}{a_{j,t}} \frac{m_j^{s,t} n_{j,t}}{n_{j+1,s}} (1-\beta_N(\sigma_j^{s,t})).
\]
Put $\gamma_{N,j}^{s,t} = \eta_{N,j}^{s,t} = 0$ when $(s,t) \notin X_j$. Note that $\sum_{s=1}^{r_{j+1}} (\gamma_{j,N}^{s,t} + \eta_{j,N}^{s,t}) = 1$, and that 
\[
\underline{\beta}(N,j) \le \sum_{s=1}^{r_{j+1}} \gamma_{j,N}^{s,t} \le \overline{\beta}(N,j).
\]
Suppose now that we have a tracial state $\tau_j$ on $A_j$ for all $j$ (given as in \eqref{eq:tau_j} above) such that $\tau_{j+1} \circ \varphi_j = \tau_j$ holds for all $j$. It then follows from iterated use of the identities established above, together with \eqref{eq:comp}, that for $i >j$, $t=1, \dots,r_j$, and $s = 1, \dots, r_i$, 
\begin{equation} \label{eq:i}
\mu_{j,t} = \sum_{s=1}^{r_i} \gamma_{N,i,j}^{s,t} \, \widehat{E}_{N,i,j}^{s,t}(\mu_{i,s}) +  \sum_{s=1}^{r_i} \eta_{N,i,j}^{s,t} \, \widehat{F}_{N,i,j}^{s,t}(\mu_{i,s}),
\end{equation}
where $E_{N,i,j}^{s,t} \in \conv\{R_k : k \mid N\}$, $F_{N,i,j}^{s,t} \in \conv\{R_k : k \nmid N\}$, and where $\gamma_{N,i,j}^{s,t}$ and $\eta_{N,i,j}^{s,t}$ are non-negative real numbers satisfying
\begin{equation} \label{eq:ii}
\sum_{s=1}^{r_i}(\gamma_{N,i,j}^{s,t} + \eta_{N,i,j}^{s,t}) = 1, \qquad \sum_{s=1}^{r_i}\gamma_{N,i,j}^{s,t} \le \prod_{k=j}^{i-1} \overline{\beta}(N,k).
\end{equation}

(ii). We first show that each tracial state $\tau$ on $B$ lifts to a tracial state $\overline{\tau}$ on $A$. Indeed, for each $j \in \N$, let $\overline{\tau}_j$ be the trace on $A_j$ given as in~\eqref{eq:tau_j} with 
\[
a_{j,t} = \tau(\varphi_{\infty,j}(e_{j,t})), \qquad \mu_{j,t} = m, \quad t=1, \dots, r_j,
\]
(where $m$ is the Lebesgue measure). Since $\widehat{R}_\ell(m) = m$ for all $\ell$ it follows from~\eqref{eq:a} that $\overline{\tau}_{j+1} \circ \varphi_j = \overline{\tau}_j$ for all $j$, and so there is a trace $\overline{\tau}$ on $A$, which satisfies $\overline{\tau} \circ \varphi_{\infty,j} = \overline{\tau}_j$ for all $j$. (The first equation in \eqref{eq:b} holds because $\tau$ is a trace on $B$.) In particular,
\[
\tau(\varphi_{\infty,j}(e_{j,t})) = a_{j,t} = \overline{\tau}_j(e_{j,t}) = \overline{\tau}(\varphi_{\infty,j}(e_{j,t})).
\]
Since $\{\varphi_{\infty,j}(e_{j,t})\}$ separate traces on $B$ we conclude that $\overline{\tau}|_B = \tau$.  

We now show that the lift constructed above is unique. Here we need our assumption that $\underline{\alpha}_N=\infty$ for all $N$. Let again $\tau$ be a tracial state on $B$ and suppose that $\widetilde{\tau}$ is (another) tracial state on $A$ that extends $\tau$. Then 
\[
a_{j,t} := \widetilde{\tau}(\varphi_{\infty,j}(e_{j,t})) = \tau(\varphi_{\infty,j}(e_{j,t})).
\]
Now, $\widetilde{\tau}_j := \widetilde{\tau} \circ \varphi_{\infty,j}$ is a trace on $A_j$ which therefore is given as in \eqref{eq:tau_j} with $a_{j,t}$ as above and with respect to some measures $\mu_{j,1}, \dots, \mu_{j,r_j} \in \cM(\T)$. We must show that $\mu_{j,t} = m$ for all $j$ and $t$. (This will show that $\widetilde{\tau} = \overline{\tau}$, cf.\ the construction of $\overline{\tau}$ above.)

The assumption that $\underline{\alpha}_N=\infty$ implies that 
\[
\lim_{i \to \infty} \prod_{k=j}^{i-1} \overline{\beta}(N,j) = 0
\]
for all $j,N \in \N$, cf.\ \eqref{eq:sum-prod}. It follows from equations \eqref{eq:i} and \eqref{eq:ii} that $\mu_{j,t}$ belongs to the norm closure of $\conv\{R_k : k \nmid N\}$ for all $N$, and hence, upon choosing $N=(n-1)!$, that $\mu_{j,t}$ belongs to the norm closure of $\conv\{R_k : k \ge n\}$ for all $n$. By Lemma~\ref{lm:m} this implies that $\mu_{j,t} = m$, as desired.

We use (i) to show that $A$ is simple. Let $n$ be a natural number and put $N=(n-1)!$. Since $\underline{\alpha}_N = \infty$, there exist $j,s,t$
such that $\beta_N(\sigma_j^{s,t}) < 1$. This implies that $\kappa(\sigma_j^{s,t}) \ge n$. Hence $\{\kappa(\sigma_j^{s,t})\}$ is unbounded.

Projections in $B$ separate traces on $B$ because $B$ is of real rank zero, being an AF-algebra. We have shown that each trace in $B$ has a unique lift to a trace on $A$.  It follows that projections in $B$ (and hence also projections in $A$) separate traces on $A$. We can therefore use \cite[Theorem~1.3]{BBEK} to conclude that $A$ has real-rank zero. 

(iii). Assume that $\overline{\alpha}_N < \infty$ for some $N$. We construct an injective mapping 
\[
T(B) \times \{\mu \in \cM(\T) : \widehat{R}_N(\mu) = \mu\} \to T(A), \qquad (\tau,\mu) \mapsto \overline{\tau}_\mu,
\]
such that each $\overline{\tau}_\mu$ extends $\tau$. 

Let $\tau \in T(B)$ and let $\mu \in \cM(\T)$ with $\widehat{R}_N(\mu)=\mu$ be given. We proceed to construct the tracial state $\overline{\tau}_\mu$ on $A$ that extends $\tau$. Let $\widetilde{\tau}_{\mu,j}$ be the trace on $A_j$ given as in~\eqref{eq:tau_j} with 
\[
a_{j,t} = \tau(\varphi_{\infty,j}(e_{j,t})), \qquad \mu_{j,t} = \mu, \quad t=1, \dots, r_j.
\]
Use~\eqref{eq:a}, \eqref{eq:b}~and~\eqref{eq:tr3} to see that 
{\allowdisplaybreaks
\begin{align}
\|\widetilde{\tau}_{\mu,i+1} \circ \varphi_i - \widetilde{\tau}_{\mu,i}\| 
&\le  \sum_{t=1}^{r_i} a_{i,t} \Big\| \sum_{s \in X_{i}(t)} \frac{a_{i+1,s}}{a_{i,t}} \frac{m_i^{s,t} n_{i,t}}{n_{i+1,s}} \sum_\ell \frac{\ell
c_\ell(\sigma_i^{s,t})}{m_i^{s,t}} \, \big(\widehat{R}_\ell(\mu) - \mu\big)\Big\| \nonumber\\
&= \sum_{t=1}^{r_i} a_{i,t} \Big\| \sum_{s \in X_i(t)} \frac{a_{i+1,s}}{a_{i,t}} \frac{m_i^{s,t} n_{i,t}}{n_{i+1,s}} \sum_{\ell \, \nmid N} \frac{\ell c_\ell(\sigma_i^{s,t})}{m_i^{s,t}} \, \big(\widehat{R}_\ell(\mu) - \mu\big)\Big\|  \nonumber\\ 
&\le \sum_{t=1}^{r_i} a_{i,t} \sum_{s \in X_i(t)} \frac{a_{i+1,s}}{a_{i,t}} \frac{m_i^{s,t} n_{i,t}}{n_{i+1,s}} \sum_{\ell \, \nmid N} \frac{\ell
c_\ell(\sigma_i^{s,t})}{m_i^{s,t}}  \nonumber\\
&\le \sum_{t=1}^{r_i} a_{i,t} \sum_{s \in X_i(t)} \frac{a_{i+1,s}}{a_{i,t}} \frac{m_i^{s,t} n_{i,t}}{n_{i+1,s}} ( 1- \underline{\beta}(N,i))  \nonumber \\
&= 1- \underline{\beta}(N,i). \label{eq:trace approximation}
\end{align}}
The hypothesis $\overline{\alpha}_N = \sum_{k=1}^\infty (1-\underline{\beta}(N,k)) < \infty$ implies that $\delta_i := \sum_{k=i}^\infty (1 - \underline{\beta}(N,k)) \to 0$ as $i \to \infty$. We deduce that $\{\widetilde{\tau}_{\mu,i} \circ \varphi_{i,j}\}_{i=j}^\infty$ is a Cauchy sequence in norm, and thus converges to a trace $\overline{\tau}_{\mu,j}$ on $A_j$ which satisfies $\|\overline{\tau}_{\mu,j} - \widetilde{\tau}_{\mu,j}\| \le \delta_j$. As $\overline{\tau}_{\mu,j+1} \circ \varphi_j = \overline{\tau}_{\mu,j}$ for all $j$, there is a tracial state $\overline{\tau}_{\mu}$ on $A$ such that $\overline{\tau}_\mu \circ \varphi_{\infty,j} = \overline{\tau}_{\mu,j}$ for all $j$. 

To show that $\overline{\tau}_\mu$ extends $\tau$ observe first that $\widetilde{\tau}_{\mu,j}(e_{j,t}) = a_{j,t} = (\tau\circ\varphi_{\infty,j})(e_{j,t})$ for all $j$ and $t$, and hence that the restriction of $\widetilde{\tau}_{\mu,j}$ to $B_j$ is equal to $\tau \circ
\varphi_{\infty,j}$. As $\widetilde{\tau}_{\mu,i} \circ \varphi_{i,j} \to \overline{\tau}_\mu \circ \varphi_{\infty,j}$, the restriction of $\overline{\tau}_{\mu} \circ \varphi_{\infty,j}$ to $B_j$ is equal to $\tau \circ \varphi_{\infty,j}$ for all $j$. Hence $\overline{\tau}_\mu$ extends $\tau$. 

Assume that $\mu, \nu \in \cM(\T)$ are such that $\widehat{R}_N(\mu)=\mu$, $\widehat{R}_N(\nu)=\nu$, and $\overline{\tau}_\mu = \overline{\tau}_\nu$. Then, using (a slightly modified version of)~\eqref{eq:tr3}, 
\begin{align*}
\|\mu-\nu\| &= \|\widetilde{\tau}_{\mu,j}-\widetilde{\tau}_{\nu,j}\| 
\\
&\le \|\widetilde{\tau}_{\mu,j}- \overline{\tau}_\mu \circ \varphi_{\infty,j}\| + \| \overline{\tau}_\mu \circ \varphi_{\infty,j} - \overline{\tau}_\nu \circ \varphi_{\infty,j}\| + \|\overline{\tau}_\nu \circ \varphi_{\infty,j} -\widetilde{\tau}_{\nu,j}\| \\ 
&\le 2 \delta_j,
\end{align*}
for all $j$, which entails that $\mu = \nu$. 

We claim that $\overline{\tau}_\mu(p) = \overline{\tau}_\nu(p)$ for every projection $p \in A$. To see this, note that each projection $p \in A$ is equivalent to $\varphi_{\infty,j}(q)$ for some projection $q$ in some $A_j$. Now, each projection $q$ in $A_j$ is equivalent to a projection $q'$ in $B_j$, so $p$ is equivalent to the projection $p'=\varphi_{\infty,j}(q')$ in $B$. This proves that $\overline{\tau}_\mu(p) = \overline{\tau}_\mu(p') = \tau(p') = \overline{\tau}_\nu(p')= \overline{\tau}_\nu(p)$, establishing the claim. 

The stable rank of any A$\T$-algebra is one, and hence its real rank must be either zero or one. The claim above and that $\overline{\tau}_\mu \ne 
\overline{\tau}_\nu$ whenever $\mu$ and $\nu$ are distinct measures fixed under $\widehat{R}_N$ show that projections in $A$ do not separate traces on $A$. Hence $A$ cannot be of real rank zero, and must therefore be of real rank one.

(iv). Suppose now that $r_n = 1$ for all $n$. Then $B$ is a direct limit of unital inclusions of simple finite-dimensional $C^*$-algebras, and hence is UHF, and in particular is simple and has unique trace $\tau^B$. It is immediate from the definitions of $\underline{\alpha}_N$ and $\overline{\alpha}_N$ that these quantities coincide for all $N$. If $\underline{\alpha}_N = \overline{\alpha}_N < \infty$, then injection of statement~(iii) depends on only one variable as $T(B) = \{\tau^B\}$. Hence we write $\tau^B_\mu$ rather than $\overline{\tau^B}_\mu$ for the trace on $A$ corresponding to a given $\mu \in \Mm(\TT)$. Since each $B_j$ has just one summand $B_{j,1}$ we will drop the second subscript henceforth, and write $B_j$ for $B_{j,1}$, $n_j$ for $n_{j,1}$, etc. The $j^{\rm th}$ approximating algebra $B_j$ has unique trace $\tr_{n_j}$, so we can use~\eqref{eq:trace approximation} and the subsequent paragraph to deduce that
\begin{equation}\label{eq:tau_mu}
\|\tau^B_\mu \circ \varphi_{\infty,j} - \tr_{n_j,\mu}\| \le \delta_j \to 0.
\end{equation}
Since $\mu \mapsto \tr_{n_j, \mu}$ is affine, we conclude that $\mu \mapsto \tau^B_\mu$ is affine. 

To see that the map $\mu \mapsto \tau^B_\mu$ is continuous, take a net $\{\mu_\alpha\}$ in $\cM(\T)$ which converges in the weak-$^*$ topology to $\mu \in \cM(\T)$. To show that $\Tt_{\mu_\alpha} \to \tau^B_\mu$ it suffices to show that $\Tt_{\mu_\alpha}(a) \to \tau^B_\mu(a)$ for all $a$ in the dense subset $\bigcup_{j=1}^\infty \varphi_{\infty,j}(A_j)$ of $A$. In other words, it suffices to show that $\Tt_{\mu_\alpha} \circ \varphi_{\infty, j}(f) \to \Tt_{\mu} \circ \varphi_{\infty, j}(f)$ for all $f \in A_j$. But if $g : A_j \to \CC$ is the function $g(z) = \tr_{n_j}(f(z))$, then
\[
\Tt_{\mu_\alpha}\circ \varphi_{\infty, j}(f) = \mu_\alpha(g) \to \mu(g) = \Tt_{\mu}\circ \varphi_{\infty, j}(f).
\]

Finally, suppose that $\underline{\alpha}_1 < \infty$. We show that $\mu \mapsto \tau^B_\mu$ is surjective, and being a continuous bijection between compact sets, it must then be a homeomorphism. 

Fix $\tau \in T(A)$. Then $\tau \circ \varphi_{\infty,j}$ is a trace on $A_j$, and is hence equal to $\tr_{n_j,\mu_j}$ for some $\mu_j \in \cM(\T)$. Since $\tr_{n_{j+1},\mu_{j+1}} \circ \varphi_j = \tr_{n_j,\mu_j}$ we can use~\eqref{eq:tr2} to estimate  
\begin{eqnarray*}
\|\mu_j - \mu_{j+1}\| 
&=&\Big\| \sum_\ell \frac{\ell c_\ell(\sigma_j)}{m_j} \widehat{R}_\ell(\mu_{j+1})-\mu_{j+1}\Big\| \\
&=& \Big\| \sum_{\ell \, >1} \frac{\ell c_\ell(\sigma_j)}{m_j} \big(\widehat{R}_\ell(\mu_{j+1})-\mu_{j+1}\big)\Big\| \\
&\le& \sum_{\ell \, >1} \frac{\ell c_\ell(\sigma_j)}{m_j} \\
&=& 1-\beta(1,j).
\end{eqnarray*}
We have assumed that $\underline{\alpha}_1 = \sum_{j=1}^\infty (1-\beta(1,j)) < \infty$. Hence $\{\mu_j\}$ is norm convergent to a measure $\mu \in \cM(\T)$. Moreover,
\begin{eqnarray*}
\|\tau^B_\mu \circ \varphi_{\infty,j} - \tau \circ \varphi_{\infty,j}\| 
& = & \|\tau^B_\mu \circ \varphi_{\infty,j} - \tr_{n_j,\mu_j}\| \\ 
&\le& \|\tau^B_\mu \circ \varphi_{\infty,j} - \tr_{n_j,\mu}\| + \|\tr_{n_j,\mu}-\tr_{n_j,\mu_j}\| \\
&=& \delta_j + \|\mu-\mu_j\| \to 0, \qquad \text{by \eqref{eq:tau_mu}}
\end{eqnarray*}
which proves that $\tau = \tau^B_\mu$.
\end{proof} 

\begin{proof}[Proof of Theorem~\ref{thm:trace simplex}]\label{pg:pf of thm}
The $1$-graph of Proposition~\ref{prp:associated BD} is $\BLUE\Lambda$ and $Q$ is equal to $P$. Hence $C^*(\BLUE\Lambda)$ is AF, and is simple if and only if $\BLUE\Lambda$ (equivalently $\Lambda$) is cofinal \cite{KPR}. The simplicity statement for $C^*(\Lambda)$ follows from Theorem~\ref{critsimple}(1). Proposition~\ref{prp:std partial inclusions} shows that the partial inclusions in the direct limit decomposition of $C^*(\Lambda)$ are standard inclusions with permutations $(\Ff_n^{i,j})^{-1}$. Moreover, the approximating subalgebras $F_n$ in $C^*(\BLUE\Lambda)$ from Proposition~\ref{prp:associated BD} are the subalgebras of constant functions in the approximating subalgebras of $C^*(\Lambda)$, so $C^*(\BLUE\Lambda)$ is the AF algebra $B$ associated to $C^*(\Lambda)$ in Theorem~\ref{thm:2}. Since each $\kappa((\Ff_n^{i,j})^{-1}) = \kappa(\Ff_n^{i,j})$ and each $c_\ell((\Ff_n^{i,j})^{-1}) = c_\ell(\Ff_n^{i,j})$, the remaining statements of the theorem now follow from Theorem~\ref{thm:2}. 
\end{proof}

\end{document}